\documentclass[reqno,10pt]{amsart}
\usepackage[foot]{amsaddr}
\usepackage{amssymb,amsmath,amsthm,amsfonts}
\usepackage{mathrsfs,dsfont,comment,mathscinet,mathtools}
\usepackage{amsmath,amssymb,amsthm,mathtools,regexpatch}
\usepackage[usenames,dvipsnames]{xcolor} 	
\usepackage{enumitem} 
\usepackage{mathrsfs,stmaryrd}
\usepackage[T1]{fontenc}
\usepackage{hyperref}

\usepackage{fancyhdr}
\usepackage[left=2.3cm,right=2.3cm,top=2.3cm,bottom=2.3cm]{geometry}
\numberwithin{equation}{section}

\theoremstyle{plain}
\newtheorem{theorem}{Theorem}[section]
\newtheorem{proposition}[theorem]{Proposition}
\newtheorem{corollary}[theorem]{Corollary}
\newtheorem{lemma}[theorem]{Lemma}
\newtheorem{example}[theorem]{Example}

\theoremstyle{remark}

\theoremstyle{definition}
\newtheorem{definition}[theorem]{Definition}

\makeatletter
\def\th@plain{%
	\thm@notefont{}
	\itshape 
}
\def\th@definition{%
	\thm@notefont{}
	\normalfont 
}
\makeatother

\makeatletter
\xpatchcmd{\proof}{\itshape}{\bfseries}{}{}
\makeatother

\usepackage{accents}
\renewcommand*{\dot}[1]{%
	\accentset{\mbox{\large\bfseries .}}{#1}}

\newcommand{\osc}{\operatorname*{osc}}
\DeclareMathOperator{\adj}{adj} 
\newcommand{\aplim}{\operatorname*{ap\:lim}}
\DeclareMathOperator{\cof}{cof}

\DeclareMathOperator{\dist}{dist}       
\renewcommand{\div}{\operatorname{div}}
\DeclareMathOperator{\loc}{loc}

\DeclareMathOperator{\sgn}{sgn}

\DeclareMathOperator*{\esssup}{ess \: sup}
\DeclareMathOperator*{\essinf}{ess \: inf}

\newcommand{\supp}{\mathrm{supp}}

\newcommand{\Lip}{\mathrm{Lip}}

\renewcommand{\d}{\mathrm{d}}
\newcommand{\N}{\mathbb{N}}
\newcommand{\Z}{\mathbb{Z}}
\newcommand{\R}{\mathbb{R}} 
\renewcommand{\S}{\mathbb{S}}      
\newcommand{\restr}[1]{|_{#1}}
\newcommand{\imt}{\mathrm{im}_{\rm T}}
\newcommand{\img}{\mathrm{im}_{\rm G}}
\newcommand{\domg}{\mathrm{dom}_{\rm G}}
\newcommand\wk{\rightharpoonup}
\newcommand{\wks}{\overset{\ast}{\rightharpoonup}}

\newcommand{\leb}{\mathscr{L}^N}
\newcommand{\haus}{\mathscr{H}^{N-1}}
\newcommand{\RNN}{\mathbb{R}^{N \times N}}
\newcommand{\RN}{\mathbb{R}^N}
\newcommand{\SN}{\mathbb{S}^{N-1}}

\newcommand*\closure[1]{\overline{#1}}

\def\Xint#1{\mathchoice
{\XXint\displaystyle\textstyle{#1}}%
{\XXint\textstyle\scriptstyle{#1}}%
{\XXint\scriptstyle\scriptscriptstyle{#1}}%
{\XXint\scriptscriptstyle\scriptscriptstyle{#1}}%
\!\int}
\def\XXint#1#2#3{{\setbox0=\hbox{$#1{#2#3}{\int}$}
\vcenter{\hbox{$#2#3$}}\kern-.5\wd0}}

\def\dashint{\Xint-}

\newcommand{\mres}{\mathbin{\vrule height 1.6ex depth 0pt width
		0.13ex\vrule height 0.13ex depth 0pt width 1.3ex}}

\newcommand{\EEE}{\color{black}}

\usepackage{stmaryrd}
\SetSymbolFont{stmry}{bold}{U}{stmry}{m}{n}

\title[Quasistatic Evolution of Orlicz-Sobolev Nematic Elastomers]{Quasistatic Evolution of Orlicz-Sobolev Nematic Elastomers}
\author[M. Bresciani]{Marco Bresciani${}^{*}$}
\address{*\,Department Mathematik, Friedrich-Alexander-Universit\"{a}t Erlangen-N\"{u}rnberg, Cauerstrasse 11, 91058 Erlangen (DE)}
\email{marco.bresciani@fau.de}
\author[B. Stroffolini]{Bianca Stroffolini${}^\dagger$}
\address{$\dagger$\,Dipartimento di Matematica e Applicazioni, Universit\'{a} di Napoli ``Federico II'', Via Cintia, 80126 Napoli (ITA)}
\email{bstroffo@unina.it}

\begin{document}

\begin{abstract} 
We investigate the variational model for nematic elastomer proposed by Barchiesi and DeSimone with the director field defined on the deformed configuration under general growth conditions on the elastic density. This leads us to consider deformations in Orlicz-Sobolev spaces. Our work builds upon a previous paper by Henao and the Second Author, and  extends their analysis to the quasistatic setting. The overall strategy parallels the one  devised by the First author in the case of Sobolev deformations for a similar model in magnetoelasticity.  We prove two existence results for energetic solutions in the rate-independent setting.  The first result concerns quasistatic evolutions driven by time-dependent applied loads. For this problem, we establish suitable  Poincar\'{e} and trace inequalities in modular form to recover the coercivity of the total energy. The second result ensures the existence of quasistatic evolution for both time-depend applied loads and boundary conditions under physical confinement. In its proof, we follow the approach advanced by Francfort and Mielke  based on a multiplicative decomposition of the deformation gradient and we implement it for energies comprising terms defined on the deformed configuration. Both existence results rely on  a compactness theorem for sequences of admissible states with uniformly bounded energy which yields the strong convergence of the composition of the nematic fields with the corresponding deformations. While proving it, we show the regular approximate differentiability of Orlicz-Sobolev maps with suitable integrability, thus generalizing a classical result  for Sobolev maps due to Goffman and Ziemer. 	\EEE 
\end{abstract}

\maketitle


\section{Introduction}

\setlength{\parindent}{0pt}

\subsection{Motivation}

A variational model describing the mechanical behaviour of nematic elastomers has been proposed by Barchiesi and DeSimone in \cite{BDesimone}. According to this model, equilibrium configurations of such an elastomer with reference configuration $\Omega\subset \R^3$ correspond to minimizers of the energy functional 
\begin{equation}
	\label{eq:intro-I}
	(\boldsymbol{y},\boldsymbol{n})\mapsto  \int_\Omega W(D\boldsymbol{y}(\boldsymbol{x}),\boldsymbol{n}( \boldsymbol{y}(\boldsymbol{x}))\,\d\boldsymbol{x}+\int_{\boldsymbol{y}(\Omega)}|D\boldsymbol{n}(\boldsymbol{\xi})|^2\,\d\boldsymbol{\xi},
\end{equation}
where  $\boldsymbol{y}\colon \Omega \to \R^3$ is the elastic deformation and  $\boldsymbol{n}\colon \boldsymbol{y}(\Omega)\to \S^2$  is the nematic field describing for the local orientation of the constituent molecules. The first term in \eqref{eq:intro-I} accounts for the elastic energy and it depends on  the deformation gradient and the director field evaluated in the reference configuration. The elastic density $W$ takes the form
\begin{equation}
	\label{eq:intro-W}
	W(\boldsymbol{F},\boldsymbol{z})\coloneqq \Phi \left( \left (\mu^{-1} \boldsymbol{z}\otimes \boldsymbol{z}+\sqrt{\mu}(\boldsymbol{I}-\boldsymbol{z}\otimes \boldsymbol{z})\right )\boldsymbol{F} \right) \quad \text{for $\boldsymbol{F}\in\R^{3\times 3}$ and $\boldsymbol{z}\in \S^2$,}
\end{equation}
where $\Phi$ is frame-indifferent density which blows up as the determinant of its argument approaches zero and $\mu>0$ is a material parameter. The second term in \eqref{eq:intro-I} stands for the Oseen-Frank energy in the so-called one-constant approximation \cite{ball.lc}, here normalized to the unit. 
The minimization problem becomes nontrivial once  the deformation is prescribed on a portion of the boundary. When applied loads are considered, their work, which needs to be subtracted from the energy in \eqref{eq:intro-I}, is accounted by the functional
\begin{equation*}
	\label{eq:intro-II}
	(\boldsymbol{y},\boldsymbol{n})\mapsto \int_\Omega \boldsymbol{f}(\boldsymbol{x})\cdot \boldsymbol{y}(\boldsymbol{x})\,\d \boldsymbol{x}+\int_{\partial \Omega} \boldsymbol{g}(\boldsymbol{x})\cdot \boldsymbol{y}(\boldsymbol{x})\,\d \mathscr{H}^2(\boldsymbol{x})+\int_{\boldsymbol{y}(\Omega)} \boldsymbol{h}(\boldsymbol{\xi})\cdot \boldsymbol{n}(\boldsymbol{\xi})\,\d\boldsymbol{\xi},
\end{equation*}
where $\boldsymbol{f}\colon \Omega \to \R^3$, $\boldsymbol{g}\colon \partial \Omega \to \R^3$, and $\boldsymbol{h}\colon \R^3\to \R^3$ are dead loads representing body and surface forces, and external fields, respectively. 

The fact that the nematic field $\boldsymbol{n}$ is defined on the unknown deformed configuration $\boldsymbol{y}(\Omega)$ represents the characteristic feature of the model in \cite{BDesimone} and it becomes particularly relevant when large deformations are expected. Similar formulations, possibly including additional lower-order terms, appear also in other context such as the modeling of magnetoelastic solids \cite{bresciani,bresciani.davoli.kruzik,kruzik.stefanelli.zeman}. The same feature also constitutes the main source of difficulties for the analysis. These mainly concern  the coercivity of the energy and the lower semicontinuity of the elastic term in view of its dependence on the composition $\boldsymbol{n}\circ \boldsymbol{y}$. 

To ensure the well-posedness of the model, one aims at establishing the existence of minimizers of the energy. Next, one might aim at investigating the model under time evolution starting from the simpler quasistatic case. In this setting, the theory of rate-independent systems \cite{MieRu} has reached a prominent role in the description of variational models in continuum mechanics. Within this framework, energetic solutions, which are formulated on the two principles of global minimality and energy-dissipation balance, represent the most basic notion of evolution which has been successfully applied to many different models (see, e.g., \cite[Chapter 4]{MieRu}).

Both questions of the existence of minimizers and energetic solutions for the model in \cite{BDesimone} have been addressed in many papers by progressively enlarging the class of admissible deformations. We provide a brief review of the most relevant literature in the next subsection. In this regard, we observe that integrability of admissible deformations correspond to  the coercivity assumptions on the elastic density $W$ (or equivalently, recalling \eqref{eq:intro-W}, the function $\Phi$). From the mechanical point of view, it is clear that the  mechanical response of the elastomer should depend as little as possible on the behaviour of the elastic density for large values. Said differently, in order for the analysis to be compatible with a broader range of constitutive models, relaxing the coercivity assumptions on the elastic energy density becomes crucial. 

With this motivation, in this paper we further develop the study in \cite{HS} and we postulate general growth conditions, possibly not of polynomial type, on the elastic density. Thus, we work with deformations in Orlicz-Sobolev spaces as already  done in  \cite{ball.convexity}. The present analysis focuses on the corresponding quasistatic model.
 
\subsection{Review of the literature}
Let $N\in\N$ with $N\geq 2$ denote the  space dimension. The existence of minimizers of the energy  has been first proved by Barchiesi and DeSimone in \cite{BDesimone} for deformations in $W^{1,N}(\Omega;\RN)$ by restricting to incompressible materials, i.e., by allowing only deformations with $\det D \boldsymbol{y}=1$. An alternative proof for deformations in $W^{1,p}(\Omega;\RN)$ with $p>N$, still in the incompressible case, has been devised by Kru\v{z}\'{i}k, Stefanelli, and Zeman in \cite{kruzik.stefanelli.zeman} for a similar model of magnetoelasticity. In the same paper,  the existence of rate-independent evolutions has also been established. Both results have been extended by  Davoli, Kru\v{z}\'{i}k, and the First Author in \cite{bresciani.davoli.kruzik} to the compressible case for the same regime $p>N$.

Subsequently, Barchiesi, Henao and Mora-Corral have proved in \cite{BHM17} the existence of minimizers for deformations $\boldsymbol{y}\in W^{1,p}(\Omega;\RN)$ with $p>N-1$ and  $\det D \boldsymbol{y}\in L^1_+(\Omega)$ excluding cavitation. This last condition is formalized by requiring that admissible deformations satisfy the divergence identities \cite{mueller.tang.yan}:
\begin{equation*}
	\div \left( (\adj D \boldsymbol{y})\boldsymbol{\psi}\circ \boldsymbol{y} \right)=(\div \boldsymbol{\psi})\circ \boldsymbol{y}\,\det D \boldsymbol{y} \quad \text{for all $\boldsymbol{\psi}\in C^\infty_{\rm c}(\RN;\RN)$.}
\end{equation*}
For the same class of admissible deformations as in \cite{BHM17}, the existence of quasistatic evolutions has been achieved by the First Author in \cite{bresciani}. There, the model in \cite{kruzik.stefanelli.zeman,bresciani.davoli.kruzik} is investigated under the constraint $|\boldsymbol{n}\circ \boldsymbol{y}|\det D \boldsymbol{y}=1$. 

In \cite{HS}, Henao and the Second Author have generalized the analysis in \cite{BHM17} by working with a larger class of deformations in the Orlicz-Sobolev space $W^{1,A}(\Omega;\RN)$, where $A\colon [0,+\infty) \to [0,+\infty)$ is a suitable convex function satisfying
\begin{equation}
	\label{eq:intro-int}
	\int_1^{+\infty} \left(\frac{s}{A(s)}\right)^{\frac{1}{N-2}}\,\d s<+\infty.
\end{equation}
Roughly speaking, in this setting, the integrability $D\boldsymbol{y}\in L^p(\Omega;\RNN)$ is replaced by the requirement
\begin{equation*}
	\int_\Omega A(|D\boldsymbol{y}(\boldsymbol{x})|)\,\d\boldsymbol{x}<+\infty.
\end{equation*}
Thus, the condition \eqref{eq:intro-int} relaxes the assumption $p>N-1$  in \cite{BHM17} to the scale of Orlicz-Sobolev spaces. By exploiting an optimal embedding for Orlicz-Sobolev maps due to Cianchi \cite{Cianchi96,cianchi.embedding} and  fine properties of weakly differentiable maps with gradients in Lorentz spaces \cite{KKM,Maly.ac,Maly,O}, together with an oscillation estimate for Orlicz-Sobolev maps on spheres from \cite{CarozzaCianchi16}, the Authors  have shown in \cite{HS} that all the monotonicity, regularity, and invertibility properties established in \cite{BHM17} for deformations in $W^{1,p}(\Omega;\RN)$ are still valid for maps in $W^{1,A}(\Omega;\RN)$. Then, with  these results at hand, they have been able to adapt the arguments in \cite{BHM17} and prove the existence of minimizers for deformations in $W^{1,A}(\Omega;\RN)$. 

Eventually, the existence result in \cite{BHM17} has been recently extended to the case of deformation in $W^{1,p}(\Omega;\RN)$ with $p>N-1$ possibly creating cavities in \cite{bresciani.friedrich.moracorral}. 
We conclude by mentioning two papers concerning the relaxation of  functionals as in \eqref{eq:intro-I} with more general nematic term  in the Sobolev setting \cite{moracorral.oliva} and in the Orlicz-Sobolev (more precisely, Zygmund-Sobolev) setting \cite{scilla.stroffolini}.

\subsection{Main results}
This work aims at extending the analysis in \cite{HS} to the quasistatic setting by adopting the framework of rate-independent systems \cite{MieRu}.
We consider the same class of deformations as in \cite{HS} given by maps in $W^{1,A}(\Omega;\RN)$,  where $A$ is an N-function (see Definition \ref{def:N-function} below) satisfying condition ($\Delta_2$) (see Definition \ref{def:doubling} below) and fulfilling \eqref{eq:intro-int} for $N\geq 3$ (no integrability condition for $A$ is required when $N=2$). To address the quasistatic setting, we follow the approach devised in \cite{bresciani}.

First, we prove a compactness result for sequences of admissible states with uniformly bounded energy that, in particular, yields the strong convergence of the compositions of the nematic fields with the corresponding deformations. This result is given Theorem \ref{thm:compactness} which generalizes \cite[Theorem 3.2]{bresciani} to the Orlicz-Sobolev case. Its proof goes along the lines of the one in \cite{bresciani} and relies on a construction involving the topological image of nested balls (see Lemma \ref{lem:nested-balls}). The latter exploits the regular approximate differentiability of deformations which we establish in Theorem \ref{thm:reg-approx-diff}, thus generalizing a classical result by Goffman and Ziemer \cite[Theorem 3.4]{goffman.ziemer} for  Sobolev maps to the Orlicz-Sobolev setting. The proof of Theorem \ref{thm:reg-approx-diff} employs an oscillation estimate for Orlicz-Sobolev maps on spheres from \cite{CarozzaCianchi16} as in \cite{HS}.
Because of the constraint $|\boldsymbol{n}|=1$ (recall that $|\boldsymbol{n}\circ \boldsymbol{y}|\det D\boldsymbol{y}=1$  in \cite{bresciani}), we need to modify the proof in \cite{bresciani}  by using the weak continuity of the Jacobian determinant of inverse deformations achieved in \cite[Proposition 4.12(c)]{HS} (to be compared with \cite[Theorem 6.3(c)]{BHM17}).   

Next, we investigate the existence of energetic solutions. Given the expression of the elastic energy in \eqref{eq:intro-I}, a natural choice for the dissipative variable is given by the composition of nematic field and deformation. 
We address two settings for each of which we prove an existence result.

The first result is Theorem \ref{thm:existence-ti} which provides the existence of quasistatic evolutions driven by time-dependent applied loads under time-independent boundary conditions. In this setting, the main point is to establish the coercivity of the total energy under Orlicz growth conditions on $W$. Given the intricate definition of the Luxemburg norm, this is not a trivial task. To overcome this issue, we establish suitable  modular Poincar\'{e} and the trace inequalities with trace term  in Proposition \ref{prop:poincare} and Proposition \ref{prop:trace}, respectively. To our knowledge, similar results have been known in the literature for homogeneous boundary conditions only (see, e.g.,  \cite[Corollary 7.4.1(a)]{hh} or \cite[Proposition 2.13]{MSZ}). Note that the existence result in \cite{HS} did not account for applied loads, so that our findings provide an extension also for the static problem.
We mention that, if $A$ satisfies an additional property, namely, the condition ($\nabla_2$), then the Luxemburg norm associated to $A$ can be controlled from above and from below by a suitable root of the modular associated to $A$ (see, e.g., \cite[Lemma 3.2.9]{hh}), so that the coercivity issue can be easily resolved. Under this additional assumption, it is actually possible to improve Theorem \ref{thm:existence-ti} by weakening the integrability assumptions on the applied forces by means of the embeddings available for Orlicz-Sobolev functions \cite{cianchi.embedding} and their traces \cite{cianchi.trace}. In this work we tried to be as general as possible and we only required $A$ to satisfy  ($\Delta_2$). As a result, our analysis is applicable,  e.g., for $A(s)=s\log(\mathrm{e}+s)$, which does not satisfy ($\nabla_2$), when $N=2$ (in this case, assumption \eqref{eq:intro-int} can be disregarded).

The second result is Theorem \ref{thm:existence-td} which concerns quasistatic evolutions for both time-dependent applied loads and boundary conditions. 
 In this case, similarly to \cite{FraMie}, we reduce to a problem with time-independent boundary conditions by considering deformations and nematic field of the form $\boldsymbol{y}=\boldsymbol{d}_t \circ \boldsymbol{u}$ and $\boldsymbol{n}=\boldsymbol{m}\circ \boldsymbol{d}_t^{-1}$, where $t\mapsto \boldsymbol{d}_t$ is the boundary datum, $\boldsymbol{u}$ is a deformation that coincides with the identity on a portion of the boundary, and $\boldsymbol{m}$ is a sphere-valued map defined on the image of $\boldsymbol{u}$. In this way, energetic solutions to our problem correspond to energetic solutions to an auxiliary problem with time-independent boundary conditions. Theorem \ref{thm:existence-a} provides the existence of quasistatic evolutions to the auxiliary problem.  For technical reasons, we need to impose a physical confinement. Our regularity assumptions on the boundary data coincide with the ones in \cite{Laz}. If we exclude applied loads and we consider evolutions driven by boundary conditions only, we can remove the confinement and we can assume the regularity $W^{1,\infty}(0,T;\RN;\RN)\cap \Lip([0,T];\Lip(\RN;\RN))$, where $T>0$ denotes the time horizon, for both the boundary datum and its (spatial) inverse.

By now, the approach to the existence theory for energetic solutions to rate-independent systems has become standard, at least for continuous dissipation distances as the one considered here.  For this reason,  in the proofs of Theorem \ref{thm:existence-ti} and Theorem \ref{thm:existence-a} we resorted to an abstract existence result  available in the literature (here recalled in Theorem \ref{thm:MR}) in order to avoid the repetitions of well-known arguments.

Eeventually, in this paper we have been working with locally injective deformations (see Proposition \ref{prop:loc-inv}). However, the restriction to globally injective deformations is possible, e.g., by means of the  Ciarlet-Ne\v{c}as condition \cite{dalmaso.lazzaroni,Laz}.

\subsection{Structure of the paper} The paper contains five sections, the first one being this introduction, and an appendix. Section \ref{sec:notation} provides the necessary background on approximately differentiable maps, Orlicz-Sobolev spaces, and rate-independent systems, and reviews some results from \cite{HS}. Section \ref{sec:compactness} is devoted to the main compactness result. Quasistatic evolutions for time-independent and time-dependent boundary data are addressed in Section \ref{sec:ti} and \ref{sec:td}, respectively. Eventually, some technical results are collected in the final appendix. 

\EEE 

\section{Preliminaries}\label{sec:notation}

\subsection{Notation}
The space dimension $N\in\N$ satisfies $N\geq 2$.  The symbol $\RNN$ denotes the set of square matrices of order $N$ with $\boldsymbol{I}$ and $\boldsymbol{O}$ indicating the identity and the null matrix, respectively.
Given $\boldsymbol{F}\in\RNN$, the adjugate matrix $\adj \boldsymbol{F}\in\RNN$ is uniquely determined by the formula $\boldsymbol{F}(\adj\boldsymbol{F})=(\det \boldsymbol{F})\boldsymbol{I}$, where  $\det \boldsymbol{F}$ stands for the determinant of $\boldsymbol{F}$. The cofactors matrix of $\boldsymbol{F}$ is defined as $\cof \boldsymbol{F}\coloneqq (\adj \boldsymbol{F})^\top$ with the superscript $\top$ denoting the transposition. We indicate by $\RNN_+$ the set of matrices $\boldsymbol{F}\in\RNN$ with $\det \boldsymbol{F}>0$. We employ the notation
\begin{equation*}
	\mathbf{M}(\boldsymbol{F})=\textstyle (\mathbf{M}_1(\boldsymbol{F}),\dots,\mathbf{M}_N(\boldsymbol{F}))\in\prod_{i=1}^N \R^{\binom{N}{i}\times \binom{N}{i}},
\end{equation*}
where $\mathbf{M}_i(\boldsymbol{F})$ collects all $i$-th order minors of $\boldsymbol{F}$. The symbol $\boldsymbol{id}$ is used for the identity map on $\RN$.

The symbol $\chi_E$ is used for the characteristic function of a set $E\subset \RN$.
We denote the boundary and the closure of $E$ as $\partial E$ and $\closure{E}$.
Given  $F\subset \RN$, we write $E \subset \subset F$ whenever $\closure{E}\subset F$. The symbol $B(\boldsymbol{x},r)$ is used for the open ball centred at $\boldsymbol{x}\in\RN$ with radius $r>0$. We set  $S(\boldsymbol{x},r)\coloneqq \partial B(\boldsymbol{x},r)$ and $\closure{B}(\boldsymbol{x},r)\coloneqq \closure{B(\boldsymbol{x},r)}$. The unit sphere centred at the origin is indicated as $\S^{N-1}\coloneqq S (\boldsymbol{0},1)$.

The Lebesgue measure and the $\alpha$-dimensional Hausdorff measure, where $\alpha>0$, on $\RN$ are denoted by $\leb$ and $\mathscr{H}^\alpha$. Analogously, the symbol $\mathscr{L}^1$ is used for the Lebesgue measure on $\R$. Otherwise differently stated, the adjective ``measurable'' and ``negligible'', and the expression ``almost everywhere'' (and all similar ones) , which is abbreviated as ``a.e.'' within the formulas, are referred to the measure $\leb$. The integration with respect to $\leb$ or $\mathscr{L}^1$ is indicated by the differential of the integration variable. The dashed integral symbol denotes the integral average.

Given a measure space $(X,\mathfrak{M},\mu)$ and $Y\in\mathfrak{M}$, the restriction of $\mu$ to $Y$ is written as $\mu \mres Y$. The symbol $L^0(\mu)$ denotes the class of $\mu$-measurable maps on $X$ taking real values. The Lebesgue spaces $L^p(\mu)$ and $L^p_{\rm loc}(\mu)$ with  $1 \leq p \leq \infty$ are defined in the usual manner. Spaces of continuous functions are indicated with the symbol $C^0$. Concerning (weakly) differentiable maps, we employ the standard notations $C^k$ with $k\in \N\cup \{\infty\}$ for continuously-differentiable maps and  their compactly supported counterpart $C^k_{\rm c}$, $\Lip$ for Lipschitz maps, and   $W^{1,p}$ for Sobolev maps. The weak gradient is indicated by $D$. Traces of Sobolev maps are denotes as restrictions.

Weak and weak-* convergence in Banach spaces are indicated by the symbols $\wk$ and $\wks$. For $X$ and $Y$ Banach spaces, the notation $X\hookrightarrow Y$ and $X\hookrightarrow \hookrightarrow Y$ are used for continuous and compact embeddings, respectively.
Given $T>0$, we write $L^p(0,T;X)$ and $W^{1,p}(0,T;X)$ for the Bochner space and the Bochner-Sobolev space of functions from $(0,T)$ to $X$. For the latter, the weak derivative with respect to the time variable $t\in (0,T)$ is indicated by a superposed dot. The symbols $AC([0,T];X)$ and $\Lip([0,T];X)$ are used for spaces of absolutely continuous and Lipschitz functions with values in $X$. In general, for functions  defined on $(0,T)$ or $[0,T]$, the value corresponding to $t$ is indicated with the subscript $t$. 

We adopt the usual convention of denoting by $C,C_1,C_2$, and so on positive constants that may change from line to line. Sometimes we highlight the dependence of this constant on some quantities or objects by using parentheses. 

\subsection{Approximately differentiable maps}

Let $E\subset \RN$ be a measurable set. The density of $E$ at a point $\boldsymbol{x}_0\in\RN$ is defined as the number
\begin{equation*}
	\Theta^N(E,\boldsymbol{x}_0)\coloneqq \lim_{r \to 0^+} \frac{\leb(E \cap B(\boldsymbol{x}_0,r))}{\leb(B( \boldsymbol{x}_0 ,r))},
\end{equation*}
whenever the limit exists. 
For more information concerning the notions of approximate limits, continuity, and differentiability, we refer to \cite[Chapter 3]{giaquinta.modica.soucek}. 

\begin{definition}[Geometric domain and image]
Let $\boldsymbol{y}\colon \Omega \to \RN$ be almost everywhere approximately differentiable	with $\det \nabla \boldsymbol{y}>0$  almost everywhere.
The {geometric domain} $\domg(\boldsymbol{y},\Omega)$ of $\boldsymbol{y}$ is defined as the set of points $\boldsymbol{x}_0\in \Omega$ such that $\boldsymbol{y}$ is approximately differentiable at $\boldsymbol{x}_0$ with $\det \nabla \boldsymbol{y}(\boldsymbol{x}_0)> 0$, and there exists a compact set $K\subset \Omega$ with $\Theta^N(K,\boldsymbol{x}_0)=1$ and a function $\boldsymbol{w}\in C^1(\RN;\RN)$ satisfying $\boldsymbol{w}\restr{K}=\boldsymbol{y}\restr{K}$ and $D\boldsymbol{w}\restr{K}=\nabla \boldsymbol{y}\restr{K}$.
Moreover, for every  $E\subset \Omega$ measurable, we set $\domg(\boldsymbol{y},E)\coloneqq \domg(\boldsymbol{y},\Omega) \cap E$ and we define  the  {geometric image} of $E$ under $\boldsymbol{y}$ as $ \img(\boldsymbol{y},E)\coloneqq \boldsymbol{y} \left( \domg(\boldsymbol{y},E) \right)$.
\end{definition}

In this work, we will employ the following version of Federer's  area formula \cite[Theorem 1, p.~220]{giaquinta.modica.soucek}. Recall that a map $\boldsymbol{y}\colon \Omega \to \RN$ is termed {almost everywhere  injective} in $E\subset \Omega$ if there exists a set $X\subset E$ with $\leb(X)=0$ such that $\boldsymbol{y}\restr{E \setminus X}$ is injective.  

\begin{proposition}[Area formula]
	\label{prop:federer}
Let $\boldsymbol{y}\colon \Omega \to \RN$ be almost everywhere approximately differentiable with $\det \nabla \boldsymbol{y}>0$ almost everywhere. Suppose that $E\subset \Omega$ is measurable and $\boldsymbol{y}$ is almost everywhere injective in $E$. Then, there holds
\begin{equation*}
	\leb\big (\img(\boldsymbol{y},E)\big )= \int_E \det  \nabla   \boldsymbol{y}(\boldsymbol{x})\,\d\boldsymbol{x}.
\end{equation*}
\end{proposition}

As a consequence Proposition~\ref{prop:federer}, $\boldsymbol{y}$ satisfies Lusin's condition (N${}^{-1}$), namely, there holds $\leb(\boldsymbol{y}^{-1}(F))=0$ for every $F\subset \RN$ with $\leb(F)=0$. Here,  $\boldsymbol{y}^{-1}(F)\subset \Omega$ denotes the preimage of $F$ under $\boldsymbol{y}$. A simple proof of this fact is given, e.g., in \cite[Remark 2.3(b)]{bresciani.friedrich.moracorral}.

The following change-of-variable formula  is a  consequence of   Proposition \ref{prop:federer}. 

\begin{corollary}[Change-of-variable formula]
	\label{cor:cov}
	Let $\boldsymbol{y}\colon \Omega \to \RN$ be almost everywhere approximately differentiable with $\det \nabla \boldsymbol{y}>0$ almost everywhere.  Suppose that $E\subset \Omega$ is measurable,  $\boldsymbol{y}$ is almost everywhere injective in $E$, and  $\boldsymbol{\psi}\colon \img(\boldsymbol{y},E)\to \R^\nu$ with $\nu\in\N$ is measurable. Then, there holds
	\begin{equation*}
		\int_{\img(\boldsymbol{y},E)} \boldsymbol{\psi}(\boldsymbol{\xi})\,\d\boldsymbol{\xi}=\int_E \boldsymbol{\psi}(\boldsymbol{y}(\boldsymbol{x})) \,\det \nabla \boldsymbol{y}(\boldsymbol{x})\,\d\boldsymbol{x},
	\end{equation*}
	whenever one of the two integrals exists. 
\end{corollary}

The next result is a reformulation of \cite[Lemma~2.25]{HeMo11} and combines it with Corollary \ref{cor:cov}. 

\begin{lemma}[Approximate differentiability of the inverse]
	\label{lem:injectivity}
Let $\boldsymbol{y}\colon \Omega \to \RN$ be almost everywhere approximately differentiable	with $\det \nabla \boldsymbol{y}>0$ almost everywhere.  Suppose that $E\subset \Omega$ is  measurable and $\boldsymbol{y}$ is almost everywhere injective in $E$.   Then,   $\boldsymbol{y}\restr{\domg(\boldsymbol{y},E)}$ is injective. Moreover, its inverse $(\boldsymbol{y}\restr{\domg(\boldsymbol{y},E)})^{-1}$ is approximately differentiable  in $\img(\boldsymbol{y},E)$  with 
\begin{equation}
	\label{eq:i}
	\nabla \boldsymbol{y}^{-1}(\boldsymbol{\xi})=(\nabla \boldsymbol{y})^{-1}(\boldsymbol{y}^{-1}(\boldsymbol{\xi})) \quad \text{for all $\boldsymbol{\xi}\in \img(\boldsymbol{y},E)$.}
\end{equation}
In particular,  for every $\boldsymbol{\varphi}\colon E\to \R^\nu$ measurable with $\nu\in\N$, there holds
\begin{equation}
	\label{eq:ii}
	\int_E \boldsymbol{\varphi}(\boldsymbol{x})\,\d\boldsymbol{x}=\int_{\img(\boldsymbol{y},E)}\boldsymbol{\varphi}(\boldsymbol{y}^{-1}(\boldsymbol{\xi}))\det \nabla \boldsymbol{y}^{-1}(\boldsymbol{\xi})\,\d\boldsymbol{\xi}
\end{equation}
whenever one of the two integrals exists. Here, in \eqref{eq:i}--\eqref{eq:ii}, we simply write $\boldsymbol{y}^{-1}$ in place of $(\boldsymbol{y}\restr{\domg(\boldsymbol{y},E)})^{-1}$.
\end{lemma}

 \subsection{Orlicz--Sobolev spaces}
 
We recall the main notions of the theory of Orlicz spaces. For a comprehensive treatment in the case of functions defined on Euclidean domains and abstract measure spaces, we refer to \cite{Adams,Kras,Kufn} and \cite{rao.ren}, respectively.

\begin{definition}[N-function]\label{def:N-function}
A function $A:[0,+\infty)\to[0,+\infty)$ is termed an {N-function} if it convex  and satisfies 
		\begin{equation*}
			A(0)=0, \qquad A(s)>0 \quad \text{for all $s>0$},\qquad  \lim_{s \to 0^+} \frac{A(s)}{s}=0, \qquad \lim_{s \to +\infty} \frac{A(s)}{s}=+\infty.
		\end{equation*}
\end{definition}
Observe that  $A$ takes finite values and, hence, it is locally Lipschitz.
Moreover,  $A$ is strictly increasing, so that  its inverse $A^{-1}$ is well defined. We  denote the left-derivative (which exists everywhere) and the derivative (which exists almost everywhere) of $A$ by $A'_-$ and $A'$, respectively. 

We introduce an order relation on N-functions. 
  
\begin{definition}[Ordering of N-functions]
Let $A$ and $B$ be two N-functions. 
\begin{enumerate}[label=(\roman*)]
	\item  $B$ is termed to {dominate} $A$ near infinity, denoted by  $A \prec B$, whenever there exist $c>0$ and $s_0\geq 0$ such that $A(s)\leq B(cs)$ for all $s\geq s_0$. 
	\item  $A$ and $B$ are termed to be {equivalent}, denoted by $A \sim B$, whenever both $A \prec B$ and $B \prec A$ hold. 
\end{enumerate}
\end{definition}

Next, we present a doubling condition for N-functions.

\begin{definition}[Condition $\boldsymbol{\Delta_2}$]\label{def:doubling}
An N-function $A$ is said to satisfy  {condition} ($\Delta_2$) if there exists a constant $\kappa_A >0$ such that 
\begin{equation}
	\label{eqn:Delta2}
	A(2s) \leq \kappa_A A(s) \quad \text{for all $s\geq 0$.}
\end{equation}
\end{definition}
By convexity,  one must have  $\kappa_A\geq 2$.
In the literature, $A$ is termed to satisfy  $(\Delta_2)$ near infinity  whenever the exists $s_0>0$ such that $A(2s)\leq\kappa_A A(s)$ holds for all $s\geq s_0$, while  $A$ is termed to satisfy $(\Delta_2)$ globally whenever \eqref{eqn:Delta2} holds.
 Henceforth, we will always refer to the global condition. However, this fact does not pose any restriction in the following sense \cite[p.~24]{Kras}:  if $A$ satisfies  $(\Delta_2)$ near infinity, then, there exists another N-function  $B$ with $A \sim B$ which satisfies  $(\Delta_2)$ globally.

\begin{lemma}\label{lem:p}
Let $A$ be an N-function satisfying  {\rm ($\Delta_2$)}. Then, 
there exists $p_A>1$ such that 
	\begin{equation}\label{eqn:p}
		\frac{s\,A'_-(s)}{A(s)} \leq p_A \quad \text{for all $s>0$}.
	\end{equation}
	Moreover
		\begin{equation}\label{eqn:pp}
			A(c s)\leq c^{\,p_A} A(s) \quad \text{for all $c>1$ and $s\geq 0$.}
		\end{equation}
\end{lemma}
\begin{proof}
The existence of $p_A>1$ satisfying \eqref{eqn:p} is proved in 	\cite[Theorem 3(1), p.~23]{rao.ren}. Using this estimate, the inequality in \eqref{eqn:pp} follows as in \cite[p.~36]{MSZ}.
\end{proof}

Next, we introduce the conjugate of an N-function. This is nothing but the usual conjugate (also known as Legendre transform) of a convex function.

\begin{definition}[Conjugate function]
Let $A$ be an N-function. The function $\bar{A}\colon [0,+\infty) \to [0,+\infty]$ defined as
\begin{align*}
	\bar{A}(s)\coloneqq\sup\{s\sigma-A(\sigma):\sigma\geq 0\}\quad \text{for all $s\geq 0$}
\end{align*}
is termed the {conjugate function}  of $A$.
\end{definition}

The conjugate  of an N-function $A$ is also an N-function and it is known that $\bar{{\bar A}}=A$.  
 From the definition of $\bar{A}$, we immediately deduce {Young's inequality}:
 \begin{equation}
 	\label{eqn:young}
 	s_1 s_2 \leq A(s_1)+\bar{A}(s_2) \quad \text{for all $s_1,s_2\geq 0$.}
 \end{equation}
 
 We now introduce Orlicz spaces.

\begin{definition}[Orlicz space]
	Let  $A$ be an N-function and  $(X,\mathfrak{M},\mu)$ be a measure space.  The corresponding {Orlicz space} is defined as
	\begin{align*}
		L^A(\mu)\coloneqq \left\{ v\in L^0(\mu):\hspace{3pt} \int_X A\left (\frac{|v|}{s}\right )\,\d\mu<+\infty \hspace{5pt} \text{for some $s >0$}   \right\}
	\end{align*}
	with the {Luxemburg norm} given by
	\begin{align}
		\label{eqn:luxemburg}
		\|v\|_{L^A(\mu)}&\coloneqq \inf \left\{ s >0: \hspace{4pt} \int_X A \left( \frac{|v|}{s } \right)\,\d\mu\leq 1 \right\}\quad \text{for all $v\in L^A(\mu)$.}
	\end{align} 
\end{definition}  

The class $L^A(\mu)$ forms a vector subspace of $L^1(\mu)$. If $A$ satisfies  $(\Delta_2)$, then one has $v\in L^A(\mu)$ if and only if $A(|v|)\in L^1(\mu)$.  
The functional $v\mapsto \int_X A(|v|)\,\d \mu$ is termed the {modular} associated to $A$. 
The expression in \eqref{eqn:luxemburg} defines  indeed a norm that makes $L^A(\mu)$ a Banach space with $L^A(\mu)\hookrightarrow L^1(\mu)$. In particular, the infimum in \eqref{eqn:luxemburg} is actually achieved, so that 
\begin{equation}
	\label{eqn:lux-min}
	\int_X A \left( \frac{|v|}{\|v\|_{L^A(\mu)}} \right)\,\d\mu \leq 1 \quad \text{for all $v\in L^A(\mu)$.}
\end{equation}
From \eqref{eqn:luxemburg}, we immediately have
\begin{equation}\label{eqn:norm-modular}
	\|v\|_{L^A(\mu)} \leq \int_X A(|v|)\,\d\mu + 1 \quad \text{for all $v\in L^A(\mu)$.}
\end{equation}
When $A$ satisfies  $(\Delta_2)$, it is  also possible to control the modular with the Luxemburg norm. Indeed, using Lemma \ref{lem:p}(ii), the monotonicity of $A$, and \eqref{eqn:lux-min}, we obtain
\begin{equation}
	\label{eqn:modular-norm}
	\begin{split}
		\int_X A(|v|)\,\d \mu  &\leq (\|v\|_{L^A(\mu)}+1)^{p_A} \int_X A\left( \frac{|v|}{1+\|v\|_{L^A(\mu)}} \right)\,\d\mu\\
		&\leq (\|v\|_{L^A(\mu)}+1)^{p_A} \int_X A\left( \frac{|v|}{\|v\|_{L^A(\mu)}}  \right)\,\d\mu  \leq (\|v\|_{L^A(\mu)}+1)^{p_A}.
	\end{split}
\end{equation}
From \eqref{eqn:norm-modular}--\eqref{eqn:modular-norm}, we see that a sequence $(v_n)$ in $L^A(\mu)$ is bounded if and only if the sequence $\big(A(|v_n|)\big)$ in $L^1(\mu)$ is also bounded. 

 The next result  ensures the equivalence between norm and modular convergence in $L^A(\mu)$ provided that  $A$ satisfies ($\Delta_2$). We refer to \cite[Theorem 12(a), p.~83]{rao.ren} for the proof. 

\begin{proposition}[Norm and modular convergence]
	\label{prop:norm-modular-conv}
	Let $A$ be an N-function satisfying {\rm ($\Delta_2$)} and let $(X,\mathfrak{M},\mu)$ be a measure space. Then, for a sequence  $(v_n)$ in $L^A(\mu)$ and $v\in L^A(\mu)$, we have that 
	\begin{equation*}
		\label{eq:norm-modular}
		\text{$v_n \to v$ in $L^A(\mu)$} \qquad \text{if and only if} \qquad \text{$\int_X A(|v_n-v|)\,\d\mu \to 0$.}
	\end{equation*} 
\end{proposition}

In the framework of Orlicz spaces, H\"{o}lder's inequality takes the following form \cite[Proposition 1, p.~58]{rao.ren}. 

\begin{proposition}[H\"{o}lder inequality]
	Let $A$ be an N-function and let $(X,\mathfrak{M},\mu)$ be a measure space. Then, we have 
	\begin{equation}\label{eqn:hoelder}
		\int_X |vw|\,\d\mu \leq 2 \|v\|_{L^A(\mu)}\|w\|_{L^{\bar{A}}(\mu)} \quad \text{for all $v\in L^A(\mu)$ and $w\in L^{\bar{A}}(\mu)$.}
	\end{equation}
\end{proposition}

Concerning the duality theory  of Orlicz spaces, we mention a few basic facts that will be relevant for our analysis. We use prime to denote dual spaces. 

\begin{proposition}[Duality of Orlicz spaces]
	\label{prop:duality}
Let $A$ be an N-function and let $(X,\mathfrak{M},\mu)$ be a measure space. Define the class
\begin{equation*}
	M^{\bar{A}}(\mu)\coloneqq  \left\{ v\in L^0(\mu):\hspace{3pt} \int_X \bar{A}\left (\frac{|v|}{s}\right )\,\d\mu<+\infty \hspace{5pt} \text{for all $s>0$}   \right\},
\end{equation*}
where $\bar{A}$ denotes the conjugate N-function of $A$.
Then, we have the following:
\begin{enumerate}[label=(\roman*)]
	\item $M^{\bar{A}}(\mu)$ is a closed subspace of $L^{\bar{A}}(\mu)$ which coincides with the closure of simple functions in $L^{\bar{A}}(\mu)$; 
	\item If $(X,\mathfrak{M},\mu)$ is separable, then  $M^{\bar{A}}(\mu)$ is separable;
	\item If  $(X,\mathfrak{M},\mu)$ is $\sigma$-finite, then $\big (M^{\bar{A}}(\mu)\big )'=L^A(\mu)$ up to an isomorphism.
\end{enumerate}
\end{proposition}
\begin{proof}
Claim (i) follows from  \cite[Corollary 4, p. 77]{rao.ren}. 
The separability of $M^{\bar{A}}(\mu)$ in (ii) is proved in \cite[Corollary~5, p.~148]{rao.ren}. The representation of the dual in (ii)	is a consequence of  \cite[Theorem 7, p. 110]{rao.ren}  taking into account that $\bar{{\bar A}}=A$. 
\end{proof}

In view of Proposition \ref{prop:duality}(ii),  we can consider the weak-* topology  on the space $L^A(\mu)$. Given a sequence $(v_n)$ in $ L^A(\mu)$ and $v\in L^A(\mu)$, we have
\begin{equation*}
	\text{$v_n \wks v$ in $L^A(\mu)$} \qquad \text{if and only if} \qquad \int_X v_n\, w\,\d\mu \to \int_X v\,w\,\d \mu \quad \text{for all $w\in M^{\bar{A}}(\mu)$.}
\end{equation*}
In particular, by claims (ii) and (iii) of Proposition \ref{prop:duality}, every bounded sequence in $L^A(\mu)$ admits a subsequence which converges in the weak-* topology.

If  $\Omega\subset \R^N$ is an open set and $S\subset \R^N$ is an hypersurface of class $C^2$, we simplify the notation by writing $L^A(\Omega)$ and $L^A(S)$ for $L^A(\leb \mres \Omega)$ and $L^A(\haus \mres S)$, respectively.  The Orlicz spaces $L^A(\Omega;\R^N)$ of vector-valued maps is defined as the set of $\boldsymbol{v}\in L^1(\Omega;\R^N)$ with $|\boldsymbol{v}|\in L^A(\Omega)$ and it is equipped with the Luxemburg norm of $|\boldsymbol{v}|$.  Analogous definitions apply for $L^A(S;\RN)$ as well as  for  spaces of tensor-valued maps.

%
%
%

We introduce Orlicz-Sobolev spaces.

\begin{definition}[Orlicz-Sobolev space]
	Let $A$ be an N-function and let $\Omega\subset \RN$ be an open set. The corresponding {Orlicz-Sobolev space} 
	is defined as 
	\begin{equation*}
		W^{1,A}(\Omega) = \{v \in L^A(\Omega): \hspace{4pt} 
		D v \in L^A(\Omega;\R^N)\},
	\end{equation*}
	 with the  {Luxemburg-Sobolev norm} given by
		\begin{equation*}
			\|v\|_{W^{1,A}(\Omega)} \coloneqq \|v\|_{L^A(\Omega)} + \|D v\|_{L^A(\Omega; \R^N)} \quad \text{for all $v\in W^{1,A}(\Omega)$.}
		\end{equation*}
\end{definition}
Clearly, $W^{1,A}(\Omega)$ is a Banach space and $W^{1,A}(\Omega)\hookrightarrow W^{1,1}(\Omega)$. The space  $W^{1,A}(\Omega;\R^N)$ of vector-valued maps is defined in an analogous manner.

Given a domain $U \subset \subset \Omega$ of class $C^2$, the space $W^{1,A}(\partial U;\R^N)$ is defined by requiring the same regularity for the compositions of functions with local charts as for classical Sobolev spaces. We refer to \cite[Subsection 3.1]{scilla.stroffolini.inv} for  details. We shall denote the tangential gradient of $\boldsymbol{v}\in W^{1,A}(\partial U;\R^N)$ by $D^{\partial U}\boldsymbol{v}\in L^A(\partial U;\RNN)$.


\subsection{Growth conditions and divergence identities}
In this paper, we will consider an N-function $A$ satisfying the following integrability condition:
\begin{equation}
	\label{eqn:growth-infinity}
	\tag{I}
	\int_1^{+\infty} \left(\frac{s}{A(s)}\right)^{\frac{1}{N-2}}\,\d  s<+\infty.
\end{equation}
Examples of functions satisfying this condition are $A(s)=s^p$ for $p>N-1$ and $A(s)=s^{N-1}\log^q (\mathrm{e}  + s)$ for  $q>N-2$. Clearly, \eqref{eqn:growth-infinity} is meaningful only for $N\geq 3$. For $N=2$, this condition can be simply disregarded. 

We highlight some results from \cite{HS}. Henceforth, $\Omega\subset \RN$ will denote a bounded  domain.  We employ standard notation for Lorentz spaces, see, e.g., \cite{Adams,Kufn}.

\begin{proposition}[Embedding]
\label{prop:embedding}	
Let  $A$ be an N-function satisfying \eqref{eqn:growth-infinity}. Then, we have  the following:
\begin{enumerate}[label=(\roman*)]
	\item We have  $L^A(\Omega) \hookrightarrow L^{N-1,1}(\Omega)$ and $L^A(\Omega) \hookrightarrow L^{N-1}(\Omega)$. 
	\item For every domain $U\subset \subset \Omega$ of class $C^2$, we have    $W^{1,A}(\partial U;\RN) \hookrightarrow\hookrightarrow C^0(\partial U;\RN)$ up to the choice of a representative. 
\end{enumerate}
\end{proposition}
A direct proof for the second embedding in (i) can be promptly deduced from \cite[Lemma~3.1]{CarozzaCianchi16}.
\begin{proof}
Up to a suitable modification of  $A$ in the interval $(0,1)$ which yields an equivalent N-function, we may assume
\begin{equation*}
	m\coloneqq \int_0^{+\infty} \left(\frac{s}{A(s)}\right)^{\frac{1}{N-2}}\,\d  s<+\infty.
\end{equation*}
Arguing as in \cite[Proposition~2.6]{HS}, we deduce the existence of a constant $C>0$ such that
\begin{equation*}
	\|v\|_{L^{N-1,1}(\Omega)} \leq C m^{\frac{N-2}{N-1}} \left( \int_{\Omega} A(|v|)\,\d \boldsymbol{x} \right)^{\frac{1}{N-1}} \quad \text{for all $v\in L^A(\Omega)$.} 
\end{equation*}
This shows that $L^A(\Omega)\subset L^{N-1,1}(\Omega)$. Then,  Proposition \ref{prop:norm-modular-conv} yields the  embedding $L^A(\Omega) \hookrightarrow L^{N-1,1}(\Omega)$. The second embedding follows as $L^{N-1,1}(\Omega)\hookrightarrow L^{N-1}(\Omega)$.  Claim (ii) was proved in \cite[Proposition~2.6]{HS}.
\end{proof}

Given $\boldsymbol{v} \in L^1_{\loc}(\Omega;\R^N)$, we define its precise representative  $\boldsymbol{v}^* \colon \Omega\to \R^N$  as
\begin{equation*}
	\boldsymbol{v}^*(\boldsymbol{x}_0)\coloneqq \limsup_{r \to 0^+} \dashint_{B(\boldsymbol{x}_0,r)} \boldsymbol{v}(\boldsymbol{x})\,\d\boldsymbol{x}.
\end{equation*}
The set $L_{\boldsymbol{v}}$ of Lebesgue points of $\boldsymbol{v}$ is formed by all $\boldsymbol{x}_0\in \Omega$ such that
\begin{equation*}
	\lim_{r\to 0^+} \dashint_{B(\boldsymbol{x}_0,r)} |\boldsymbol{v}(\boldsymbol{x})-\boldsymbol{v}^*(\boldsymbol{x}_0)|\,\d\boldsymbol{x}=0.
\end{equation*}
In that case, the superior limit in the definition of $\boldsymbol{v}^*(\boldsymbol{x}_0)$ is actually a limit.
Concerning the Hausdorff dimension of the complement of the set of Lebesgue points, we have the following. 

\begin{lemma}[Exceptional points]
	\label{lem:exceptional}
Let $A$ be an N-function satisfying \eqref{eqn:growth-infinity}. Then, $\mathscr{H}^1(\Omega \setminus L_{\boldsymbol{v}})=0$ for all $\boldsymbol{v}\in W^{1,A}(\Omega;\RN)$. 	
\end{lemma}
\begin{proof}
By  Proposition~\ref{prop:embedding}(i), $D\boldsymbol{v}\in L^{N-1,1}(\Omega;\RNN)$. Thus, \cite[Theorem 3.7]{Maly} yields the conclusion.
\end{proof}

We need to introduce some notation. Let $U\subset \subset \Omega$ be a domain of class $C^2$. Then, there exists $\delta_U>0$ such that the signed distance function $d_U\colon \RN \to \R$ given by 
\begin{equation*}
	d_U(\boldsymbol{x})\coloneqq \begin{cases}
		-\dist(\boldsymbol{x};\partial U) & \text{if $\boldsymbol{x}\in \RN \setminus \closure{U}$,}\\
		\dist(\boldsymbol{x};\partial U) & \text{if $\boldsymbol{x}\in  \closure{U}$,}
	\end{cases}
\end{equation*}
is of class $C^2$  in the tubular neighborhood $\{ \boldsymbol{x}\in \RN:\: -\delta_U < d_U(\boldsymbol{x})<\delta_U \}\subset \subset \Omega$. In this case, the set $U_l\coloneqq \{ \boldsymbol{x}\in \RN:\: d_U(\boldsymbol{x})>l \}\subset \subset \Omega$ is a domain of class $C^2$ for all $l\in (-\delta_U,\delta_U)$. 

We recast \cite[Lemma~2.11]{HS} by considering precise representatives.

\begin{proposition}
	\label{prop:boundary}
Let $A$ be an N-function satisfying \eqref{eqn:growth-infinity} and let $U\subset \subset \Omega$ be a domain of class $C^2$. Also, let $(\boldsymbol{v}_n)_n$ be a sequence in $W^{1,A}(\Omega;\RN)$ and $\boldsymbol{v}\in W^{1,A}(\Omega;\RN)$ be such that $\boldsymbol{v}_n \wks \boldsymbol{v}$ in $W^{1,A}(\Omega;\RN)$. Then,  for almost every $l\in (-\delta_U,\delta_U)$, we have 
\begin{equation}
	\label{eqn:boundary}
	\partial U_l \subset L_{\boldsymbol{v}} \cap \bigcap_{n=1}^\infty L_{\boldsymbol{v}_n}, \qquad  \boldsymbol{v}_n^*\restr{\partial U_l}\in C^0(\partial U_l;\RN)\text{ for all $n\in\N$}, \qquad \boldsymbol{v}^*\restr{\partial U_l}\in  C^0(\partial U_l;\RN). 
\end{equation}
Moreover, along a not relabeled  subsequence possibly depending on $l$, we have $\boldsymbol{v}_n^* \to \boldsymbol{v}^*$ uniformly on $\partial U_l$.
\end{proposition}
\begin{proof}
Consider $\delta_U>0$ such that $d_U$ is of class $C^2$ in the tubular neighborhood of thickness $\delta_U$. 	
The first claim in \eqref{eqn:boundary} follows from Lemma~\ref{lem:exceptional}. 
As in \cite[Lemma~2.11]{HS}, we have $\boldsymbol{v}_n\in W^{1,A}(\partial U_l;\RN)$ for all $n\in\N$ and $\boldsymbol{v}\in W^{1,A}(\partial U_l;\RN)$ for almost all $l\in (-\delta_U,\delta_U)$. By Proposition~\ref{prop:embedding}(ii), each of these maps admits a continuous representative. By a mollification argument, we check that, for almost all $l\in (-\delta_U,\delta_U)$, this coincides with the restriction of the precise representative to $\partial U_l$. By arguing as in \cite[Lemma 2.11]{HS}, we show that $\boldsymbol{y}_n\wk \boldsymbol{y}$ in $W^{1,A}(\partial U;\RN)$ and the uniform convergence follows by Proposition \ref{prop:embedding}(ii).
\end{proof}

We recall the divergence identities.

\begin{definition}[Divergence identities]
A map $\boldsymbol{y}\in W^{1,A}(\Omega;\R^N)$ with  $\det D \boldsymbol{y}\in L^1(\Omega)$ is termed to satisfy the {divergence identities } if	
\begin{equation*}
	\tag{DIV}\label{DIV}
	\quad \mathrm{div}\left( (\adj D\boldsymbol{y})\boldsymbol{\psi}\circ \boldsymbol{y} \right)=\left( \div \boldsymbol{\psi} \right)\circ \boldsymbol{y}\det D \boldsymbol{y} \quad \text{for all $\boldsymbol{\psi}\in C^\infty_{\rm c}(\RN;\RN)$,}
\end{equation*}
where the divergence on the left-hand side  is understood in the sense of distributions, namely
\begin{equation*}
	- \int_\Omega \boldsymbol{\psi}\circ \boldsymbol{y} \cdot (\cof D \boldsymbol{y})D\varphi\,\d\boldsymbol{x}=
	\int_\Omega \left( \div \boldsymbol{\psi}\right)\circ \boldsymbol{y}  (\det D \boldsymbol{y})\,\varphi \,\d\boldsymbol{x} \quad \text{for all $\varphi \in C^\infty_{\rm c}(\Omega)$ and $\boldsymbol{\psi}\in C^\infty_{\rm c}(\RN
		;\RN)$.}
\end{equation*}
\end{definition}

Note that $\adj D \boldsymbol{y}\in L^1(\Omega;\RNN)$  by Proposition~\ref{prop:embedding}(i), so that  (\ref{DIV}) is meaningful.
We present a few consequences of the divergence identities. The first result concerns the weak continuity of Jacobian minors.

\begin{proposition}[Weak continuity of the Jacobian minors]\label{prop:weak-det}
	Let $(\boldsymbol{y}_n)$ be a sequence in $W^{1,A}(\Omega;\R^N)$ with each $\boldsymbol{y}_n$ satisfying {\rm (\ref{DIV})} and $\det D \boldsymbol{y}_n\in L^1_+(\Omega)$. Also, let $\boldsymbol{y}\in W^{1,A}(\Omega;\R^N)$, $\boldsymbol{\Theta}\in L^1(\Omega;\RNN)$,  and $\vartheta\in L^1(\Omega)$ be such that
	\begin{equation*}
		\text{$\boldsymbol{y}_n \wks \boldsymbol{y}$  in $W^{1,A}(\Omega;\RN)$,} \quad \text{$\cof D\boldsymbol{y}_n \wk \boldsymbol{\Theta}$ in $L^1(\Omega;\RNN)$,}\quad \text{$\det D \boldsymbol{y}_n\wk \vartheta$ in $L^1(\Omega)$.}
	\end{equation*}
	Then, $\boldsymbol{y}$ satisfies {\rm (\ref{DIV})}, $\boldsymbol{\Theta}=\cof D \boldsymbol{y}$, and $\vartheta=\det D \boldsymbol{y}$. Therefore
	\begin{equation*}
		\label{eqn:weak-minors}
		\text{$\mathbf{M}(D\boldsymbol{y}_n)\wk \mathbf{M}(D\boldsymbol{y})$ in $L^1(\Omega;\textstyle \prod_{i=1}^N\R^{\binom{N}{i}\times \binom{N}{i}})$.}
	\end{equation*}
\end{proposition}
\begin{proof}
As a consequence of Proposition~\ref{prop:embedding}(i), we have  $\boldsymbol{y}_n \wk \boldsymbol{y}$ in $W^{1,N-1}(\Omega;\RN)$. Thus
\begin{equation*}
\text{$\mathbf{M}_i(D\boldsymbol{y}_n)\wk \mathbf{M}_i(D\boldsymbol{y})$ in $L^1(\Omega;\R^{\binom{N}{i}\times \binom{N}{i}})$ \quad for $i=1,\dots,N-2$}	
\end{equation*}
 by  \cite[Theorem~8.20, Part~4]{Daco}, while \cite[Theorem~8.20, Part~5]{Daco} entails $\cof D \boldsymbol{y}_n \wks \cof D \boldsymbol{y}$ as $\RNN$-valued distributions over $\Omega$.  We deduce $\boldsymbol{\Theta}=\cof D \boldsymbol{y}$. By the Sobolev embedding, up to subsequences, $\boldsymbol{y}_n \to \boldsymbol{y}$ almost everywhere in $\Omega$. Therefore, from \cite[Theorem~4]{mueller}, we deduce that $\boldsymbol{y}$ satisfies \eqref{DIV} and  $\vartheta=\det D \boldsymbol{y}$. 
\end{proof}

We define the class of admissible deformations:
\begin{equation}\label{eqn:admissible-deformation}
	\mathcal{Y}\coloneqq \left \{ \boldsymbol{y}\in W^{1,A}(\Omega;\RN): \hspace{2pt}\det D \boldsymbol{y}\in L^1_+(\Omega), \hspace{4pt} \text{$\boldsymbol{y}$ satisfies (\ref{DIV})}  \right  \}.
\end{equation}

For technical reasons, as in \cite[Definition~2.10]{HS}, we introduce the class of good domains.

\begin{definition}[Good domains]
\label{def:good}
Let $\boldsymbol{y}\in \mathcal{Y}$. We define the class $\mathcal{U}_{\boldsymbol{y}}$ of {good domains} to be formed by all domains $U\subset \subset \Omega$ of class $C^2$ satisfying the following conditions:
\begin{enumerate}[label=(\roman*)]
	\item $\boldsymbol{y}\restr{\partial U}\in W^{1,A}(\partial U;\RN)$ and $(\cof D \boldsymbol{y})\restr{\partial U}\in L^1(\partial U;\RNN)$;
	\item $\haus(\partial U \setminus \domg(\boldsymbol{y},\partial U))=0$ and there holds
	\begin{equation*}
		\text{$D^{\partial U} \boldsymbol{y}\restr{\partial U}=  (D \boldsymbol{y})\restr{\partial U}(\boldsymbol{I}-\boldsymbol{\nu}_U \otimes \boldsymbol{\nu}_U)$ $\haus$-a.e. on $\partial U$;}
	\end{equation*}
	\item $\partial U \subset L_{\boldsymbol{y}}$ and $\boldsymbol{y}^*\restr{\partial U}\in C^0(\partial U;\RN)$;
	\item There holds
	\begin{equation*}
		\displaystyle \lim_{\delta \to 0^+} \dashint_{0}^\delta \left |\int_{\partial U_l} |\cof D \boldsymbol{y}|\,\d\haus- \int_{\partial U} |\cof D \boldsymbol{y}|\,\d\haus\right |\,\d l=0;
	\end{equation*}
	\item For every $\boldsymbol{\psi}\in C^1_{\rm c}(\R^N;\R^N)$, there holds
	\begin{equation*}
		\begin{split}
			\lim_{\delta \to 0^+} \dashint_{0}^\delta \bigg | \int_{\partial U_l} \boldsymbol{\psi}\circ \boldsymbol{y}\cdot \left ((\cof D\boldsymbol{y})\boldsymbol{\nu}_{U_l}\right)\,\d\haus
			- \int_{\partial U} \boldsymbol{\psi}\circ \boldsymbol{y}\cdot \left ( (\cof D\boldsymbol{y})\boldsymbol{\nu}_{U}\right)\,\d\haus          \bigg |\,\d l=0.
		\end{split}
	\end{equation*}
\end{enumerate}
\end{definition}
Definition \ref{def:good} coincides with \cite[Definition~2.10]{HS} except that we included condition (iii). Note that $\boldsymbol{y}\in W^{1,A}(\partial U;\RN)$ and $\partial U \subset L_{\boldsymbol{y}}$ automatically yield $\boldsymbol{y}^*\in C^0(\partial U;\RN)$. Property (iii) is added only for notational simplicity as it unambiguously fixes the continuous representative of $\boldsymbol{y}\restr{\partial U}$ to be $\boldsymbol{y}^*\restr{\partial U}$. As a consequence of \cite[Proposition~4.5]{HS}, the set $NC$ introduced in \cite[Definition~4.4]{HS} is contained in $\Omega \setminus L_{\boldsymbol{y}}$, so that $\partial U \cap NC=\emptyset$ for all $U\in \mathcal{U}_{\boldsymbol{y}}$. Therefore, $\mathcal{U}_{\boldsymbol{y}}$ is contained in the class  $\mathcal{U}_{\boldsymbol{y}}^{NC}$ given in \cite[Definition~4.6]{HS}.

For a comprehensive account on the topological degree, we refer to \cite{FoGa95book}. In the next definition, we take advantage of the known property of the degree to depend solely on boundary values.

\begin{definition}[Topological degree and topological image]
Let $\boldsymbol{y}\in\mathcal{Y}$. 
\begin{enumerate}[label=(\roman*),topsep=0pt]
	\item For all $U\in\mathcal{U}_{\boldsymbol{y}}$, the {topological degree} $\deg(\boldsymbol{y},U,\cdot)\colon \RN \setminus \boldsymbol{y}^*(\partial U)\to \Z$ of $\boldsymbol{y}$ over $U$ is defined as the topological degree of any extension of $\boldsymbol{y}^*\restr{\partial U}$ to $\closure{U}$. Then, the {topological image} of $U$ under $\boldsymbol{y}$ is defined as
	\begin{equation*}
		\imt(\boldsymbol{y},U)\coloneqq \left\{
		\boldsymbol{\xi}\in\RN\setminus \boldsymbol{y}^*(\partial U):\:\deg(\boldsymbol{y},U,\boldsymbol{\xi})\neq 0	\right\}.
	\end{equation*} 
	\item The {topological image} of $\Omega$ under $\boldsymbol{y}$ is defined as
	\begin{equation*}
		\imt(\boldsymbol{y},\Omega)\coloneqq \bigcup_{U\in\mathcal{U}_{\boldsymbol{y}}} \imt(\boldsymbol{y},U).
	\end{equation*}
\end{enumerate}
\end{definition}
In view of the continuity of the degree, $\imt(\boldsymbol{y},U)$ is open and bounded with $\partial \hspace{1pt}\imt(\boldsymbol{y},U)=\boldsymbol{y}^*(\partial U)$. Thus, $\imt(\boldsymbol{y},\Omega)$ is also open. In particular, as argued in \cite[p. 17]{HS}, the set $\imt(\boldsymbol{y},\Omega)$ does not dependent on representatives of $\boldsymbol{y}$.

The next result asserts the coincidence of geometric and topological image for admissible deformations. We refer to \cite[Proposition~4.2(c)]{HS} for the proof. 

\begin{proposition}[Excluding cavitation]
	\label{prop:imt-img}
Let $\boldsymbol{y}\in \mathcal{Y}$. Then, the following hold:
\begin{enumerate}[label=(\roman*),topsep=0pt]
	\item $\img(\boldsymbol{y},U)= \imt(\boldsymbol{y},U)$ up to negligible sets for all $U\in \mathcal{U}_{\boldsymbol{y}}$.
	\item $\img(\boldsymbol{y},\Omega)= \imt(\boldsymbol{y},\Omega)$ up to negligible sets.
\end{enumerate} 
\end{proposition}

We recall some consequences of \eqref{DIV} concerning the local invertibility of admissible deformations. The next definition is well posed in view of Lemma~\ref{lem:injectivity} and Proposition~\ref{prop:imt-img}(i).

\begin{definition}[Good domains of injectivity and local inverse]
Let $\boldsymbol{y}\in \mathcal{Y}$.
\begin{enumerate}[label=(\roman*)]
	\item We define the class $\mathcal{U}_{\boldsymbol{y}}^{\rm inj}$ of {good domains of injectivity} to be formed by all $U\in \mathcal{U}_{\boldsymbol{y}}$ such that $\boldsymbol{y}$ is almost everywhere injective in $U$;
	\item For all $U\in \mathcal{U}_{\boldsymbol{y}}^{\rm inj}$, we define the {local inverse} $\boldsymbol{y}_U^{-1}\colon \imt(\boldsymbol{y},U)\to \RN$ of $\boldsymbol{y}$ on $U$ by setting
	\begin{equation*}
		\boldsymbol{y}_U^{-1}(\boldsymbol{\xi})\coloneqq \begin{cases}
			(\boldsymbol{y}\restr{\domg(\boldsymbol{y},\Omega)})^{-1}(\boldsymbol{\xi}) & \text{if $\boldsymbol{\xi} \in \img(\boldsymbol{y},\Omega) \cap \imt(\boldsymbol{y},\Omega)$,}\\
			\boldsymbol{0} & \text{if $\boldsymbol{\xi}\in \imt(\boldsymbol{y},\Omega)\setminus \img(\boldsymbol{y},\Omega)$.}
		\end{cases}
	\end{equation*}	
\end{enumerate}
\end{definition}

We combine \cite[Proposition~4.9 and Proposition~4.11]{HS} as follows.

\begin{proposition}[Local invertibility and regularity of the inverse]
\label{prop:loc-inv}
Let $\boldsymbol{y}\in \mathcal{Y}$. Then, we have the following:
\begin{enumerate}[label=(\roman*)]
	\item For almost all $\boldsymbol{x}_0\in \Omega$, there exists $r>0$ such that $B(\boldsymbol{x}_0,r)\in \mathcal{U}_{\boldsymbol{y}}^{\rm inj}$. 
	\item For any $U\in \mathcal{U}_{\boldsymbol{y}}^{\rm inj}$, there holds $\boldsymbol{y}_U^{-1}\in W^{1,1}(\imt(\boldsymbol{y},U);\RN)$ with $D\boldsymbol{y}_U^{-1}=(D\boldsymbol{y})^{-1}\circ \boldsymbol{y}_U^{-1}$ almost everywhere in $\imt(\boldsymbol{y},U)$. 
\end{enumerate}
\end{proposition}

The next result has been given in \cite[Proposition 4.12]{HS}. 

\begin{proposition}[Weak continuity of the inverse]\label{prop:weak-continuity-inverse}
Let $(\boldsymbol{y}_n)$ be a sequence in $\mathcal{Y}$ and let $\boldsymbol{y}\in \mathcal{Y}$ be such that $\boldsymbol{y}_n \wks \boldsymbol{y}$ in $W^{1,A}(\Omega;\R^N)$. Then, we have the following:
\begin{enumerate}[label=(\roman*)]
	\item For every $U\in \mathcal{U}_{\boldsymbol{y}}$ and for every compact set $H\subset \imt(\boldsymbol{y},U)$, we have $H\subset \imt(\boldsymbol{y}_n,\Omega)$ for $n\gg 1$ depending on both $H$ and $U$.
	\item For every $U\in \mathcal{U}_{\boldsymbol{y}}$ and for every open set $V \subset \subset \imt(\boldsymbol{y},U)$, we have
	\begin{align*}
		\boldsymbol{y}_{n,U}^{-1} \to  \boldsymbol{y}_U^{-1} &\text{ in $L^1(V;\R^N)$}.
	\end{align*}
	Additionally, if the sequence $(\det D \boldsymbol{y}_{n,U}^{-1}\restr{V})$ in $L^1(V)$ is equi-integrable, there holds
	\begin{equation*}
		\text{$\mathbf{M}(D\boldsymbol{y}_{n,U}^{-1})\wk \mathbf{M}(D\boldsymbol{y}_U^{-1})$ in $L^1(V;\textstyle \prod_{i=1}^{N}\R^{\binom{N}{i}\times \binom{N}{i}})$.}
	\end{equation*} 
	\item If the sequence $(\det D \boldsymbol{y}_n)$ in $L^1(\Omega)$ is equi-integrable, then, up to subsequences,  we have
	\begin{equation*}
		\text{$\chi_{\imt(\boldsymbol{y}_n,\Omega)} \to \chi_{\imt(\boldsymbol{y},\Omega)}$ in $L^1(\RN)$.}
	\end{equation*} 
\end{enumerate}	
\end{proposition}

\subsection{Rate-independent systems}\label{subsec:RIS}
We recall the definition of rate-independent system. For a comprehensive treatment of the subject, we refer to \cite{MieRu}. 

Let $T>0$ be a time horizon. A rate-independent system is determined by a triplet $(\mathcal{Q},\mathcal{E},\mathcal{D})$, where 
$\mathcal{Q}$ is the state space, $\mathcal{E}\colon [0,T]\times \mathcal{Q} \to(-\infty,+\infty]$ is the energy functional, and $\mathcal{D}\colon \mathcal{Q} \times \mathcal{Q} \to [0,+\infty]$ is the dissipation distance. 

In this paper, we will consider rate-independent systems $(\mathcal{Q},\mathcal{E},\mathcal{D})$ satisfying the following assumptions:
\begin{enumerate}[label=(Q),topsep=0pt]
	\item \label{it:Q} \textit{Topology of the state space:} $\mathcal{Q}$ is an Hausdorff topological space such that the restriction of its topology to compact sets is separable and metrizable.
\end{enumerate}
\begin{enumerate}[label=(E\arabic*),topsep=0pt]
	\item \label{it:E1}\textit{Compact sublevels:} For all $t\in[0,T]$ and $M>0$, the set $\{\boldsymbol{q}\in\mathcal{Q}:\: \mathcal{E}(t,\boldsymbol{q})\leq M\}$ is compact.
	\item \label{it:E2} \textit{Energetic control of the power:} 
	\begin{enumerate}[label=(\roman*),leftmargin=22pt,topsep=0pt]
		\item \label{it:E21} For all $\boldsymbol{q}\in\mathcal{Q}$, we have $\mathcal{E}(0,\boldsymbol{q})<+\infty$ if and only if $\mathcal{E}(t,\boldsymbol{q})<+\infty$ for all $t\in [0,T]$.
		\item \label{it:E22} For all $\boldsymbol{q}\in\mathcal{Q}$ with $\mathcal{E}(0,\boldsymbol{q})<+\infty$, the map $t\mapsto \mathcal{E}(t,\boldsymbol{q})$ belongs to $AC([0,T])$.
		\item  \label{it:E23}There exist a set $P\subset (0,T)$ with $\mathscr{L}^1(P)=0$, a function $\lambda\in L^1(0,T)$, and a constant $K>0$ such that, for all $\boldsymbol{q}\in\mathcal{Q}$ with $\mathcal{E}(0,\boldsymbol{q})<+\infty$, the function $t\mapsto\mathcal{E}(t,\boldsymbol{q})$ is differentiable at all $t\in (0,T)\setminus P$ and its derivative $\partial_t \mathcal{E}(t,\boldsymbol{q})$ satisfies
		\begin{equation*}
			|\partial_t \mathcal{E}(t,\boldsymbol{q})|\leq \lambda(t) \left( \mathcal{E}(t,\boldsymbol{q}) + K \right) \quad \text{for all $t\in (0,T)\setminus P$.}
		\end{equation*}
	\end{enumerate}
	\item \label{it:E3} \textit{Modulus of continuity of the power:} For all $M>0$, $\varepsilon>0$, and $t\in (0,T)\setminus P$, there exists a constant $\delta=\delta(M,\varepsilon,t)>0$ such that
	\begin{equation*}
		\left | \frac{\mathcal{E}(t+h,\boldsymbol{q})-\mathcal{E}(t,\boldsymbol{q})}{h} - \partial_t \mathcal{E}(t,\boldsymbol{q}) \right |<\varepsilon \quad \text{for all $\boldsymbol{q}\in \mathcal{Q}$ with $\mathcal{E}(t,\boldsymbol{q})\leq M$ and $h\in (-\delta,\delta)$.}
	\end{equation*}
\end{enumerate}	
\begin{enumerate}[label=(D\arabic*),topsep=0pt]
	\item \label{it:D1} \textit{Triangle inequality:} 
			\begin{equation*}
				\mathcal{D}(\boldsymbol{q},\widehat{\boldsymbol{q}})\leq \mathcal{D}(\boldsymbol{q},\widetilde{\boldsymbol{q}})+\mathcal{D}(\widetilde{\boldsymbol{q}},\widehat{\boldsymbol{q}}) \quad \text{for all $\boldsymbol{q},\widehat{\boldsymbol{q}},\widetilde{\boldsymbol{q}}\in\mathcal{Q}$.} 
			\end{equation*}
	\item \label{it:D2}\textit{Positivity:}
			\begin{equation*}
				\text{$\mathcal{D}(\boldsymbol{q},\widehat{\boldsymbol{q}})=0$ if and only if $\boldsymbol{q}=\widehat{\boldsymbol{q}}$ \quad for all $\boldsymbol{q},\widehat{\boldsymbol{q}}\in\mathcal{Q}$.}
			\end{equation*}
		\item \label{it:D3} \textit{Continuity:} 
		\begin{equation*}
			\text{$\mathcal{D}\restr{\{\boldsymbol{q}\in\mathcal{Q}:\:\mathcal{E}(t,\boldsymbol{q})\leq M  \} \times \{  \boldsymbol{q}\in\mathcal{Q}:\:\mathcal{E}(t,\boldsymbol{q})\leq M  \}}$   is continuous for  all   $t\in [0,T]$ and  $M>0$.}
		\end{equation*}  
\end{enumerate}

We briefly comment these assumptions. By \ref{it:Q}, topological and sequential notions are equivalent on  compact subsets of $\mathcal{Q}$, while separability is a technical assumption. Owing to  \ref{it:E1},  for all $t\in [0,T]$, the  functional $\boldsymbol{q}\mapsto \mathcal{E}(t,\boldsymbol{q})$ is coercive and lower semicontinuous and, in turn, admits minimizers. The properties in \ref{it:E2} provide a control of the variation of $\mathcal{E}$ with respect to time by means of Gronwall's inequality.
  By   \ref{it:D1}--\ref{it:D2}, the function  $\mathcal{D}$ is like a distance except that it may take value $+\infty$  and it is required to be symmetric. Eventually, assumption \ref{it:D3} is clear.

We recall the notion of energetic solution to rate-independent systems. 

\begin{definition}[Energetic solution]\label{def:energetic-solution}
A function $\boldsymbol{q}\colon t\mapsto \boldsymbol{q}_t$ from $[0,T]$ to $\mathcal{Q}$ is termed an {energetic solution} to the rate-indepdendent system $(\mathcal{Q},\mathcal{E},\mathcal{D})$ if  the map $t\mapsto \partial_t \mathcal{E}(t,\boldsymbol{q}_t)$ is well defined in $L^1(0,T)$ and the following two conditions are satisfied:
\begin{enumerate}[label=(\roman*)]
	\item \textit{Global stability:} 
	\begin{equation*}
		\mathcal{E}(t,\boldsymbol{q}_t)\leq \mathcal{E}(t,\widehat{\boldsymbol{q}})+\mathcal{D}(\boldsymbol{q}_t,\widehat{\boldsymbol{q}}) \quad \text{for all $t\in [0,T]$ and $\widehat{\boldsymbol{q}}\in \mathcal{Q}$.}
	\end{equation*}
	\item \textit{Energy-dissipation balance:}
	\begin{equation*}
		\mathcal{E}(t,\boldsymbol{q}_t)+\mathrm{Var}_{\mathcal{D}}(\boldsymbol{q};[0,t])=\mathcal{E}(0,\boldsymbol{q}_0)+ \int_0^t \partial_t \mathcal{E}(\tau,\boldsymbol{q}_\tau)\,\d\tau \quad \text{for all $t\in [0,T]$.}
	\end{equation*}
\end{enumerate} 
\end{definition}

In the previous definition, the variation of $\boldsymbol{q}$ with respect to $\mathcal{D}$ on the interval $[t',t'']\subset [0,T]$ is defined as
\begin{equation}\label{eqn:var}
	\mathrm{Var}_{\mathcal{D}}(\boldsymbol{q};[t',t''])\coloneqq \sup \left\{ \sum_{i=1}^m \mathcal{D}(\boldsymbol{q}_{t_i},\boldsymbol{q}_{t_{i-1}}): \hspace{2pt} m\in \N, \hspace{2pt} t'=t_0<t_1<\dots<t_m=t'' \right\}.
\end{equation}
This corresponds to the energy dissipated along the evolution  $t\mapsto \boldsymbol{q}_t$ in the interval $[t',t'']$.

The following result ensures the existence of energetic solutions for rate-independent processes. More general  results are available in the literature \cite[Theorem~2.1.16]{MieRu}, but the following version is sufficient for our purposes, see \cite[Theorem 2.1.11]{MieRu}.

\begin{theorem}[Existence of energetic solutions]
	\label{thm:MR}
Suppose that the assumptions {\rm \ref{it:Q}}, {\rm \ref{it:E1}--\ref{it:E3}}, and {\rm \ref{it:D1}--\ref{it:D3}} are satisfied. Then, for every $\boldsymbol{q}^0\in \mathcal{Q}$ such that
\begin{equation*}
	\mathcal{E}(0,\boldsymbol{q}^0)\leq \mathcal{E}(0,\widehat{\boldsymbol{q}})+\mathcal{D}(\boldsymbol{q}^0,\widehat{\boldsymbol{q}}) \quad \text{for all $\widehat{\boldsymbol{q}}\in \mathcal{Q}$,}
\end{equation*}
there exists an energetic solution $\boldsymbol{q}\colon t\mapsto \boldsymbol{q}_t$ to the rate-independent system $(\mathcal{Q},\mathcal{E},\mathcal{D})$ which is measurable and satisfies  the initial condition $\boldsymbol{q}_0=\boldsymbol{q}^0$.
\end{theorem}

In the previous theorem, the measurability of  $\boldsymbol{q}$ is understood with respect to the Borel $\sigma$-algebra of  $\mathcal{Q}$.

\section{Compactness result}
\label{sec:compactness}

In this section, we present a compactness result which will be systematically employed in the following.   
For the rest of the paper, otherwise differently stated, we assume that
\begin{center}
	$\Omega \subset \R^N$ is a bounded Lipschitz domain, \quad 
	$A$ is an N-function satisfying  ($\Delta_2$)  and \eqref{eqn:growth-infinity}.
\end{center}

As already mentioned,  the integrability assumption \eqref{eqn:growth-infinity} can be dropped when $N=2$.

\subsection{Regular approximate differentiability}
The notion of regular approximate differentiability  has been introduced in \cite{goffman.ziemer}. We recall below its definition.

\begin{definition}[Regular approximate differentiability]
	\label{def:regular-approximate-differentiability}
	A measurable map $\boldsymbol{v}\colon \Omega \to \RN$  is termed {regularly approximately differentiable}  at $\boldsymbol{x}_0\in \Omega$ if there exists a matrix $\nabla \boldsymbol{v}(\boldsymbol{x}_0)\in \R^{N \times N}$ such that
	\begin{equation}
		\label{eqn:regular-approximate-differentiability-aplim}
		\aplim_{r \to 0^+} \sup_{\boldsymbol{x}\in \partial B(\boldsymbol{x}_0,r)} \frac{|\boldsymbol{v}(\boldsymbol{x})-\boldsymbol{v}(\boldsymbol{x}_0)- \nabla \boldsymbol{v}(\boldsymbol{x}_0)(\boldsymbol{x}-\boldsymbol{x}_0)|}{|\boldsymbol{x}-\boldsymbol{x}_0|}=0.
	\end{equation}
	In this case,  $\nabla \boldsymbol{v}(\boldsymbol{x}_0)$  is termed the regular approximate gradient of $\boldsymbol{v}$ at $\boldsymbol{x}_0$.
\end{definition}

More explicitly, \eqref{eqn:regular-approximate-differentiability-aplim} requires the existence of a set $P\subset (0,\dist(\boldsymbol{x}_0;\partial \Omega))$ with 
\begin{equation*}
	\Theta^1_+(P,0)\coloneqq \lim_{r \to 0^+} \frac{\mathscr{L}^1(P\cap (0,r))}{r}=1
\end{equation*}
such that
\begin{equation*}
	\lim_{\substack{r\to 0\\ r\in P}} \sup_{\boldsymbol{x}\in \partial B(\boldsymbol{x}_0,r)} \frac{|\boldsymbol{v}(\boldsymbol{x})-\boldsymbol{v}(\boldsymbol{x}_0)- \nabla \boldsymbol{v}(\boldsymbol{x}_0)(\boldsymbol{x}-\boldsymbol{x}_0)|}{|\boldsymbol{x}-\boldsymbol{x}_0|}=0.
\end{equation*}

The regular approximate gradient coincides with the approximate gradient whenever it exists.

 It  has been proved in \cite[Theorem 3.4]{goffman.ziemer} that the precise representative of maps in the space $W^{1,p}(\Omega;\R^N)$ with $p>N-1$ is almost everywhere regulary approximately differentiable. In the next theorem, we extend this result to the scale of Orlicz-Sobolev spaces.

\begin{theorem}[Regular approximate differentiability]\label{thm:reg-approx-diff}
	Let  $A$ an N-function satisfying  {\rm ($\Delta_2$)}  and \eqref{eqn:growth-infinity}. Then,  the precise representative  of any $\boldsymbol{v}\in W^{1,A}(\Omega;\R^N)$  is almost everywhere regularly approximately differentiable  with regular approximate gradient given by $D\boldsymbol{v}$. 
\end{theorem}
In the following, we will be interested in the regular approximate differentiabiliy of maps $\boldsymbol{y}\in \mathcal{Y}$  given in \eqref{eqn:admissible-deformation}. By \cite[Proposition~4.5]{HS}, any such  $\boldsymbol{y}$ admits a representative $\widehat{\boldsymbol{y}}$ which is almost everywhere differentiable. From this, the almost everywhere regular approximate differentiability of $\boldsymbol{y}$ easily follows. The result in Theorem~\ref{thm:reg-approx-diff} is more general as it  holds for every map in $W^{1,A}(\Omega;\RN)$. 

\begin{proof}[Proof of Theorem \ref{thm:reg-approx-diff}]
	The proof follows the same strategy in \cite[Theorem 3.4]{goffman.ziemer}.  
	Without loss of generality, we can look at scalar-valued maps. 
	Let $v\in W^{1,A}(\Omega)$. By \cite[Theorem~1.4]{alberico.cianchi},  we have
	\begin{equation*}
		\lim_{r \to 0^+}\dashint_{B(\boldsymbol{x}_0,r)} A_N \left( \frac{|v(\boldsymbol{x})-v(\boldsymbol{x}_0)-Dv(\boldsymbol{x}_0)\cdot(\boldsymbol{x}-\boldsymbol{x}_0)|}{r} \right)\,\d \boldsymbol{x} = 0, \quad \text{for almost all $\boldsymbol{x}_0\in\Omega$,}
	\end{equation*}
	where $A_N$ is an N-function with $A_N \succ A$. In turn, 
	\begin{equation}
		\label{eqn:leb1}
		\lim_{r \to 0^+}\dashint_{B(\boldsymbol{x}_0,r)} A \left( \frac{|v(\boldsymbol{x})-v(\boldsymbol{x}_0)-Dv(\boldsymbol{x}_0)\cdot(\boldsymbol{x}-\boldsymbol{x}_0)|}{r} \right)\,\d \boldsymbol{x} = 0, \quad \text{for almost all $\boldsymbol{x}_0\in\Omega$.}
	\end{equation}
	Additionally
	\begin{equation}
		\label{eqn:leb2}
		\lim_{r \to 0^+}\dashint_{B(\boldsymbol{x}_0,r)} A \left( |Dv(\boldsymbol{x})-Dv(\boldsymbol{x}_0)| \right)\,\d \boldsymbol{x} = 0, \quad \text{for almost all $\boldsymbol{x}_0\in\Omega$}
	\end{equation}
	thanks to  \cite[Lemma 3.1]{alberico.cianchi}. 
	
	Fix any $\boldsymbol{x}_0\in L_{v}$ with $v^*(\boldsymbol{x}_0)=v(\boldsymbol{x}_0)$ satisfying \eqref{eqn:leb1}--\eqref{eqn:leb2}. We prove that 
	\begin{equation}
		\label{eqn:aplimgamma}
		\aplim_{r \to 0^+} h(r)=0,
	\end{equation}
	where 
	\begin{equation*}
		h(r)\coloneqq \sup_{\boldsymbol{z}\in S(\boldsymbol{0},1)} \left | \frac{v^*(\boldsymbol{x}_0+r\boldsymbol{z})-v(\boldsymbol{x}_0)}{r}- Dv(\boldsymbol{x}_0)\cdot \boldsymbol{z} \right |.
	\end{equation*}
	Let $\delta_{\boldsymbol{x}_0}\coloneqq \dist(\boldsymbol{x}_0;\partial \Omega)$. For $r\in (0,\delta_{\boldsymbol{x}_0})$,  
 define $w_r \colon B(\boldsymbol{0},1) \to \R$ by setting
	\begin{equation*}
		w_r(\boldsymbol{z})\coloneqq \frac{v(\boldsymbol{x}_0+r\boldsymbol{z})-v(\boldsymbol{x}_0)}{r}-Dv(\boldsymbol{x}_0)\cdot \boldsymbol{z} \quad \text{for all $\boldsymbol{z}\in B(\boldsymbol{0},1)$.}
	\end{equation*} 
	Then, $w_r \in  W^{1,A}(B(\boldsymbol{0},1))$ with $Dw_r(\boldsymbol{z})=Dv(\boldsymbol{x}_0+r\boldsymbol{z})-Dv(\boldsymbol{x}_0)$  for all $\boldsymbol{z}\in B(\boldsymbol{0},1)$. 
	Define the function $f \colon (0,\delta_{\boldsymbol{x}_0}) \to (0,+\infty)$ as 
	\begin{equation*}
		f(r)\coloneqq \int_{B(\boldsymbol{0},1)} \left\{ A(|w_r|)+A(|Dw_r|) \right\}\,\d\boldsymbol{z}.
	\end{equation*}
	From \eqref{eqn:leb1}--\eqref{eqn:leb2}, we deduce
	\begin{equation}
		\label{eqn:alpha}
		\lim_{r\to 0^+}f(r) =0.
	\end{equation}
	By 	the coarea formula, 
	\begin{equation*}
		f(r)=\int_0^1 g_r(t)\,\d t,
	\end{equation*}
	where $g_r \colon [0,1] \to [0,+\infty)$ is given by
	\begin{equation*}
		g_r(t)\coloneqq \int_{S(\boldsymbol{0},t)} \left\{ A(|w_r|)+A(|Dw_r|)\right\}\,\d\haus.
	\end{equation*}
	Define $X_r\coloneqq \left\{ t\in [1/2,1]: \hspace{3pt} g_r(t)<\sqrt{ f(r)}  \right\}$. We have
	\begin{equation*}
		f(r)\geq \int_{[1/2,1]\setminus X_r} g_r(t)\,\d t\geq \sqrt{f(r)}\,\mathscr{L}^1([1/2,1]\setminus X_r),
	\end{equation*}
	so that
	\begin{equation*}
		\mathscr{L}^1([1/2,1]\setminus X_r)\leq \sqrt{f(r)}.
	\end{equation*}
	
	Now, denote by $R$ the set of radii $r\in (0,\delta_{\boldsymbol{x}_0})$ with $S(\boldsymbol{x}_0,r)\subset L_v$ and $v\restr{S(\boldsymbol{x}_0,r)}\in W^{1,A}(S(\boldsymbol{x}_0,r))$ satisfying
	\begin{equation*}
		  \quad D^{S(\boldsymbol{x}_0,r)}v(\boldsymbol{x})= (\boldsymbol{I}-\boldsymbol{\nu}_{S(\boldsymbol{x}_0,r)}(\boldsymbol{x})\otimes \boldsymbol{\nu}_{S(\boldsymbol{x}_0,r)}(\boldsymbol{x}))Dv(\boldsymbol{x})\quad \text{for $\haus$-almost all  $\boldsymbol{x}\in S(\boldsymbol{x}_0,r)$}.
	\end{equation*}
	Fix $r\in R$ and denote by $Y_r$ the set of $t\in [0,1]$ such that $t r\in R$. Then, $\mathscr{L}^1([0,1]\setminus Y_r)=0$ by Lemma~\ref{lem:exceptional} and Proposition~\ref{prop:boundary}. For all $t\in Y_r$, we have
	\begin{equation*}
		D^{S(\boldsymbol{0},t)}w_r(\boldsymbol{z}) = (\boldsymbol{I}-\boldsymbol{\nu}_{S(\boldsymbol{0},t)}(\boldsymbol{z})\otimes \boldsymbol{\nu}_{S(\boldsymbol{0},t)})(\boldsymbol{z}) D w_r(\boldsymbol{z}) \quad \text{for $\haus$-almost all $\boldsymbol{z}\in S(\boldsymbol{0},t)$.}
	\end{equation*}
	Additionally, $S(\boldsymbol{0},t)\subset L_{w_r}$ and, from $v^*\restr{S(\boldsymbol{x}_0,tr)} \in C^0(S(\boldsymbol{x}_0,tr))$, we deduce
	 $w_r^*\restr{S(\boldsymbol{0},t)}\in C^0(S(\boldsymbol{0},t))$
	 with
	\begin{equation*}
		w^*_r(\boldsymbol{z})=\frac{v^*(x_0+r\boldsymbol{z})-v(\boldsymbol{x}_0)}{r}-Dv(\boldsymbol{x}_0)\cdot \boldsymbol{z} \quad \text{for all $\boldsymbol{z}\in S(\boldsymbol{0},t)$.}
	\end{equation*}
	
	Let $t\in X_r \cap Y_r$. By Theorem~\ref{thm:embedding-appendix}, for $A_{N-1}$ defined as in \eqref{eqn:A_N-1}, we have
	\begin{equation*}
		\begin{split}
			\sup_{S(\boldsymbol{0},t)} |w_r^*|& \leq {C_{\rm em}} A_{N-1}^{-1}\left( \dashint_{S(\boldsymbol{0},t)} A(|D^{S(\boldsymbol{0},t)}w_r|)\,\d\haus + \dashint_{S(\boldsymbol{0},t)} A(|w_r|)\,\d\haus  \right)\\ &\leq {C_{\rm em}} A_{N-1}^{-1}\left( \dashint_{S(\boldsymbol{0},t)} A(|Dw_r|)\,\d\haus + \dashint_{S(\boldsymbol{0},t)} A(|w_r|)\,\d\haus  \right)
		\end{split}
	\end{equation*}
	for some constant $C_{\rm em}=C_{\rm em}(N)>0$.
	As $t\in [1/2,1]$, we obtain
	\begin{equation*}
		\sup_{S(\boldsymbol{0},t)} |w_r^*| \leq   C_{\rm em} A_{N-1}^{-1}\left( c(N)  g_r(t) \right)<    C_{\rm em} A_{N-1}^{-1}\left( c(N)  \sqrt{f(r)} \right),
	\end{equation*}
	where $c(N)\coloneqq 2^{N-1}/\haus(\SN)$.
	Equivalently
	\begin{equation*}
		h(tr)\leq  C_{\rm em} A_{N-1}^{-1} \left( c(N)  \sqrt{f(r)}\right) \quad  \text{for all $r\in R$ and $t\in X_r \cap Y_r$.}
	\end{equation*}
	Note that, in the previous equation, the right-hand side goes to zero, as $r\to 0^+$, because of \eqref{eqn:alpha} and \eqref{eq:aa}.

	At this point, we can conclude as in \cite[Theorem 3.4]{goffman.ziemer}. Choose a sequence $(r_j)$ in $R$ which goes to zero, as $j\to \infty$. Proceeding as before, we choose $X_{r_j}\subset [1/2,1]$ with $\mathscr{L}^1([1/2,1]\setminus X_{r_j})\leq \sqrt{f(r_j)}$ and $Y_{r_j}\subset [0,1]$ with $\mathscr{L}^1([0,1] \setminus Y_{r_j})=0$  for which 
	\begin{equation*}
		h(tr_j)\leq  C_{\rm em}A_{N-1}^{-1} \left( c(N) \sqrt{f(r_j)}\right) \quad \text{for all $t\in X_{r_j} \cap Y_{r_j}$.} 
	\end{equation*}
	Setting $P\coloneqq \bigcup_{j=1}^\infty r_j (X_{r_j} \cap Y_{r_j})$, we check that $\Theta^1_+(P,0)=1$ and the previous estimate yields
	\begin{equation*}
		\lim_{\substack{r \to 0 \\ r \in P}  } h(r)=0.
	\end{equation*} 
	This proves \eqref{eqn:aplimgamma}.
\end{proof}

 \subsection{Compactness result}

The following compactness result extends \cite[Theorem 3.2]{bresciani} to the Orlicz-Sobolev setting. 
Here, we refer to the class of admissible deformations  introduced in \eqref{eqn:admissible-deformation}. 

\begin{theorem}[Compactness] \label{thm:compactness}
Let $(\boldsymbol{y}_n)$ be a sequence in $\mathcal{Y}$ and let  $(\boldsymbol{n}_n)$ be a sequence with
\begin{equation*}
	\text{$\boldsymbol{n}_n\in W^{1,2}(\imt(\boldsymbol{y}_n,\Omega);\S^{N-1})$ for every $n\in \N$.}
\end{equation*}
Suppose that
\begin{equation}\label{eqn:bdd}
	\sup_{n \in \N} \left\{\|\boldsymbol{y}_n\|_{W^{1,A}(\Omega;\R^N)}+\|\Gamma(|\cof D \boldsymbol{y}_n|)\|_{L^1(\Omega)}+  \|\gamma(\det D \boldsymbol{y}_n) \|_{L^1(\Omega)}  + \|D\boldsymbol{n}_n\|_{L^2(\imt(\boldsymbol{y}_n,\Omega);\R^{N \times N})} \right\}<+\infty,
\end{equation}
where $\Gamma,\gamma \colon (0,+\infty)\to [0,+\infty]$  are two Borel functions satisfying
\begin{equation}\label{eqn:eta}
	\lim_{\vartheta \to +\infty}\frac{\Gamma(\vartheta)}{\vartheta}= \lim_{\vartheta \to +\infty} \frac{\gamma(\vartheta)}{\vartheta}=+\infty,\qquad   \lim_{\vartheta \to 0^+} \gamma(\vartheta)=+\infty.
\end{equation}
Then, there exist $\boldsymbol{y}\in \mathcal{Y}$ and $\boldsymbol{n}\in W^{1,2}(\imt(\boldsymbol{y},\Omega);\S^{N-1})$ such that, up to subsequences, we have:
\begin{align}
	\label{eqn:compactness-deformation}
	\text{$\boldsymbol{y}_n \wks \boldsymbol{y}$}& \quad \text{in  $W^{1,A}(\Omega;\R^N)$,}\\
	\label{eqn:compactness-determinant}
	\text{$\mathbf{M}(D \boldsymbol{y}_n) \wk \mathbf{M}(D \boldsymbol{y})$}&\quad \text{in $\textstyle L^1(\Omega;\prod_{i=1}^N\R^{\binom{N}{i}\times \binom{N}{i}})$},&\\
	\label{eqn:compactness-director}
	\text{$\chi_{\imt(\boldsymbol{y}_n,\Omega)}\boldsymbol{n}_n \to \chi_{\imt(\boldsymbol{y},\Omega)} \boldsymbol{n}$}&\quad \text{in $L^2(\R^N;\R^N)$,}&\\
	\label{eqn:compactness-director-gradient}
	\text{$\chi_{\imt(\boldsymbol{y}_n,\Omega)}D\boldsymbol{n}_n \wk  \chi_{\imt(\boldsymbol{y},\Omega)} D\boldsymbol{n}$}& \quad \text{in $L^2(\R^N;\R^{N \times N})$,}\\
	\label{eqn:compactness-composition}
	\text{$\boldsymbol{n}_n \circ \boldsymbol{y}_n \to \boldsymbol{n}\circ \boldsymbol{y}$}& \quad \text{in $L^1(\Omega;\R^N)$.}
\end{align}
\end{theorem}

For the proof of the convergence of compositions, we need a preliminary result which is adapted from \cite[Proposition~2.32]{bresciani}. Its proof is analogous to the one in \cite{bresciani} and, hence, we omit it.

\begin{lemma}[Topological image of neasted balls]\label{lem:nested-balls}
Let $\boldsymbol{y}\in \mathcal{Y}$ and $\boldsymbol{x}_0\in L_{\boldsymbol{y}}$. Suppose that $\boldsymbol{y}^*$ is regularly approximately differentiable at $\boldsymbol{x}_0$ with $\det \nabla \boldsymbol{y}^*(\boldsymbol{x}_0)>0$ and $\boldsymbol{y}^*(\boldsymbol{x}_0)=\boldsymbol{y}(\boldsymbol{x}_0)$. Then, there exists $r_{\boldsymbol{x}_0}>0$ with the following property: for every sequence $(\boldsymbol{y}_n)$ in $\mathcal{Y}$ such that $\boldsymbol{y}_n \wks \boldsymbol{y}$ in $W^{1,A}(\Omega)$, there exist $r,r',r''>0$ with $r_{\boldsymbol{x}_0}<r''<r'<r<\dist(\boldsymbol{x}_0;\partial \Omega)$, depending on $(\boldsymbol{y}_n)$, and a subsequence $(\boldsymbol{y}_{n_k})$, depending on $r,r',r''$, for which we have
\begin{equation}\label{eqn:nested-balls1}
	B(\boldsymbol{x}_0,r),B(\boldsymbol{x}_0,r'),B(\boldsymbol{x}_0,r'')\in \mathcal{U}_{\boldsymbol{y}}^{\rm inj} \cap \bigcap_{k=1}^\infty \mathcal{U}_{\boldsymbol{y}_{n_k}}^{\rm inj},
\end{equation}
\begin{equation}\label{eqn:nested-balls2}
	\text{$\boldsymbol{y}_{n_k}^* \to \boldsymbol{y}^*$ uniformly on $S(\boldsymbol{x}_0,r)\cup S(\boldsymbol{x}_0,r')\cup S(\boldsymbol{x}_0,r'')$,}
\end{equation}
\begin{equation}\label{eqn:nested-balls3}
	\imt(\boldsymbol{y},B(\boldsymbol{x}_0,r''))\subset \subset \imt(\boldsymbol{y},B(\boldsymbol{x}_0,r'))\subset \subset \imt(\boldsymbol{y},B(\boldsymbol{x}_0,r)),
\end{equation}
\begin{equation}\label{eqn:nested-balls4}
	\imt(\boldsymbol{y}_{n_k},B(\boldsymbol{x}_0,r''))\subset \subset \imt(\boldsymbol{y},B(\boldsymbol{x}_0,r'))\subset \subset \imt(\boldsymbol{y}_{n_k},B(\boldsymbol{x}_0,r)) \quad \text{for all $k\in \N$.}
\end{equation}
\end{lemma} 

We now move to the proof of the compactness result.

\begin{proof}[Proof of Theorem \ref{thm:compactness}]
The proof proceeds as the one of \cite[Theorem 3.2]{bresciani} and it is subdivided into three steps. 

\textit{Step 1 (Compactness of deformations).} From \eqref{eqn:bdd}-\eqref{eqn:eta}, thanks to De la Vall\'{e}e Poussin criterion and the uniform boundedness of $(\boldsymbol{n}_n)$,  we deduce the existence of functions
\begin{equation*}
	\boldsymbol{y}\in W^{1,A}(\Omega;\R^N), \quad \boldsymbol{\Theta}\in L^1(\Omega;\RNN), \quad \vartheta\in L^1(\Omega), \quad \boldsymbol{b}\in L^\infty(\R^N;\R^N), \quad \boldsymbol{B}\in L^2(\R^N;\R^{N\times N})
\end{equation*}
such that, for a not relabeled subsequence,
\begin{equation}\label{eqn:def}
	\text{$\boldsymbol{y}_n \wks \boldsymbol{y}$ in $W^{1,A}(\Omega;\R^N)$,} \qquad \text{$\cof D \boldsymbol{y}_n \wk \boldsymbol{\Theta}$ in $L^1(\Omega;\RNN)$,} \qquad \text{$\det D \boldsymbol{y}_n \wk \vartheta$ in $L^1(\Omega)$,}
\end{equation}
and
\begin{equation}\label{eqn:nem}
	\text{$\chi_{\imt(\boldsymbol{y}_n,\Omega)}\boldsymbol{n}_n \wks \boldsymbol{b}$ in $L^\infty(\R^N;\R^{N})$}, \qquad \text{$\chi_{\imt(\boldsymbol{y}_n,\Omega)}D\boldsymbol{n}_n \wk \boldsymbol{B}$ in $L^2(\R^N;\R^{N \times N})$.}
\end{equation}
By means of a standard contradiction argument based on \eqref{eqn:eta}, we see that $\vartheta>0$ almost everywhere in $\Omega$. Then, by  Proposition~\ref{prop:weak-det}, we find 
 $\boldsymbol{\Theta}=\cof D \boldsymbol{y}$ and   $\vartheta=\det D \boldsymbol{y}$. In particular,  $\boldsymbol{y}\in \mathcal{Y}$ and \eqref{eqn:compactness-deformation}--\eqref{eqn:compactness-determinant} are proved. As $(\det D \boldsymbol{y}_n)$ is equi-integrable in view of \eqref{eqn:def} and the Dunford-Pettis theorem,  Proposition \ref{prop:weak-continuity-inverse} yields
 \begin{equation}
 	\label{eqn:imt-convergence}
 	\text{$\chi_{\imt(\boldsymbol{y}_n,\Omega)} \to \chi_{\imt(\boldsymbol{y},\Omega)}$ in $L^1(\RN)$,}
 \end{equation}
 again for a not relabeled subsequence.

\textit{Step 2 (Compactness of nematic directors).}
This step of the proof follows \cite[Proposition 7.1]{BHM17}. 
Let $U\in \mathcal{U}_{\boldsymbol{y}}$ and $V\subset \subset \imt(\boldsymbol{y},U)$ be open. By Proposition \ref{prop:weak-continuity-inverse}(i), up to subsequences, we have $V\subset \subset \imt(\boldsymbol{y}_n,\Omega)$ for all $n\in \N$, while  from \eqref{eqn:bdd} we see that $(D \boldsymbol{n}_n)$ is bounded in $L^2(V;\R^{N\times N})$. As the sequence $(\boldsymbol{n}_n)$ is uniformly bounded, we deduce the existence of  $\boldsymbol{n}\in W^{1,2}(V;\R^N)$ such that $\boldsymbol{n}_n \wk \boldsymbol{n}$ in $W^{1,2}(V;\R^N)$ for a not relabeled subsequence. Thanks to the Sobolev embedding, we actually have $\boldsymbol{n}\in W^{1,2}(V;\S^{N-1})$. Moreover, from \eqref{eqn:nem}, we see that $\boldsymbol{n}=\boldsymbol{b}$ and $D\boldsymbol{n}=\boldsymbol{B}$ in $V$.  

By means of a diagonal argument, we show that $\boldsymbol{n}\in W^{1,2}(\imt(\boldsymbol{y},\Omega);\S^{N-1})$ and we select a not relabeled subsequence for which we have
\begin{equation}\label{eqn:diagonal-nematic}
	\text{$\boldsymbol{n}_n \to \boldsymbol{n}$ a.e. in $V$ \quad and \quad  $\boldsymbol{n}_n \wk \boldsymbol{n}$ in $W^{1,2}(V;\R^N)$ \quad for all open sets $V\subset \subset \imt(\boldsymbol{y},\Omega)$.}
\end{equation}
In this way, we find that $\boldsymbol{b}=\boldsymbol{n}$ and $\boldsymbol{B}=D\boldsymbol{n}$ in $\imt(\boldsymbol{y},\Omega)$.
By arguing as in  \cite[Proposition 7.1]{BHM17}, we show that   $\boldsymbol{b}=\boldsymbol{0}$ and $\boldsymbol{B}=\boldsymbol{O}$ in $\R^N\setminus \imt(\boldsymbol{y},\Omega)$. Therefore,  \eqref{eqn:compactness-director-gradient} is established and \eqref{eqn:compactness-director} follows by combining  \eqref{eqn:imt-convergence}--\eqref{eqn:diagonal-nematic}. 

\textit{Step 3 (Convergence of compositions).} 
As in \cite[Proposition 7.8]{BHM17}, define $\widehat{\gamma}\colon (0,+\infty)\to [0,+\infty]$ as $\widehat{\gamma}(z)\coloneqq z\gamma(1/z)$. By \eqref{eqn:eta}, we have 
\begin{equation*}
	\lim_{z \to +\infty} \frac{\widehat{\gamma}(z)}{z}=+\infty.
\end{equation*}
Let $\boldsymbol{x}_0\in L_{\boldsymbol{y}}$ be such that $\boldsymbol{y}^*$ is regularly approximately differentiable at $\boldsymbol{x}_0$ with $\det \nabla \boldsymbol{y}^*(\boldsymbol{x}_0)>0$ and $\boldsymbol{y}^*(\boldsymbol{x}_0)=\boldsymbol{y}(\boldsymbol{x}_0)$. By Theorem \ref{thm:reg-approx-diff}, almost every point in $\Omega$ satisfies these properties. Applying Lemma \ref{lem:nested-balls}, we select $r_{\boldsymbol{x}_0}<r''<r'<r<\dist(\boldsymbol{x}_0;\partial \Omega)$ and a subsequence $(\boldsymbol{y}_{n_k})$ for which \eqref{eqn:nested-balls1}--\eqref{eqn:nested-balls4} hold.
For simplicity, set $B\coloneqq B(\boldsymbol{x}_0,r)$, $B'\coloneqq B(\boldsymbol{x}_0,r')$, and $B''\coloneqq B(\boldsymbol{x}_0,r'')$. Using Lemma~\ref{lem:injectivity} with $\boldsymbol{y}^{-1}_{n_k,B}$ and Proposition~\ref{prop:imt-img}, we compute
\begin{equation*}
	\begin{split}
		\int_{\imt(\boldsymbol{y},B')} \widehat{\gamma} \left(\det D \boldsymbol{y}_{n_k,B}^{-1}\right)\,\d\boldsymbol{\xi}&\leq 
		\int_{\imt(\boldsymbol{y}_{n_k},B)} \widehat{\gamma} \left(\det D \boldsymbol{y}_{n_k,B}^{-1}\right)\,\d\boldsymbol{\xi}\\
		&= \int_{\imt(\boldsymbol{y}_{n_k},B)} \gamma \left( \frac{1}{\det D \boldsymbol{y}_{n_k,B}^{-1}} \right)\,\det D \boldsymbol{y}_{n_k,B}^{-1} \,\d\boldsymbol{\xi}\\
		&=\int_{\imt(\boldsymbol{y}_{n_k},B)} \gamma \left( \det D \boldsymbol{y}_{n_k} \circ \boldsymbol{y}_{n_k,B}^{-1}  \right)\,\det D \boldsymbol{y}_{n_k,B}^{-1} \,\d\boldsymbol{\xi}\\
		&=\int_B \gamma(\det D \boldsymbol{y}_{n_k})\,\d\boldsymbol{x},
	\end{split}
\end{equation*}
where the right-hand side is uniformly bounded in view of \eqref{eqn:bdd}. Thus, the sequence $(\det D \boldsymbol{y}_{n_k,B}^{-1}\restr{\imt(\boldsymbol{y},B')})$  is equi-integrable by \eqref{eqn:nested-balls4} and the De la Vall\'{e}e Poussin criterion. Hence, from Proposition \ref{prop:weak-continuity-inverse}, we obtain
\begin{equation}\label{eqn:convergence-inverses}
	\text{$\boldsymbol{y}_{n_k,B}^{-1}\to \boldsymbol{y}_B^{-1}$ in $L^1(\imt(\boldsymbol{y},B');\R^N)$,} \qquad \text{$\det D \boldsymbol{y}_{n_k,B}^{-1} \wk \det D \boldsymbol{y}_B^{-1}$ in $L^1(\imt(\boldsymbol{y},B'))$.}
\end{equation} 
Now, let $\boldsymbol{\varphi}\in C^1_{\rm c}(\Omega;\R^N)$ with $K\coloneqq \supp\,\boldsymbol{\varphi}\subset B''$. Using Lemma~\ref{lem:injectivity} with $\boldsymbol{y}^{-1}_{n_k,B}$ and Proposition~\ref{prop:imt-img},  observing that $\boldsymbol{\varphi}\circ \boldsymbol{y}_{n_k,B}^{-1}=\boldsymbol{0}$ on $\img(\boldsymbol{y},B')\setminus \img(\boldsymbol{y}_{n_k},K)=\imt(\boldsymbol{y},B')\setminus \img(\boldsymbol{y}_{n_k},K)$, we write
\begin{equation*}
	\begin{split}
		\int_{B''} \boldsymbol{n}_{n_k}\circ \boldsymbol{y}_{n_k}\cdot \boldsymbol{\varphi}\,\d\boldsymbol{x}&=\int_{\img(\boldsymbol{y}_{n_k},B'')} \boldsymbol{n}_{n_k}\cdot \boldsymbol{\varphi}\circ \boldsymbol{y}_{n_k,B}^{-1}\,\det D \boldsymbol{y}_{n_k,B}^{-1}\,\d\boldsymbol{\xi}\\
		&=\int_{\imt(\boldsymbol{y},B')} \boldsymbol{n}_{n_k}\cdot \boldsymbol{\varphi}\circ \boldsymbol{y}_{n_k,B}^{-1}\,\det D \boldsymbol{y}_{n_k,B}^{-1}\,\d\boldsymbol{\xi}.
	\end{split}
\end{equation*} 
Given \eqref{eqn:diagonal-nematic}--\eqref{eqn:convergence-inverses}, we pass to the limit, as $k\to \infty$, on the right-hand side of the previous equation with the aid of \cite[Proposition 2.61]{FoLe} and we obtain
\begin{equation*}
	\lim_{k\to \infty} \int_{B''} \boldsymbol{n}_{n_k}\circ \boldsymbol{y}_{n_k}\cdot \boldsymbol{\varphi}\,\d\boldsymbol{x}=\int_{\imt(\boldsymbol{y},B')} \boldsymbol{n}\cdot \boldsymbol{\varphi}\circ \boldsymbol{y}_{B}^{-1}\,\det D \boldsymbol{y}_{B}^{-1}\,\d\boldsymbol{\xi}.
\end{equation*}
By Corollary~\ref{cor:cov} with $\boldsymbol{y}^{-1}_B$, the right-hand side of the previous equation equals
\begin{equation*}
	\begin{split}
		\int_{\imt(\boldsymbol{y},B')} \boldsymbol{n}\cdot \boldsymbol{\varphi}\circ \boldsymbol{y}_{B}^{-1}\,\det D \boldsymbol{y}_{B}^{-1}\,\d\boldsymbol{\xi}&=\int_{B'} \boldsymbol{n}\circ \boldsymbol{y}\cdot \boldsymbol{\varphi}\,\d\boldsymbol{x}=\int_{B''} \boldsymbol{n}\circ \boldsymbol{y}\cdot \boldsymbol{\varphi}\,\d\boldsymbol{x}.
	\end{split}
\end{equation*}
Thus, we have shown that $\boldsymbol{n}_{n_k}\circ \boldsymbol{y}_{n_k} \wks \boldsymbol{n}\circ \boldsymbol{y}$ as $\RN$-valued distributions on $B''$. Taking any $1<q< \infty$ and observing that 
\begin{equation*}
	\int_{B''} |\boldsymbol{n}_{n_k}\circ \boldsymbol{y}_{n_k}|^q\,\d\boldsymbol{x}=\leb(B'')=	\int_{B''} |\boldsymbol{n}\circ \boldsymbol{y}|^q\,\d\boldsymbol{x} \qquad \text{for all $k\in \N$,}
\end{equation*}
we deduce that $\boldsymbol{n}_{n_k}\circ \boldsymbol{y}_{n_k}\to \boldsymbol{n}\circ \boldsymbol{y}$ in $L^q(B'';\R^N)$ and, in particular, in $L^q(B(\boldsymbol{x}_0,r_{\boldsymbol{x}_0});\R^N)$. Since almost every point of $\Omega$ satisfies the same properties of $\boldsymbol{x}_0$ by Theorem \ref{thm:reg-approx-diff} and the radius $r_{\boldsymbol{x}_0}$ depends only on $\boldsymbol{y}$, by covering $\Omega$ with balls centred at any such point and applying a diagonal argument, we conclude that $\boldsymbol{n}_{n_k}\circ \boldsymbol{y}_{n_k}\to \boldsymbol{n}\circ \boldsymbol{y}$ in $L^q(\Omega;\R^N)$. We refer to Step 4 in the proof of \cite[Theorem~3.2]{bresciani} for details.
\end{proof}

\section{Quasistatic evolution for time-independent boundary data}
\label{sec:ti}

In this section, we study a quasistatic model for the evolution of nematic elastomers driven by time-dependent applied loads under time-independent boundary conditions.

\subsection{Setting of the problem}\label{subsec:setting-ti}
We begin by describing the quasistatic setting and by listing all the assumptions on the elastic energy density and the applied loads. 

\textbf{\em State space.} Recalling \eqref{eqn:admissible-deformation}, we introduce the class of admissible states 
\begin{equation}\label{eqn:admissible-state}
\mathcal{Q}\coloneqq \left\{ (\boldsymbol{y},\boldsymbol{n}): \hspace{4pt} \boldsymbol{y}\in \mathcal{Y}, \hspace{3pt} \boldsymbol{n}\in W^{1,2}(\imt(\boldsymbol{y},\Omega);\SN)  \right\}.
\end{equation}
Consider the map from $\mathcal{Q}$ to $W^{1,A}(\Omega;\RN)\times L^2(\RN;\RN) \times   L^2(\RN;\RNN)$ given by
\begin{equation}\label{eqn:topology-Q}
	(\boldsymbol{y},\boldsymbol{n})\mapsto \big(\boldsymbol{y}, \chi_{\imt(\boldsymbol{y},\Omega)}\boldsymbol{n}, \chi_{\imt(\boldsymbol{y},\Omega)} D\boldsymbol{n}  \big),
\end{equation}
where the codomain  is equipped with the product topology determined by the weak-* topology of $W^{1,A}(\Omega;\RN)$ and the weak topologies of $L^2(\RN;\RN)$ and    $L^2(\RN;\RNN)$. We endow $\mathcal{Q}$ with the topology that makes the map in \eqref{eqn:topology-Q} a homeomorphism onto its image. Thus, given a sequence $(\boldsymbol{q}_n)$ in $\mathcal{Q}$ and $\boldsymbol{q}\in\mathcal{Q}$ with $\boldsymbol{q}_n=(\boldsymbol{y}_n,\boldsymbol{n}_n)$ and $\boldsymbol{q}=(\boldsymbol{y},\boldsymbol{n})$, we have $\boldsymbol{q}_n \to \boldsymbol{q}$ in $\mathcal{Q}$ whenever
\begin{align}\label{eqn:convergence-Q}
	\begin{split}
	\boldsymbol{y}_n \wks \boldsymbol{y} \quad  &\text{in $W^{1,A}(\Omega;\RN)$,}\\
	\chi_{\imt(\boldsymbol{y}_n,\Omega)}\boldsymbol{n}_n \wk \chi_{\imt(\boldsymbol{y},\Omega)}\boldsymbol{n}  \quad  &\text{in $L^2(\RN;\RN)$,}\\ 
	\chi_{\imt(\boldsymbol{y}_n,\Omega)}D\boldsymbol{n}_n \wk \chi_{\imt(\boldsymbol{y},\Omega)}D\boldsymbol{n}  \quad  &\text{in $L^2(\RN;\RNN)$.}
\end{split}
\end{align}

\textbf{\em Boundary data.} Let $\Lambda \subset \partial \Omega$ be $\haus$-measurable   with $\haus(\Lambda)>0$ and $\boldsymbol{d}\in \mathcal{Y}$. We impose Dirichlet boundary conditions on admissible deformations by considering 
\begin{equation}
	\label{eqn:admissible-state-bc}
	\mathcal{Y}_{\boldsymbol{d}}\coloneqq \left \{\boldsymbol{y}\in \mathcal{Y}:\hspace{2pt}\boldsymbol{y}\restr{\Lambda}=\boldsymbol{d} \restr{\Lambda}  \right  \}, \qquad \mathcal{Q}_{\boldsymbol{d}}\coloneqq \left \{(\boldsymbol{y},\boldsymbol{n})\in \mathcal{Q}: \hspace{2pt} \boldsymbol{y}\in \mathcal{Y}_{\boldsymbol{d}} \right \},
\end{equation}
where we recall \eqref{eqn:admissible-deformation} and \eqref{eqn:admissible-state}. The space $\mathcal{Q}_{\boldsymbol{d}}$ is endowed with the topology induced by $\mathcal{Q}$.

\textbf{\em Nematoelastic energy.} We define the energy functional $I\colon \mathcal{Q}\to [0,+\infty]$ as
\begin{equation}
\label{eqn:energy-I}
I(\boldsymbol{q})\coloneqq  \int_{\Omega} W(D\boldsymbol{y},\boldsymbol{n}\circ \boldsymbol{y})\,\d\boldsymbol{x}+\int_{\imt(\boldsymbol{y},\Omega)} |D\boldsymbol{n}|^2\,\d\boldsymbol{\xi} \quad \text{for all $\boldsymbol{q}=(\boldsymbol{y},\boldsymbol{n})\in\mathcal{Q}$.}
\end{equation}
The  first term on the right-hand side of the previous equation accounts for the elastic energy of the system and involves the density $W\colon \RNN \times \S^{N-1}\to [0,+\infty]$. The second term on the right-hand side  \eqref{eqn:energy-I} represents the nematic term in the so-called one-parameter approximation. 

On the density $W$ we make the following assumptions:
\begin{enumerate}[label*=(W\arabic*),topsep=0pt]
	\item \label{it:W1}\textit{Finiteness and continuity:} The function $W$ is continuous on $\RNN_+$ and $W(\boldsymbol{F})=+\infty$ for all $\boldsymbol{F}\in \RNN$ with $\det \boldsymbol{F}\leq 0$;
	\item \label{it:W2}\textit{Coercivity:} There exist a constant $c_W>0$ and two Borel functions $\Gamma,\gamma\colon (0,+\infty) \to [0,+\infty]$  satisfying \eqref{eqn:eta} such that
	\begin{equation*}
	W(\boldsymbol{F},\boldsymbol{z})\geq c_W A(|\boldsymbol{F}|)+\Gamma(|\adj \boldsymbol{F} |)+\gamma(\det \boldsymbol{F}) \quad \text{for all $\boldsymbol{F}\in \RNN_+$ and $\boldsymbol{z}\in \SN$;}
	\end{equation*}
	\item \label{it:W3}\textit{Polyconvexity:} There exists a function $ \widehat{W}\colon \prod_{i=1}^{N-1} \R^{\binom{N}{i}\times \binom{N}{i}}\times (0,+\infty)\times \SN \to [0,+\infty]$  such that
	\begin{equation*}
	\text{$(\boldsymbol{M}_1,\dots,\boldsymbol{M}_{N-1},\boldsymbol{M}_N)\mapsto \widehat{W}(\boldsymbol{M}_1,\dots,\boldsymbol{M}_{N-1},\boldsymbol{M}_N,\boldsymbol{z})$ is convex for all $\boldsymbol{z}\in \SN$}  
	\end{equation*}
	and
	\begin{equation*}
	W(\boldsymbol{F},\boldsymbol{z})=\widehat{W}(\mathbf{M}(\boldsymbol{F}),\boldsymbol{z}) \quad \text{for all $\boldsymbol{F}\in \RNN_+$ and $\boldsymbol{z}\in \SN$.}
	\end{equation*}
\end{enumerate}

Assumptions \ref{it:W1}--\ref{it:W3} are standard.  Condition \ref{it:W2} requires  superlinear growth in both the  adjugate  and the determinant of the matrix variable of $W$. The growth condition with respect to $\adj \boldsymbol{F}$ determined by $\Gamma$ is unnecessary whenever $\lim_{s \to +\infty} A(s)/s^{N-1}=+\infty$ as for our prototypical example $A(s)=s^{N-1}\log^q(\mathrm{e}+s)$ with $q>N-2$ and $N\geq 3$.

\begin{example}\label{ex:W}
As in \cite{BDesimone,BHM17,HS}, our prototypical  elastic energy density takes the form
\begin{equation}
	\label{eqn:ex-W}
	W(\boldsymbol{F},\boldsymbol{z})\coloneqq \Phi \left( \boldsymbol{N}^{-1}(\boldsymbol{z})\boldsymbol{F}\right), \qquad \boldsymbol{N}(\boldsymbol{z})\coloneqq \mu^{-1}\boldsymbol{z}\otimes \boldsymbol{z}+\mu^{\frac{1}{N-1}} (\boldsymbol{I}-\boldsymbol{z}\otimes \boldsymbol{z}),
\end{equation}
where $\mu>0$ is a material parameter and $\Phi \colon \RNN_+ \to [0,+\infty]$ is an auxiliary density. Specifically, we consider
\begin{equation*}
	\Phi(\boldsymbol{X})\coloneqq  A(|\boldsymbol{X}|)+|\adj \boldsymbol{X} |^\zeta+\sigma(\det \boldsymbol{X}),	
\end{equation*}
where $\zeta>1$, and  $\sigma\colon (0,+\infty) \to [0,+\infty]$ is a convex function satisfying
\begin{equation}
	\label{eqn:gamma}
	\lim_{v \to +\infty} \frac{\sigma(v)}{v}=+\infty, \qquad \lim_{v \to 0^+} \sigma(v)=+\infty.
\end{equation}
The functions $\sigma(v)\coloneqq a v^\alpha- b \log v + c $ or $\sigma(v)\coloneqq a v^\alpha + b  v^{-\beta}$ for suitable  $a,b>0$, $c \in \R$, and $\alpha,\beta>1$ are possible choices.  

Conditions {\rm \ref{it:W1} } and {\rm  \ref{it:W3} } are promptly checked.  To see {\rm  \ref{it:W2}}, observe that
\begin{equation}
	\label{eqn:N}
	\text{$\det \boldsymbol{N}(\boldsymbol{z})=1$, \quad   $|\boldsymbol{N}(\boldsymbol{z})|=\mu_1$, \quad $|\boldsymbol{N}^{-1}(\boldsymbol{z})|=\mu_2$ \quad  for all $\boldsymbol{z}\in \SN$.}
\end{equation} 
where we set
\begin{equation}
	\label{eqn:m}
		\mu_1\coloneqq \sqrt{\mu^{-2}+\mu^{\frac{2}{N-1}}}, \qquad \mu_2\coloneqq \sqrt{\mu^{2}+\mu^{-\frac{2}{N-1}}}.
\end{equation}
Using \eqref{eqn:pp} and \eqref{eqn:N}, we estimate
\begin{equation*}
	A(|\boldsymbol{F}|)\leq  A(\mu_1|\boldsymbol{N}^{-1}(\boldsymbol{z})\boldsymbol{F}|)\leq (\mu_1 +1)^{p_A} A(|\boldsymbol{N}^{-1}(\boldsymbol{z})\boldsymbol{F}|). 
\end{equation*}
Then, from \eqref{eqn:N} and the identity
\begin{equation*}
	\adj \left( \boldsymbol{N}^{-1}(\boldsymbol{z})\boldsymbol{F} \right)=(\adj \boldsymbol{N}^{-1}(\boldsymbol{z}))(\adj \boldsymbol{F})=\boldsymbol{N}(\boldsymbol{z})\,\adj \boldsymbol{F},	
\end{equation*}
we get
\begin{equation*}
	|\adj \boldsymbol{F}|\leq |\boldsymbol{N}^{-1}(\boldsymbol{z})|\,|\adj \left(\boldsymbol{N}^{-1}(\boldsymbol{z})\boldsymbol{F} \right)|=\mu_2|\adj \left(\boldsymbol{N}^{-1}(\boldsymbol{z})\boldsymbol{F} \right)|
\end{equation*}
and
\begin{equation*}
\mu_2^{-\zeta}|\adj \boldsymbol{F}|^\zeta\leq |\adj \left(\boldsymbol{N}^{-1}(\boldsymbol{z})\boldsymbol{F} \right)|^\zeta.
\end{equation*}
Eventually,  \eqref{eqn:N}  yields
\begin{equation*}
\sigma\big (\det(\boldsymbol{N}^{-1}(\boldsymbol{z})\boldsymbol{F}) \big )=\sigma\big( (\det \boldsymbol{N}^{-1}(\boldsymbol{z})) (\det \boldsymbol{F}) \big)=	\sigma(\det \boldsymbol{F}).
\end{equation*}
Altogether, we obtain
\begin{equation*} 
	W(\boldsymbol{F},\boldsymbol{z})\geq (\mu_1+1)^{-p_A} A(|\boldsymbol{F}|)+\,\mu_2^{-\zeta}|\adj \boldsymbol{F}|^\zeta+\sigma(\det \boldsymbol{F}),
\end{equation*}
which is {\rm \ref{it:W2} } for $c_W=(\mu_1+1)^{-p_A}$, $\Gamma(\vartheta)=(\vartheta/\mu_2)^\zeta$, and $\gamma(\vartheta)=\sigma(\vartheta)$.
\end{example}

\textbf{\em Applied loads.} We fix a time horizon $T>0$ and we consider applied loads
\begin{equation}
\label{eqn:applied-loads}
\boldsymbol{f}\in AC([0,T];M^{\bar{A}}(\Omega;\RN)), \quad \boldsymbol{g}\in AC([0,T]; M^{\bar{A}}(\Sigma;\RN)), \quad \boldsymbol{h}\in AC([0,T];L^1(\RN;\RN)),
\end{equation}
representing body forces, surfaces forces, and external  fields, respectively.
In \eqref{eqn:applied-loads}, $\bar{A}$ denotes the conjugate N-function of $A$, and $\Sigma\subset \partial \Omega$ is $\haus$-measurable. 
The work of the applied loads is accounted by the functional $\mathcal{L}\colon [0,T] \times \mathcal{Q} \to \R$ defined as
\begin{equation}
	\label{eqn:functional-L}
\mathcal{L}(t,\boldsymbol{q})\coloneqq\int_{\Omega} \boldsymbol{f}_t\cdot \boldsymbol{y}\,\d\boldsymbol{x}+\int_\Sigma \boldsymbol{g}_t\cdot \boldsymbol{y}\restr{\Sigma}\,\d\haus+\int_{\imt(\boldsymbol{y},\Omega)}\boldsymbol{h}_t\cdot \boldsymbol{n}\,\d\boldsymbol{\xi} \qquad \text{for all $\boldsymbol{q}=(\boldsymbol{y},\boldsymbol{n})\in \mathcal{Q}$.}
\end{equation}
Thus, the total energy $\mathcal{E}\colon [0,T] \times  \mathcal{Q} \to (-\infty,+\infty]$ reads
\begin{equation}
	\label{eqn:functional-E}
	\mathcal{E}(t,\boldsymbol{q})\coloneqq I(\boldsymbol{q})-\mathcal{L}(t,\boldsymbol{q})\quad \text{for all $t\in[0,T]$ and $\boldsymbol{q}\in\mathcal{Q}$}.
\end{equation}

\textbf{\em Dissipation distance.} We define the dissipation distance  $\mathcal{D}\colon \mathcal{Q}\times \mathcal{Q} \to [0,+\infty)$ by setting
\begin{equation}\label{eqn:dissipation}
	\mathcal{D}(\boldsymbol{q},\widehat{\boldsymbol{q}})\coloneqq \int_\Omega |\boldsymbol{n}\circ \boldsymbol{y}-\widehat{\boldsymbol{n}}\circ \widehat{\boldsymbol{y}}|\,\d\boldsymbol{x} \qquad \text{for all $\boldsymbol{q}=(\boldsymbol{y},\boldsymbol{n}),\widehat{\boldsymbol{q}}=(\widehat{\boldsymbol{y}},\widehat{\boldsymbol{n}})\in \mathcal{Q}$.}
\end{equation}
The variation of a function $\boldsymbol{q}\colon [0,T]\to \mathcal{Q}$ with respect to $\mathcal{D}$ is defined as in \eqref{eqn:var}.

\textbf{\em Quasistatic evolution.}  The next theorem constitutes the main result of the section. It asserts the existence of energetic solutions for time-independent boundary data, where we resort to Definition \ref{def:energetic-solution}.

\begin{theorem}[Existence of energetic solutions for time-independent boundary data]
	\label{thm:existence-ti}
	Assume that $W$ satisfies {\rm (W1)--(W3)}. 
	Let $\boldsymbol{d}\in \mathcal{Y}$ and  $\boldsymbol{f}$, $\boldsymbol{g}$, and $\boldsymbol{h}$ be as in \eqref{eqn:applied-loads}. Then,  for every $\boldsymbol{q}^0\in \mathcal{Q}_{\boldsymbol{d}}$ satisfying 
	\begin{equation*}
		\mathcal{E}(0,\boldsymbol{q}^0)\leq \mathcal{E}(0,\widehat{\boldsymbol{q}})+\mathcal{D}(\boldsymbol{q}^0,\widehat{\boldsymbol{q}}) \quad \text{for all $\widehat{\boldsymbol{q}}\in \mathcal{Q}_{\boldsymbol{d}}$,}
	\end{equation*}
	there exists an energetic solution $\boldsymbol{q}\colon t\mapsto \boldsymbol{q}_t$  to the rate-independent system $(\mathcal{Q}_{\boldsymbol{d}},\mathcal{E}\restr{\mathcal{Q}_{\boldsymbol{d}}},\mathcal{D}\restr{\mathcal{Q}_{\boldsymbol{d}} \times \mathcal{Q}_{\boldsymbol{d}}})$  which is measurable and  satisfies the initial condition $\boldsymbol{q}_0=\boldsymbol{q}^0$. 
\end{theorem}

\subsection{Modular inequalities}
In this subsection, we establish two modular inequalities for Orlicz-Sobolev maps which will be instrumental for the analysis in this section.
As these preliminary results  are of independent interest, we state them for N-functions which do not necessarily satisfy \eqref{eqn:growth-infinity}.
We begin with the  Poincar\'{e} inequality.

\begin{proposition}[Poincar\'{e} inequality with trace term in modular form]
	\label{prop:poincare}
	Let  $A$ be an $N$-function satisfying  {\rm ($\Delta_2$)}  and let $\Lambda\subset \partial \Omega$ be $\haus$-measurable  with $\haus(\Lambda)>0$. 
	Then, there exists a constant $C_{\rm P}=C_{\rm P}(\Omega,\Lambda,A)>0$ such that
	\begin{equation*}
		 \int_{\Omega} A(|\boldsymbol{v}|)\,\d \boldsymbol{x} \leq C_{\rm P} \left\{ \int_{\Omega} A(|D\boldsymbol{v}|)\,\d \boldsymbol{x}+\int_{\Lambda} A(|\boldsymbol{v}\restr{\Lambda}|)\,\d \haus \right\} \quad \text{for all $\boldsymbol{v}\in W^{1,A}(\Omega;\RN)$.}
	\end{equation*}
\end{proposition}
\begin{proof}
	We consider scalar-valued maps only. The vector-valued case can be easily deduced from the scalar case thanks to the convexity of $A$. 
	
	First, let us consider $v\in W^{1,A}(\Omega)\cap L^\infty(\Omega)$. Define $w\coloneqq A(|v|)$. As $A$ is locally Lipschitz and $v$ is bounded,  we know that $w\in W^{1,1}(\Omega)$ with  $Dw=A'(|v|)\sgn(v)Dv$ by  \cite[Lemma 8.31]{Adams}. Here, we adopt the convention $\sgn(0)=0$. By applying the classical Poincaré inequality with trace term in $W^{1,1}(\Omega)$ to $w$ (see, e.g., \cite[Theorem B.3.15]{kruzik.roubicek}), we obtain
	\begin{equation}
		\label{eqn:pt}
		\int_{\Omega} A(|v|)\,\d\boldsymbol{x}\leq \widetilde{C}_{\rm P} \left\{ \int_{\Omega} A'(|v|)|D v|\,\d\boldsymbol{x}+\int_{\Lambda} A(|v\restr{\Lambda}|)\,\d\haus \right\}
	\end{equation} 
	for some constant $\widetilde{C}_{\rm P}=\widetilde{C}_{\rm P}(\Omega,\Lambda)>0$. At this point, we argue as in \cite[Proposition 2.13]{MSZ}. 
	Let $0<\varepsilon<1$. Exploiting \eqref{eqn:p}--\eqref{eqn:pp} and the monotonicity of $A$, we estimate
	\begin{equation}
		\label{eqn:pt-gradient}
		\begin{split}
			\int_{\Omega} A'(|v|)|Dv|\,\d\boldsymbol{x}&=\int_{\{|Dv|<\varepsilon |v|\}} A'(|v|)|Dv|\,\d\boldsymbol{x} + \int_{\{|Dv|\geq\varepsilon |v|\}} A'(|v|)|Dv|\,\d\boldsymbol{x} \\
			&\leq \varepsilon \left\{\int_\Omega A'(|v|)|v|\,\d \boldsymbol{x}+\int_\Omega A'\left(\frac{|Dv|}{\varepsilon}\right) \frac{|Dv|}{\varepsilon}\,\d \boldsymbol{x} \right\}\\
			&\leq \varepsilon p_A \left\{ \int_{\Omega} A(|v|)\,\d\boldsymbol{x} + \int_{\Omega} A \left( \frac{|Dv|}{\varepsilon}\right)\,\d \boldsymbol{x} \right\} \\
			&\leq \varepsilon p_A \int_{\Omega} A(|v|)\,\d\boldsymbol{x} + \frac{p_A}{\varepsilon^{p_A-1}}\int_{\Omega} A \left( |Dv|\right)\,\d\boldsymbol{x}.
		\end{split}
	\end{equation}
	Therefore, by combining \eqref{eqn:pt}--\eqref{eqn:pt-gradient} and choosing $\varepsilon \ll 1$,  we obtain
	\begin{equation*}
		\int_\Omega A(|v|)\d\boldsymbol{x}\leq C_{\rm P} \left\{ \int_\Omega A(|Dv|)\,\d\boldsymbol{x}+\int_\Lambda A(|v\restr{\Lambda}|)\,\d\haus  \right\}
	\end{equation*}
	  with a constant ${C}_{\rm P}=C_{\rm P}(\widetilde{C}_{\rm P},p_A)>0$.
	
For $v\in W^{1,A}(\Omega)$ not necessarily bounded, we consider $v_m\coloneqq (-m)\vee (v \wedge m)$ for $m>0$. Then, we have $v_m\in W^{1,A}(\Omega)\cap L^\infty(\Omega)$ with $Dv_m=\chi_{\{ -m\leq v\leq m \}}Dv$ and $v_m\restr{\Lambda}=(-m)\vee (v\restr{\lambda} \wedge m)$. The previous argument yields the estimate
	\begin{equation*}
		\int_\Omega A(|v_m|)\,\d\boldsymbol{x}\leq C_{\rm P} \left\{ A(|Dv_m|)\,\d\boldsymbol{x}+\int_\Lambda A(|v_m\restr{\Lambda}|)\,\d\haus \right\}.
	\end{equation*}
Given the monotonicity of $A$, we can pass to the limit, as $m\to +\infty$, in each of the integrals by means of the monotone convergence theorem, so that 
the desired inequality follows.
\end{proof}

Next, we present the trace inequality.

\begin{proposition}[Trace inequality in modular form]
	\label{prop:trace}
	Let  $A$ be an $N$-function satisfying  {\rm ($\Delta_2$)}. 
	Then, there exists a constant $C_{\rm tr}=C_{\rm tr}(\Omega,A)>0$ such that	
	\begin{equation*}
		 \int_{\partial \Omega} A(|\boldsymbol{v}\restr{\partial \Omega}|)\,\d \haus \leq C_{\rm tr} \left\{ \int_{\Omega} A(|D\boldsymbol{v}|)\,\d\boldsymbol{x}+ \int_{\Omega} A(|\boldsymbol{v}|)\,\d \boldsymbol{x} \right\} \quad \text{for all $\boldsymbol{v}\in W^{1,A}(\Omega;\RN)$.}
	\end{equation*}
\end{proposition}
\begin{proof}
Also here, by convexity and truncation, it is sufficient to consider $v\in W^{1,A}(\Omega)\cap L^\infty(\Omega)$. We define $w\coloneqq A(|v|)$ and we applying the classical trace inequality in $W^{1,1}(\Omega;\RN)$ to $w$. This gives 	
	\begin{equation*}
		\int_{\partial \Omega} A(|v\restr{\partial \Omega}|)\,\d\haus \leq \widetilde{C}_{\rm tr}\left\{ \int_{\Omega} A'(|v|)|D v|\,\d\boldsymbol{x}+\int_{\Omega} A(|v|)\,\d \boldsymbol{x} \right\},
	\end{equation*} 
 for some constant $\widetilde{C}_{\rm tr}=\widetilde{C}_{\rm tr}(\Omega)>0$. Arguing as in  \eqref{eqn:pt-gradient} with $\varepsilon=1$, we find
	\begin{equation*}
		\int_{\Omega} A'(|v|)|Dv|\,\d\boldsymbol{x}\leq p_A \left\{  \int_{\Omega} A(|v|)\,\d\boldsymbol{x} + \int_{\Omega} A \left( |Dv|\right)\,\d\boldsymbol{x} \right\}, 
	\end{equation*}
	so that 
	\begin{equation*}
		\int_{\partial \Omega} A(|v\restr{\partial \Omega}|)\,\d \haus\leq \widetilde{C}_{\rm tr} p_A \int_{\Omega} A(|Dv|)\,\d\boldsymbol{x}+\widetilde{C}_{\rm tr}(p_A+1)\int_{\Omega} A(|v|)\,\d\boldsymbol{x}. 
	\end{equation*}
	Therefore, the desired inequality holds with $C_{\rm tr}\coloneqq \widetilde{C}_{\rm tr}(p_A+1)>0$. 
\end{proof}

\subsection{Proof of Theorem \ref{thm:existence-ti}}

To prove Theorem \ref{thm:existence-ti}, we will resort to Theorem \ref{thm:MR}. Thus, we need to show that  all the properties listed in Subsection \ref{subsec:RIS} are fulfilled. We begin by investigating its coercivity and lower semicontinuity of the energy.

\begin{proposition}[Coercivity and lower semicontinuity of the energy]\label{prop:energy-ti}
Assume that $W$ satisfies {\rm \ref{it:W1}--\ref{it:W3}}. 
Let $\boldsymbol{d}\in \mathcal{Y}$, and  $\boldsymbol{f}$, $\boldsymbol{g}$, and $\boldsymbol{h}$ be as in \eqref{eqn:applied-loads}.	Then, the following hold:	
\begin{enumerate}[label=(\roman*)]
	\item \emph{Coercivity:} There exist two constants $K_1,K_2>0$  such that
	\begin{equation}
		\label{eqn:coercivity-ti}
		\mathcal{E}(t,\boldsymbol{q})\geq K_1 \left( \int_{\Omega} \big \{  A(|D\boldsymbol{y}|) +  \Gamma(|\adj D \boldsymbol{y}|)+  \gamma(\det D \boldsymbol{y})  \big \} \,\d\boldsymbol{x} + \int_{\imt(\boldsymbol{u},\Omega)} |D\boldsymbol{n}|^2\,\d\boldsymbol{\xi} \right)  - K_2
	\end{equation}
	for all $t\in [0,T]$ and  all $\boldsymbol{q}=(\boldsymbol{y},\boldsymbol{n})\in\mathcal{Q}_{\boldsymbol{d}}$.
	\item \emph{Compactness:} Let $t\in [0,T]$ and let  $(\boldsymbol{q}_n)$ be a sequence in $ \mathcal{Q}_{\boldsymbol{d}}$ with $\sup_{n\in\N} \mathcal{E}(t,\boldsymbol{q}_n) <+\infty$. Then,
	there exists a not relabeled subsequence of $(\boldsymbol{q}_n)$  and $\boldsymbol{q}\in\mathcal{Q}_{\boldsymbol{d}}$ such that $\boldsymbol{q}_n \to \boldsymbol{q}$ in $\mathcal{Q}$ and, in addition, \eqref{eqn:compactness-determinant}--\eqref{eqn:compactness-director} and \eqref{eqn:compactness-composition} hold true for $\boldsymbol{q}_n=(\boldsymbol{y}_n,\boldsymbol{n}_n)$ and $\boldsymbol{q}=(\boldsymbol{y},\boldsymbol{n})$. 
	\item \emph{Lower semicontinuity:} Let $t\in [0,T]$ and let  $(\boldsymbol{q}_n)$ be a sequence in $ \mathcal{Q}$ such that  $\boldsymbol{q}_n \to \boldsymbol{q}$ in $\mathcal{Q}$ for some $\boldsymbol{q}\in \mathcal{Q}$. Then, there holds
	\begin{equation*}
		\mathcal{E}(t,\boldsymbol{q})\leq \liminf_{n\to \infty} \mathcal{E}(t,\boldsymbol{q}_n).
	\end{equation*}
\end{enumerate}
\end{proposition}
	Using Proposition \ref{prop:energy-ti}(ii)--(iii), the existence of minimizers of the functional $\boldsymbol{q}\mapsto \mathcal{E}(t,\boldsymbol{q})$ within $\mathcal{Q}_{\boldsymbol{d}}$ is easily established by means of the direct method. Thus, we recover \cite[Theorem 5.1]{HS} by including also applied loads.	
\begin{proof}
\emph{Claim (i).} First, we show \eqref{eqn:coercivity-ti}. 
Let $t\in [0,T]$ and $\boldsymbol{q}=(\boldsymbol{y},\boldsymbol{n})\in \mathcal{Q}_{\boldsymbol{d}}$. Also, let $\varepsilon>0$ be arbitrary. 
Applying Young's inequality \eqref{eqn:young} and Proposition \ref{prop:poincare}, we estimate
\begin{equation*}
\begin{split}
\left | \int_\Omega \boldsymbol{f}_t\cdot \boldsymbol{y}\,\d\boldsymbol{x} \right | &\leq \int_\Omega \frac{|\boldsymbol{f}_t|}{\varepsilon}\,\varepsilon|\boldsymbol{y}|\,\d\boldsymbol{x}\leq \int_\Omega \bar{A} \left( \frac{|\boldsymbol{f}_t|}{\varepsilon} \right) \,\d\boldsymbol{x}+ \int_\Omega A(\varepsilon|\boldsymbol{y}|)\,\d\boldsymbol{x}\\
&\leq C(\boldsymbol{f},\varepsilon) + \varepsilon \int_\Omega A(|\boldsymbol{y}|)\,\d\boldsymbol{x} \leq C(\boldsymbol{d},\boldsymbol{f},\varepsilon)+\varepsilon \int_\Omega A (|D\boldsymbol{y}|)\,\d\boldsymbol{x}.
\end{split}
\end{equation*}
	Similarly, using \eqref{eqn:young} and Proposition \ref{prop:trace}, we obtain
	\begin{equation*}
	\begin{split}
	\left | \int_\Sigma \boldsymbol{g}_t\cdot \boldsymbol{y}\,\d\boldsymbol{x}  \right | &\leq \int_\Sigma  \frac{|\boldsymbol{g}_t|}{\varepsilon}\,\varepsilon|\boldsymbol{y}|\,\d\boldsymbol{x}\leq \int_\Sigma \bar{A} \left( \frac{|\boldsymbol{g}_t|}{\varepsilon} \right) \,\d\boldsymbol{x} + \int_\Sigma  A(\varepsilon|\boldsymbol{y}|)\,\d\boldsymbol{x}\\
	&\leq C(\boldsymbol{g},\varepsilon) + \varepsilon\int_{\partial \Omega} A(|\boldsymbol{y}|)\,\d\boldsymbol{x}
	\leq C(\boldsymbol{g},\varepsilon) +\varepsilon C_{\rm tr} \left( \int_\Omega A (|D\boldsymbol{y}|)\,\d\boldsymbol{x} + \int_\Omega A(|\boldsymbol{y}|)\,\d\boldsymbol{x} \right),
	\end{split}
	\end{equation*}
	then, applying Proposition \ref{prop:poincare}, we get
	\begin{equation*}
	\left | \int_\Sigma \boldsymbol{g}_t\cdot \boldsymbol{y}\,\d\haus  \right | \leq C_1(C_{\rm P},C_{\rm tr},\boldsymbol{d},\boldsymbol{g},\varepsilon) + \varepsilon C_2(C_{\rm P},C_{\rm tr})  \int_\Omega A (|D\boldsymbol{y}|)\,\d\boldsymbol{x}.
	\end{equation*}
	Eventually, by H\"{o}lder's inequality,  we have 
	\begin{equation*}
	\begin{split}
	\left | \int_{\imt(\boldsymbol{y},\Omega)} \boldsymbol{h}_t\cdot \boldsymbol{n} \,\d\boldsymbol{\xi}\right | \leq \int_{\imt(\boldsymbol{y},\Omega)} |\boldsymbol{h}_t|\,\d\boldsymbol{\xi} \leq \|\boldsymbol{h}_t\|_{L^1(\RN;\RN)}\leq C(\boldsymbol{h}).
	\end{split}
	\end{equation*}
	Altogether, we obtain
	\begin{equation}\label{eqn:L-coercivity}
	|\mathcal{L}(t,\boldsymbol{q})|\leq C_1(C_{\rm P},C_{\rm tr},\boldsymbol{d},\boldsymbol{f},\boldsymbol{g},\boldsymbol{h},\varepsilon) +  \varepsilon C_2(C_{\rm P},C_{\rm tr}) \int_\Omega A (|D\boldsymbol{y}|)\,\d\boldsymbol{x}.
	\end{equation}
	At this point, from \ref{it:W2} and \eqref{eqn:L-coercivity}, the estimate \eqref{eqn:coercivity-ti} follows with constants $K_1=K_1(c_W,C_{\rm P},C_{\rm tr},\varepsilon)>0$ and $K_2=K_2(C_{\rm P},C_{\rm tr},\boldsymbol{d},\boldsymbol{f},\boldsymbol{g},\boldsymbol{h},\varepsilon)>0$ by choosing $\varepsilon \ll 1$. 
	
	\emph{Claim (ii).} By \eqref{eqn:coercivity-ti}, $(D\boldsymbol{y}_n)$ is uniformly bounded in $L^A(\Omega;\RNN)$. Using \eqref{eqn:norm-modular} and Proposition \ref{prop:poincare}, we see that $(\boldsymbol{y}_n)$ is bounded in $W^{1,A}(\Omega;\RN)$. Thus \eqref{eqn:bdd} holds, so that Theorem \ref{thm:compactness} ensures the existence of a not relabeled subsequence of $(\boldsymbol{q}_n)$ for which \eqref{eqn:compactness-deformation}--\eqref{eqn:compactness-composition} hold for some $\boldsymbol{q}=(\boldsymbol{y},\boldsymbol{n})\in\mathcal{Q}$. In particular,  $\boldsymbol{q}\in \mathcal{Q}_{\boldsymbol{d}}$ by trace theory. 
	
	\emph{Claim (iii).} Without loss of generality, we may assume that  $\lim_{n\to \infty} \mathcal{E}(t,\boldsymbol{q}_n)$ exists finite. From \eqref{eqn:coercivity-ti} and Theorem~\ref{thm:compactness}, we find that 
	\eqref{eqn:compactness-deformation}--\eqref{eqn:compactness-determinant}, \eqref{eqn:compactness-director-gradient}, and \eqref{eqn:compactness-composition} hold for $\boldsymbol{q}_n=(\boldsymbol{y}_n,\boldsymbol{n}_n)$ and $\boldsymbol{q}=(\boldsymbol{y},\boldsymbol{n})$. 
	By applying \cite[Theorem 5.4]{ball.currie.olver}, we obtain
\begin{equation*}
	I(\boldsymbol{q})\leq \liminf_{n\to \infty} I(\boldsymbol{q}_n).
\end{equation*}	
Given \eqref{eqn:applied-loads} and the weak continuity of the trace operator, thanks to \eqref{eqn:compactness-deformation} and \eqref{eqn:compactness-director}, we get
\begin{equation*}
	\mathcal{L}(t,\boldsymbol{q})=\lim_{n\to \infty} \mathcal{L}(t,\boldsymbol{q}_n).
\end{equation*}	
Thus, the conclusion follows.
\end{proof}

We now prove the estimates for $\partial_t \mathcal{E}$. 

\begin{proposition}[Estimates for the power]
	\label{prop:power-ti}
	Assume that $W$ satisfies {\rm \ref{it:W1}--\ref{it:W3}}. 
	Let $\boldsymbol{d}\in \mathcal{Y}$, and  $\boldsymbol{f}$, $\boldsymbol{g}$, and $\boldsymbol{h}$ be as in \eqref{eqn:applied-loads}. Denote by $P \subset [0,T]$ the complement of the  set of times at which the functions $\boldsymbol{f}$, $\boldsymbol{g}$, and $\boldsymbol{h}$ are all differentiable.  
	Then, the following hold:	
\begin{enumerate}[label=(\roman*)]
	\item \emph{Control of the power:} There exist a  function $\lambda \in L^1(0,T)$ and a constant $K>0$ such that 
	\begin{equation}
		\label{eqn:power-ti}
		|\partial_t \mathcal{E}(t,\boldsymbol{q})|\leq \lambda(t) \left( \mathcal{E}(t,\boldsymbol{q}) + K \right) \quad \text{for all $t\in (0,T)\setminus P$ and all   $\boldsymbol{q}\in\mathcal{Q}_{\boldsymbol{d}}$.}
	\end{equation}
	\item \emph{Modulus of continuity of the power:}
	Let $M>0$ and $t\in (0,T)\setminus P$. Then, for each $\varepsilon>0$,  there exists a constant $\delta=\delta(M,t,\varepsilon)>0$ such that
	\begin{equation*}
		\left  | \frac{\mathcal{E}(t+h,\boldsymbol{q})-\mathcal{E}(t,\boldsymbol{q})}{h}-\partial_t \mathcal{E}(t,\boldsymbol{q}) \right |\leq \varepsilon \quad \text{for all $\boldsymbol{q}\in \mathcal{Q}$ with $\mathcal{E}(t,\boldsymbol{q})\leq M$ and all $h\in (-\delta,\delta)$.}
	\end{equation*}
\end{enumerate}
\end{proposition}
\begin{proof}
	\emph{Claim (i).} Let $t\in (0,T)\setminus P$ and $\boldsymbol{q}=(\boldsymbol{y},\boldsymbol{n})\in\mathcal{Q}_{\boldsymbol{d}}$. In view of \eqref{eqn:applied-loads}, we compute
	\begin{equation*}
		\partial_t \mathcal{E}(t,\boldsymbol{q})=-\int_{\Omega} \dot{\boldsymbol{f}}_t\cdot \boldsymbol{y}\,\d\boldsymbol{x}-\int_\Sigma \dot{\boldsymbol{g}}_t\cdot \boldsymbol{y}\restr{\Sigma}\,\d\haus-\int_{\imt(\boldsymbol{y},\Omega)}\dot{\boldsymbol{h}}_t\cdot \boldsymbol{n}\,\d\boldsymbol{\xi}. 
	\end{equation*}
	We prove \eqref{eqn:power-ti}. Using \eqref{eqn:norm-modular}, \eqref{eqn:hoelder}, and Proposition \ref{prop:poincare}, we estimate
	\begin{equation*}
		\begin{split}
			\left | \int_{\Omega} \dot{\boldsymbol{f}}_t \cdot \boldsymbol{y}\,\d\boldsymbol{x}  \right |&\leq 2 \|\dot{\boldsymbol{f}}_t\|_{L^{\bar{A}}(\Omega;\RN)} \|\boldsymbol{y}\|_{L^A(\Omega;\RN)} \leq 2 \|\dot{\boldsymbol{f}}_t\|_{L^{\bar{A}}(\Omega;\RN)} \left ( \int_{\Omega}A(|\boldsymbol{y}|)\,\d \boldsymbol{x} +1  \right )\\
			&\leq 2 \|\dot{\boldsymbol{f}}_t\|_{L^{\bar{A}}(\Omega;\RN)} \left( C_{\rm P} \int_{\Omega}A(|D\boldsymbol{y}|)\,\d \boldsymbol{x} +  C_{\rm P} \int_{\Lambda} A(|\boldsymbol{d}\restr{\Lambda}|)\,\d\haus +1\right).
		\end{split}
	\end{equation*}
	Similarly, this time by employing  Proposition  \ref{prop:trace}, we get
	\begin{equation*}
		\begin{split}
			\left | \int_{\Sigma} \dot{\boldsymbol{g}}_t \cdot \boldsymbol{y}\,\d\haus  \right |&\leq 2 \|\dot{\boldsymbol{g}}_t\|_{L^{\bar{A}}(\Sigma;\RN)} \|\boldsymbol{y}\|_{L^A(\Sigma;\RN)} \leq 2 \|\dot{\boldsymbol{g}}_t\|_{L^{\bar{A}}(\Sigma;\RN)} \left ( \int_{\Sigma}A(|\boldsymbol{y}|)\,\d \haus +1  \right )\\
			&\leq 2 \|\dot{\boldsymbol{g}}_t\|_{L^{\bar{A}}(\Sigma;\RN)} \left( C_{\rm tr} \int_{\Omega}A(|D\boldsymbol{y}|)\,\d \boldsymbol{x} +  C_{\rm tr} \int_{\Lambda} A(|\boldsymbol{d}\restr{\Lambda}|)\,\d\haus +1\right).
		\end{split}
	\end{equation*}
	Eventually, applying H\"{o}lder  inequality,  we obtain  
	\begin{equation*}
		\begin{split}
			\left | \int_{\imt(\boldsymbol{y},\Omega)} \dot{\boldsymbol{h}}_t \cdot \boldsymbol{n}\,\d\boldsymbol{\xi} \right |\leq \int_{\imt(\boldsymbol{y},\Omega)} |\dot{\boldsymbol{h}}_t| \,\d\boldsymbol{\xi} \leq \|\dot{\boldsymbol{h}}_t\|_{L^1(\RN;\RN)}. 
		\end{split}
	\end{equation*}
	Altogether, setting
	\begin{equation*}
		\widehat{\lambda}(t)\coloneqq C_1(C_{\rm P},C_{\rm tr}) \left(2 \|\dot{\boldsymbol{f}}_t\|_{L^{\bar{A}}(\Omega;\RN)}+2 \|\dot{\boldsymbol{g}}_t\|_{L^{\bar{A}}(\Sigma;\RN)} +   \|\dot{\boldsymbol{h}}_t\|_{L^1(\RN;\RN)} \right),
	\end{equation*}
	we find  
	\begin{equation*}
		\begin{split}
				|\partial_t \mathcal{E}(t,\boldsymbol{q})|&\leq \widehat{\lambda}(t)   \left(  \int_{\Omega}A(|D\boldsymbol{y}|)\,\d \boldsymbol{x}  +  C_2(C_{\rm P}, C_{\rm tr}, \Omega,\boldsymbol{d})  \right)\\
				&\leq \widehat{\lambda}(t) \left( \int_\Omega \left\{ A(|D\boldsymbol{y}|)+\Gamma(|\adj D \boldsymbol{y}|) + \gamma(\det D \boldsymbol{y}) \right\} \,\d\boldsymbol{x}+ \int_{\imt(\boldsymbol{y},\Omega)} |D \boldsymbol{n}|^2\,\d\boldsymbol{\xi} +C_2(C_{\rm P}, C_{\rm tr}, \Omega,\boldsymbol{d}) \right).
		\end{split}
	\end{equation*}
	Thus, by applying \eqref{eqn:coercivity-ti}, we see that \eqref{eqn:power-ti} holds for $\lambda(t) = \widehat{\lambda}(t)/K_1$ and $K= K_2+C_2/K_1$.
	
	\emph{Claim (ii).} Let  $\boldsymbol{q}\in \mathcal{Q}_{\boldsymbol{d}}$ with $\mathcal{E}(t,\boldsymbol{q})\leq M$. By \eqref{eqn:coercivity-ti},  we have
	\begin{equation*}
		\int_{\Omega} A(|D\boldsymbol{y}|)\,\d\boldsymbol{x}\leq C(M,K_1,K_2).
	\end{equation*}
	From this, using  \eqref{eqn:norm-modular} together with Proposition \ref{prop:poincare} and Proposition \ref{prop:trace}, we obtain
	\begin{equation}
		\label{eqn:C1}
		\|\boldsymbol{y}\|_{L^A(\Omega;\RN)}+ \|\boldsymbol{y}\restr{\Sigma}\|_{L^A(\Sigma;\RN)}\leq C(C_{\rm P},C_{\rm tr},M,\boldsymbol{d}).
	\end{equation}
	Let $h \in \R$ be such that $t+h\in (0,T)$. 
	Using \eqref{eqn:hoelder} and \eqref{eqn:C1},  we estimate
	\begin{equation*}
		\begin{split}
			\left | \int_\Omega \left (\frac{\boldsymbol{f}_{t+h}-\boldsymbol{f}_t}{h}-\dot{\boldsymbol{f}}_t \right ) \cdot \boldsymbol{y}\,\d\boldsymbol{x}  \right | 
			&\leq 2C(C_{\rm P},C_{\rm tr},M,\boldsymbol{d})\left \| \frac{\boldsymbol{f}_{t+h}-\boldsymbol{f}_t}{h}-\dot{\boldsymbol{f}}_t \right \|_{L^{\bar{A}}(\Omega;\RN)}.
		\end{split}
	\end{equation*}
	Analogously 
	\begin{equation*}
		\begin{split}
			\left | \int_\Sigma \left (\frac{\boldsymbol{g}_{t+h}-\boldsymbol{g}_t}{h}-\dot{\boldsymbol{g}}_t \right ) \cdot \boldsymbol{y}\restr{\Sigma}\,\d\haus \right |
			&\leq 2C(C_{\rm P},C_{\rm tr},M,\boldsymbol{d})\left \| \frac{\boldsymbol{g}_{t+h}-\boldsymbol{g}_t}{h}-\dot{\boldsymbol{g}}_t \right \|_{L^{\bar{A}}(\Sigma;\RN)}.
		\end{split}
	\end{equation*}
	Similarly, we obtain
	\begin{equation*}
		\begin{split}
			\left | \int_{\imt(\boldsymbol{y},\Omega)} \left (\frac{\boldsymbol{h}_{t+h}-\boldsymbol{h}_t}{h}-\dot{\boldsymbol{h}}_t \right ) \cdot \boldsymbol{n}\,\d\boldsymbol{\xi} \right | 
			&\leq \int_{\imt(\boldsymbol{y},\Omega)} \left | \frac{\boldsymbol{h}_{t+h}-\boldsymbol{h}_t}{h}-\dot{\boldsymbol{h}}_t \right | \,\d\boldsymbol{\xi} 
			\leq \left \| \frac{\boldsymbol{h}_{t+h}-\boldsymbol{h}_t}{h}-\dot{\boldsymbol{h}}_t \right \|_{L^1(\RN;\RN)}.
		\end{split}
	\end{equation*}
	Therefore 
	\begin{equation*}
		\left | \frac{\mathcal{E}(t+h,\boldsymbol{q})-\mathcal{E}(t,\boldsymbol{q})}{h} - \partial_t \mathcal{E}(t,\boldsymbol{q})\right |\leq (4C(C_{\rm P},C_{\rm tr},M,\boldsymbol{d})+1) \omega_t(h)
	\end{equation*}
	where
	\begin{equation*}
		\omega_t(h)\coloneqq \left \| \frac{\boldsymbol{f}_{t+h}-\boldsymbol{f}_t}{h}-\dot{\boldsymbol{f}}_t \right \|_{L^{\bar{A}}(\Omega;\RN)}+ \left \| \frac{\boldsymbol{g}_{t+h}-\boldsymbol{g}_t}{h}-\dot{\boldsymbol{g}}_t \right \|_{L^{\bar{A}}(\Sigma;\RN)}+ \left \| \frac{\boldsymbol{h}_{t+h
			}-\boldsymbol{h}_t}{h}-\dot{\boldsymbol{h}}_t \right \|_{L^1(\RN;\RN)}.
	\end{equation*}
	Since $\omega_t(h)\to 0$, as $h \to 0$, the desired claim follows. 
\end{proof}

At this point, we are ready to prove the main result of the section.

\begin{proof}[Proof of Theorem \ref{thm:existence-ti}]
We consider the rate-independet system $(\mathcal{Q}_{\boldsymbol{d}},\mathcal{E}\restr{\mathcal{Q}_{\boldsymbol{d}}},\mathcal{D}\restr{\mathcal{Q}_{\boldsymbol{d}}\times \mathcal{Q}_{\boldsymbol{d}}})$ and we check that all  the  assumptions stated in Subsection \ref{subsec:RIS} are satisfied.   The space $\mathcal{Q}$ fulfils condition \ref{it:Q}, so that the same holds for $\mathcal{Q}_{\boldsymbol{d}}$. By Proposition \ref{prop:energy-ti}(ii), for $t\in [0,T]$ and $M>0$, the set $\{\boldsymbol{q}\in\mathcal{Q}_{\boldsymbol{d}}:\:\mathcal{E}(t,\boldsymbol{q})\leq M\}$  is sequentially compact and, by claim (i) of the same proposition and Proposition \ref{prop:poincare}, it is contained in a metrizable subset of $\mathcal{Q}_{\boldsymbol{d}}$. Hence, each of these sublevel sets is actually compact and 
 \ref{it:E1} is satisfied. Conditions \ref{it:E21}--\ref{it:E22} in \ref{it:E2} are easily checked thanks to \ref{it:W1} and \eqref{eqn:applied-loads}, while property \ref{it:E23} in 
 \ref{it:E2} and \ref{it:E3} follow from Proposition~\ref{prop:power-ti}.  Assumptions \ref{it:D1}--\ref{it:D2} are immediate. Eventually, \ref{it:D3} holds in view of Proposition~\ref{prop:energy-ti}(ii). 
Therefore,  Theorem~\ref{thm:existence-ti} is established by invoking	Theorem~\ref{thm:MR}.
\end{proof}

\section{Quasistatic evolution for time-dependent boundary data}

\label{sec:td}

In this section, we extend the quasistatic model of the previous section by accounting also for time-dependent boundary conditions. 

\subsection{Setting of the problem}
\label{subsec:setting-td}
We begin by specifying  all the assumptions, some of which are maintained from Section \ref{sec:ti}. 

\textbf{\em State space.}
For technical reasons,   we impose an a confinement condition on admissible deformation. We fix $O\subset\RN$ open and bounded, and we set
\begin{equation*}
	\mathcal{Y}^O\coloneqq \left\{ \boldsymbol{y}\in \mathcal{Y}: \hspace{2pt} \boldsymbol{y}(\boldsymbol{x})\in O \hspace{4pt} \text{for almost all $\boldsymbol{x}\in \Omega$}  \right\}, \qquad \mathcal{Q}^O\coloneqq \left\{ (\boldsymbol{y},\boldsymbol{n})\in \mathcal{Q}:\,\boldsymbol{y}\in \mathcal{Y}^O\right\},
\end{equation*}
where we recall \eqref{eqn:admissible-deformation} and \eqref{eqn:admissible-state}. The class $\mathcal{Q}^O\subset \mathcal{Q}$ is equipped with the subspace topology.
As a consequence of Proposition~\ref{prop:imt-img}, we have $\imt(\boldsymbol{y},\Omega)\subset \closure{O}$ for all $\boldsymbol{y}\in\mathcal{Y}^O$. 

\textbf{\em Boundary data.}
For $T>0$, we consider a time-dependent boundary datum $\boldsymbol{d}\colon t\mapsto \boldsymbol{d}_t$ with regularity
\begin{equation}
	\label{eqn:bd}
	\boldsymbol{d}\in W^{1,\infty}(0,T;\Lip(\closure{O};\closure{O}))\cap \Lip([0,T];\Lip(\closure{O};\closure{O})).
\end{equation} 
For all $t\in [0,T]$, the map  $\boldsymbol{d}_t$ is assumed to be invertible with $\det D \boldsymbol{d}_t>0$ almost everywhere and its inverse $\boldsymbol{d}_t^{-1}$ is assumed to belong to $\Lip(\closure{O};\closure{O})$. We define the function ${\boldsymbol{d}}^{-1}\colon t\mapsto \boldsymbol{d}_t^{-1}$ and we additionally require 
\begin{equation}
	\label{eqn:bd-inverse}
	{\boldsymbol{d}}^{-1}\in W^{1,\infty}(0,T;\Lip(\closure{O};\closure{O})) \cap \Lip([0,T];\Lip(\closure{O};\closure{O})).
\end{equation}
The regularity assumptions on the boundary datum are the same as in \cite{Laz}. In particular, we recall the observation made in \cite[Example 1.2]{Laz}. 

From the identity $D\boldsymbol{d}_t^{-1}=(D\boldsymbol{d}_t)^{-1}\circ \boldsymbol{d}_t$, we see that the map $ t\mapsto (D\boldsymbol{d}_t)^{-1}$ belongs to $ L^{\infty}(0,T;L^\infty({O};\RNN))$. Also, we note that 
\begin{equation}
	\label{eq:md}
	\inf_{t\in [0,T]} \inf_{\mathstrut \boldsymbol{w}\in O} \det D \boldsymbol{d}_t(\boldsymbol{w})\geq \frac{1}{\|\boldsymbol{d}^{-1}\|_{C^0([0,T];\Lip(\closure{O};\closure{O}))}}.
\end{equation}

From \eqref{eqn:bd}, by arguing as in \cite[Remark~1.10]{Laz}, we deduce the existence of a set $Z_{\boldsymbol{d}}\subset O$ with $\leb(Z_{\boldsymbol{d}})=0$ such that
\begin{equation}
	\label{eq:Z}
	\text{$\boldsymbol{d}_t$ is differentiable on $O \setminus Z_{\boldsymbol{d}}$ for all $t\in [0,T]$.}
\end{equation}
For every $t\in [0,T]$, we have $\dot{\boldsymbol{d}}_t\in \Lip(\closure{O};\closure{O})$. Thus, by Rademacher's theore, there exists a set $\widetilde{Z}_{\boldsymbol{d},t}\subset O$ with $\leb(\widetilde{Z}_{\boldsymbol{d},t})=0$ such that
\begin{equation}
	\label{eq:Z-tilde}
	\text{$\dot{\boldsymbol{d}}_t$ is differentiable on $O \setminus \widetilde{Z}_{\boldsymbol{d},t}$ for all $t\in [0,T]$.}
\end{equation} 

Similarly to \eqref{eqn:admissible-state-bc}, given a $\haus$-measurable subset $\Lambda \subset \partial \Omega$ with $\haus(\Lambda)>0$, we employ the notation
\begin{equation*}
	\mathcal{Y}^O_{\boldsymbol{d}_t}\coloneqq  \left\{ \boldsymbol{y}\in \mathcal{Y}^O: \hspace{2pt}\boldsymbol{y}\restr{\Lambda}=\boldsymbol{d}_t\restr{\Lambda} \right\}, \qquad \mathcal{Q}^O_{\boldsymbol{d}_t}\coloneqq \left\{ (\boldsymbol{y},\boldsymbol{n})\in \mathcal{Q}: \hspace{2pt} \boldsymbol{y}\in \mathcal{Y}^O_{\boldsymbol{d}_t} \right\}.
\end{equation*}

\textbf{\em Nematoelastic energy.}
The nematoelastic energy $I\colon \mathcal{Q}\to [0,+\infty]$ is  defined as in \eqref{eqn:energy-I}. In this section, in addition to \ref{it:W1}--\ref{it:W3},  the density  $W\colon \RNN \times \S^{N-1}\to [0,+\infty]$  is required to satisfy: 
\begin{enumerate}[label*=(W\arabic*),topsep=0pt,start=4]
	\item \label{it:W4}\textit{Regularity:}  The map $\boldsymbol{F}\mapsto  W(\boldsymbol{F},\boldsymbol{z})$ is of class $C^1$ on $\RNN_+$ for every fixed $\boldsymbol{z}\in \S^{N-1}$.
	\item \label{it:W5} \textit{Monotonicity:} The map $W$ satisfies \ref{it:W2} and there exist $\vartheta_1,\vartheta_2>0$ with $\vartheta_1\leq \vartheta_2$ such that
	\begin{equation*}
		\label{eqn:W5}
		\text{$\Gamma$ is decreasing on $(0,\vartheta_1)$, \qquad $\Gamma$ and $\gamma$ are both increasing on $(\vartheta_2,+\infty)$.}
	\end{equation*}
\end{enumerate}	
In view of \ref{it:W4}, denoting by $\partial_{\boldsymbol{F}} W(\boldsymbol{F},\boldsymbol{z})$ the gradient of the function $\boldsymbol{F}\mapsto  W(\boldsymbol{F},\boldsymbol{z})$, we  define the Kirchhoff stress tensor $\boldsymbol{K}\colon \RNN_+ \times \S^{N-1} \to \RNN$ by setting $\boldsymbol{K}(\boldsymbol{F},\boldsymbol{z})\coloneqq \partial_{\boldsymbol{F}} W(\boldsymbol{F},\boldsymbol{z})\boldsymbol{F}^\top$.   
On this tensor, we  make the following assumptions:
\begin{enumerate}[label=(W\arabic*), start=6,topsep=0pt]
	\item \label{it:W6} \textit{Kirchhoff stress estimate:} there exist two constants $a_W>0$ and $b_W\geq0$ such that
	\begin{equation*}
		\text{$|\boldsymbol{K}(\boldsymbol{F},\boldsymbol{z})|\leq a_W \left ( W(\boldsymbol{F},\boldsymbol{z}) + b_W \right)$ \quad for all $\boldsymbol{F}\in \RNN_+$ and $\boldsymbol{z}\in \S^{N-1}$. }
	\end{equation*}
	\item \label{it:W7} \textit{Continuity of Kirchhoff stress:} for every $\varepsilon>0$ there exists $\delta_W=\delta_W(\varepsilon)>0$ such that
	\begin{equation*} 
		\text{$
			|\boldsymbol{K}(\boldsymbol{G}\boldsymbol{F},\boldsymbol{z})-\boldsymbol{K}(\boldsymbol{F},\boldsymbol{z})|\leq \varepsilon \left( W(\boldsymbol{F},\boldsymbol{z}) + b_W \right)$} 
	\end{equation*}
	for all $\boldsymbol{F},\boldsymbol{G}\in \RNN_+$ with $|\boldsymbol{G}-\boldsymbol{I}|<\delta_W$ and $\boldsymbol{z}\in \S^{N-1}$.
\end{enumerate} 

While \ref{it:W4} is clear, condition \ref{it:W5} usually follows from the global monotonicity of $\Gamma$ and the convexity of $\gamma$ . The estimates \ref{it:W6}--\ref{it:W7} for the Kirchhoff tensor have been already assumed in \cite{FraMie}. In particular, these  are compatible with the blow up of the elastic energy under extreme compressions.

\begin{example}
Let $W$, $\boldsymbol{N}$, $\mu$, $\zeta$, $\mu_1$, $\mu_2$, and $\sigma$ be as in Example \ref{ex:W}. Additionally, suppose that $A$ and $\sigma$ are both twice continuously differentiable, and there exist $Q,L>0$  such that  
\begin{equation}
	\label{eqn:A-derivatives}
	A''(s)s^2\leq Q (A(s)+1) \quad \text{for all $s>0$}
\end{equation} 
and
\begin{equation}
	\label{eqn:gamma-derivatives}
	\sigma''(v)v^2+\sigma'(v)v\leq L (\sigma(v)+1) \quad \text{for all $v>0$.}
\end{equation}
For example, our prototypical N-function $A(s)\coloneqq s^{N-1} \log^q(\mathrm{e}+s)$ for $q>N-2$ satisfies \eqref{eqn:A-derivatives}, while \eqref{eqn:gamma-derivatives} is fulfilled by $\sigma(v)\coloneqq a v^\alpha-b \log v +c$ and $\sigma(v)\coloneqq a v^\alpha + b v^{-\beta}$ for $a,b>0$, $c \in\R$, and $\alpha,\beta>1$.

We claim that $W$ as defined in \eqref{eqn:ex-W} satisfies all the requirements. Conditions {\rm \ref{it:W3}}--{\rm \ref{it:W4}} are immediate. As shown in Example \ref{ex:W}, $W$ satisfies {\rm \ref{it:W2}}  for  $	c_W=(\mu_1+1)^{-p_A}$, $\Gamma(\vartheta)=(\vartheta/\mu_2)^{\zeta}$, and $\gamma(\vartheta)=\sigma(\vartheta)$, where  we recall \eqref{eqn:m}. In particular, by  \eqref{eqn:gamma} and the convexity of $\sigma$, the function $\gamma$ admits a global minimum  $\bar{\vartheta}>0$, so that condition  {\rm \ref{it:W5} } is satisfied for $\vartheta_1=\vartheta_2=\bar{\vartheta}$.  

We are left to show {\rm \ref{it:W4}} and {\rm \ref{it:W5}}. Let $\boldsymbol{L}\colon \RNN_+ \to \RNN$ be defined as $\boldsymbol{L}(\boldsymbol{X})\coloneqq D\Phi(\boldsymbol{X})\boldsymbol{X}^\top$.  First, with the aid of \eqref{eqn:adj} and \eqref{eqn:det}, for all $\boldsymbol{X}\in\RNN_+$  we compute
\begin{equation}
	\label{eqn:HH}
	\boldsymbol{L}(\boldsymbol{X})=A'(|\boldsymbol{X}|)\frac{\boldsymbol{X}\boldsymbol{X}^\top }{|\boldsymbol{X}|}+\zeta\left( |\adj \boldsymbol{X}|^\zeta\boldsymbol{I} -|\adj \boldsymbol{X}|^{\zeta-2}(\cof \boldsymbol{X})(\adj \boldsymbol{X}) \right)+\sigma'(\det \boldsymbol{X}) (\det \boldsymbol{X})\boldsymbol{I}.
\end{equation} 
Then, using \eqref{eqn:p} and \eqref{eqn:gamma-derivatives}, we estimate
\begin{equation*}
	\begin{split}
		|\boldsymbol{L}(\boldsymbol{X})| &\leq A'(|\boldsymbol{X}|)|\boldsymbol{X}|+(\sqrt{N}+1) \zeta|\adj \boldsymbol{X}|^\zeta+\sqrt{N}\,|\sigma'(\det \boldsymbol{X})|\det \boldsymbol{X}\\
		& \leq  p_A\,A(|\boldsymbol{X}|)+(\sqrt{N}+1) \zeta |\adj \boldsymbol{X}|^\zeta +\sqrt{N}\, L (\sigma(\det \boldsymbol{X})+1)
	\end{split}
\end{equation*} 
which gives
\begin{equation}
	\label{eqn:H}
	|\boldsymbol{L}(\boldsymbol{X})|\leq C_1(N,p_A,\zeta,L) (\Phi(\boldsymbol{X}) + 1) \qquad \text{for all $\boldsymbol{X}\in \RNN_+$.}
\end{equation}
Next, by means of \eqref{eqn:adj}--\eqref{eqn:HH},  for all $\boldsymbol{X}\in \RNN_+$ and $\boldsymbol{Y}\in\RNN$ we find (cf. \cite[p.~264]{dalmaso.lazzaroni}):
\begin{equation*}
	\begin{split}
		\langle \d  \boldsymbol{L}(\boldsymbol{X}),\boldsymbol{Y}\boldsymbol{X}\rangle \coloneqq & \left( \frac{A''(|\boldsymbol{X}|)}{|\boldsymbol{X}|^2}-\frac{A'(|\boldsymbol{X}|)}{|\boldsymbol{X}|^3} \right) \left ((\boldsymbol{X}\boldsymbol{X}^\top) :\boldsymbol{Y}\right ) \boldsymbol{X}\boldsymbol{X}^\top +\frac{A'(|\boldsymbol{X}|)}{|\boldsymbol{X}|} \left( \boldsymbol{X}\boldsymbol{X}^\top\boldsymbol{Y}^\top + \boldsymbol{Y}\boldsymbol{X}\boldsymbol{X}^\top  \right)\\
		&+  \zeta^2 \left( |\adj \boldsymbol{X}|^\zeta \mathrm{tr}\boldsymbol{Y}-|\adj \boldsymbol{X}|^{\zeta-2} \big ((\cof \boldsymbol{X})(\adj \boldsymbol{X})\big  ):\boldsymbol{Y}  \right) \boldsymbol{I}\\
		&-\zeta(\zeta-2)|\adj \boldsymbol{X}|^{\zeta-2} (\mathrm{tr} \boldsymbol{Y}) (\cof \boldsymbol{X})(\adj \boldsymbol{X})\\
		&+\zeta(\zeta-2) |\adj \boldsymbol{X}|^{\zeta-4} \left( \big((\cof \boldsymbol{X})(\adj \boldsymbol{X})\big ):\boldsymbol{Y} \right) (\cof \boldsymbol{X})(\adj \boldsymbol{X})\\
		&-\zeta |\adj \boldsymbol{X}|^{\zeta-2} \left( 2(\mathrm{tr} \boldsymbol{Y})\boldsymbol{I}-\boldsymbol{Y}-\boldsymbol{Y}^\top  \right)(\cof \boldsymbol{X})(\adj \boldsymbol{X})\\
		& +\left( \sigma''(\det \boldsymbol{X})(\det \boldsymbol{X})^2 + \sigma'(\det \boldsymbol{X})(\det \boldsymbol{X})  \right)(\mathrm{tr}\boldsymbol{Y})\boldsymbol{I}.
	\end{split}
\end{equation*}
Then, using \eqref{eqn:p} and \eqref{eqn:A-derivatives}--\eqref{eqn:gamma-derivatives}, we obtain
\begin{equation}
	\label{eqn:LL}
	|	\langle \d  \boldsymbol{L}(\boldsymbol{X}),\boldsymbol{Y}\boldsymbol{X}\rangle|\leq C_2(Q,L,\zeta,N)(\Phi(\boldsymbol{X})+1)|\boldsymbol{Y}| \quad \text{for all $\boldsymbol{X}\in\RNN_+$ and $\boldsymbol{Y}\in\RNN$.}
\end{equation}
Now, recalling \eqref{eqn:ex-W}, for all $\boldsymbol{F}\in \RNN_+$ and $\boldsymbol{z}\in \SN$ we have
\begin{equation}
	\label{eqn:KH}
	\boldsymbol{K}(\boldsymbol{F},\boldsymbol{z})= \boldsymbol{N}^{-1}(\boldsymbol{z})D\Phi(\boldsymbol{N}^{-1}(\boldsymbol{z})\boldsymbol{F})\boldsymbol{F}^\top= \boldsymbol{N}^{-1}(\boldsymbol{z}) \boldsymbol{L}(\boldsymbol{N}^{-1}(\boldsymbol{z})\boldsymbol{F}) \boldsymbol{N}(\boldsymbol{z}).
\end{equation}
Thus, from \eqref{eqn:m}, \eqref{eqn:H}, and \eqref{eqn:KH}, we see that  {\rm \ref{it:W4}} holds true for $a_W=\mu_1\mu_2C_1$ and $b_W=1$. Eventually, let us denote by $\d_{\boldsymbol{F}} \boldsymbol{K}(\boldsymbol{F},\boldsymbol{z})$  the differential of the function $\boldsymbol{F}\mapsto W(\boldsymbol{F},\boldsymbol{z})$. Given \eqref{eqn:KH}, for all $\boldsymbol{F}\in\RNN_+$, $\boldsymbol{G}\in\RNN$, and $\boldsymbol{z}\in\SN$, we have
\begin{equation*}
	\langle \d_{\boldsymbol{F}}\boldsymbol{K}(\boldsymbol{F},\boldsymbol{z}), \boldsymbol{G}\boldsymbol{F} \rangle =\boldsymbol{N}^{-1}(\boldsymbol{z}) \langle \d \boldsymbol{L}(\boldsymbol{N}^{-1}(\boldsymbol{z})\boldsymbol{F}),  \boldsymbol{N}^{-1}(\boldsymbol{z}) \boldsymbol{G}\boldsymbol{F}\rangle \boldsymbol{N}(\boldsymbol{z}).
\end{equation*} 
From the previous equation, by applying \eqref{eqn:LL} with $\boldsymbol{X}=\boldsymbol{N}^{-1}(\boldsymbol{z})\boldsymbol{F}$ and $\boldsymbol{Y}=\boldsymbol{N}^{-1}(\boldsymbol{z})\boldsymbol{G}\boldsymbol{N}(\boldsymbol{z})$, we deduce
\begin{equation*}
	\langle \d_{\boldsymbol{F}} \boldsymbol{K}(\boldsymbol{F},\boldsymbol{z}),\boldsymbol{G}\boldsymbol{F}\rangle \leq C_3\left( W(\boldsymbol{F},\boldsymbol{z}) + 1 \right) |\boldsymbol{G}| \quad \text{for all $\boldsymbol{F}\in \RNN_+$, $\boldsymbol{G}\in \RNN$, and $\boldsymbol{z}\in \SN$,}
\end{equation*}
where $C_3\coloneqq (\mu_1\,\mu_2)^2 C_2$. Therefore, we conclude that $W$ satisfies {\rm \ref{it:W5}} by applying \cite[Proposition~5.2(c)]{FraMie}.
\end{example}

 \textbf{\textit{Applied loads.}}
 Given the boundedness of admissible deformations, we can assume the following regularity
 \begin{equation}
 	\label{eqn:forces}
 	\boldsymbol{f}\in AC([0,T];L^1(\Omega;\RN)), \quad \boldsymbol{g}\in AC([0,T];L^1(\Lambda;\RN))
 \end{equation}
 for the applied forces.
Instead, for the external field we assume
\begin{equation}
	\label{eqn:field}
 	\boldsymbol{h}\in AC([0,T];L^1(O;\RN)) \cap L^1(0,T;W^{1,1}(O;\RN)).
\end{equation}
The functionals $\mathcal{L}\colon [0,T]\times \mathcal{Q}^O\to \R$ and $\mathcal{E}\colon [0,T]\times \mathcal{Q}^O\to (-\infty,+\infty]$ are defined as in \eqref{eqn:functional-L} and \eqref{eqn:functional-E}, respectively.

\textbf{\textit{Dissipation distance.}}
The  dissipation distance  $\mathcal{D}\colon \mathcal{Q}\times \mathcal{Q} \to [0,+\infty)$ and the corresponding variation  are defined as in \eqref{eqn:dissipation} and \eqref{eqn:var}, respectively.

\textbf{\textit{Quasistatic evolution.}} In the setting of Subsection \ref{subsec:setting-td}, it is not possible to directly apply Theorem~\ref{thm:MR} to establish the existence of energetic solutions as in   Definition~\ref{def:energetic-solution}. Indeed, because of the time-dependent boundary conditions, the structural assumptions in Subsection~\ref{subsec:RIS} are not satisfied. 

We follow the approach in \cite{dalmaso.lazzaroni,FraMie,Laz} and we introduce an ad hoc notion of solution.
In order to state our main result, we define $\mathcal{P}\colon [0,T]\times \mathcal{Q}^O \to \R$ by setting
\begin{equation}
	\label{eqn:P}
	\begin{split}
			\mathcal{P}(t,\boldsymbol{q})&\coloneqq \int_\Omega \partial_{\boldsymbol{F}} W(D\boldsymbol{y},\boldsymbol{n\circ \boldsymbol{y}}):D(\dot{\boldsymbol{d}}_t \circ \boldsymbol{d}_t^{-1}\circ \boldsymbol{y})\,\d\boldsymbol{x}\\
			&+ \int_{\imt(\boldsymbol{y},\Omega)}  \left(|D\boldsymbol{n}|^2 \boldsymbol{I}-2 (D\boldsymbol{n})^\top D\boldsymbol{n}  \right)  : \left ( (D\dot{\boldsymbol{d}}_t \circ \boldsymbol{d}_t^{-1})(D\boldsymbol{d}_t^{-1}) \right) \,\d\boldsymbol{\xi}\\
			&+ \int_\Omega \boldsymbol{f}_t \cdot (\dot{\boldsymbol{d}}_t \circ \boldsymbol{d}_t^{-1}\circ \boldsymbol{y})\,\d\boldsymbol{x} 
			+ \int_\Lambda \boldsymbol{g}_t \cdot (\dot{\boldsymbol{d}}_t \circ \boldsymbol{d}_t^{-1}\circ \boldsymbol{y}\restr{\Lambda})\,\d\haus\\
			&+ \int_{\imt(\boldsymbol{y},\Omega)} \mathrm{tr} \left( (D\dot{\boldsymbol{d}}_t \circ \boldsymbol{d}_t^{-1})(D\boldsymbol{d}_t^{-1}) \right) \boldsymbol{h}_t  \cdot \boldsymbol{n}\,\d\boldsymbol{\xi}
	\end{split}
\end{equation}
for all $t\in (0,T)$ and all $\boldsymbol{q}=(\boldsymbol{y},\boldsymbol{n})\in\mathcal{Q}^O$. The functional $\mathcal{P}$ corresponds to the power of the displacement loadings. 

The main result of the section reads as follows.

\begin{theorem}[Existence of energetic solutions for time-dependent boundary data]
	\label{thm:existence-td}
Assume that $W$ satisfies {\rm \ref{it:W1}--\rm \ref{it:W7}}.  Let $\boldsymbol{d}$, $\boldsymbol{f}$, $\boldsymbol{g}$, and $\boldsymbol{h}$ be as in \eqref{eqn:bd}--\eqref{eqn:bd-inverse} and \eqref{eqn:forces}--\eqref{eqn:field}, respectively. 
Then, for every initial datum $\boldsymbol{q}^0\in \mathcal{Q}^O_{\boldsymbol{d}_0}$ satisfying 
\begin{equation}
	\label{eqn:ic-td}
	 \mathcal{E}(0,\boldsymbol{q}^0)\leq \mathcal{E}(0,\widehat{\boldsymbol{q}})+\mathcal{D}(\boldsymbol{q}^0,\widehat{\boldsymbol{q}}) \quad \text{for all $\widehat{\boldsymbol{q}}\in \mathcal{Q}^O_{\boldsymbol{d}_0}$,}
\end{equation}
there exists a  function $\boldsymbol{q}\colon t\mapsto \boldsymbol{q}_t$ from $[0,T]$ to $\mathcal{Q}^O$ which is measurable and satisfies the initial condition $\boldsymbol{q}_0=\boldsymbol{q}^0$ together with the following properties:
\begin{enumerate}[label=(\roman*)]
	\item \emph{Boundary conditions:} $\boldsymbol{q}_t \in \mathcal{Q}^O_{\boldsymbol{d}_t}$ for all $t\in [0,T]$;
	\item \emph{Global stability:} 
	\begin{equation*}
		\mathcal{E}(t,\boldsymbol{q}_t)\leq \mathcal{E}(t,\widehat{\boldsymbol{q}})+\mathcal{D}(\boldsymbol{q}_t,\widehat{\boldsymbol{q}}) \quad \text{for all $t\in [0,T]$ and all $\widehat{\boldsymbol{q}}\in \mathcal{Q}^O_{\boldsymbol{d}_t}$;}
	\end{equation*}
	\item \emph{Energy-dissipation balance:}
	\begin{equation*}
		\mathcal{E}(t,\boldsymbol{q}_t)+\mathrm{Var}_{\mathcal{D}}(\boldsymbol{q};[0,t])=\mathcal{E}(0,\boldsymbol{q}^0)+ \int_0^t \partial_t \mathcal{E}(\tau,\boldsymbol{q}_\tau)\,\d\tau + \int_0^t  \mathcal{P}(\tau,\boldsymbol{q}_\tau)\,\d\tau \quad \text{for all $t\in [0,T]$.}
	\end{equation*}
\end{enumerate}
\end{theorem}

The previous result is achieved by introducing an auxiliary formulation as described in the next subsection. 

\subsection{Auxiliary formulation}

We employ the approach pioneered in \cite{FraMie} and subsequently adopted  in \cite{dalmaso.lazzaroni,Laz} to reduce to a formulation with time-independent boundary data.
For convenience, we set
\begin{equation*}
	\mathcal{Y}^O_{\boldsymbol{id}}\coloneqq \left\{ \boldsymbol{u}\in \mathcal{Y}^O: \hspace{2pt} \boldsymbol{u}\restr{\Lambda}=\boldsymbol{id}\restr{\Lambda} \right\}, \qquad \mathcal{Q}^O_{\boldsymbol{id}}\coloneqq \left\{ (\boldsymbol{u},\boldsymbol{m})\in \mathcal{Q}: \hspace{2pt} \boldsymbol{u}\in \mathcal{Y}^O_{\boldsymbol{id}} \right\}.
\end{equation*}
The space $\mathcal{Q}^O_{\boldsymbol{id}}$ is endowed with the topology induced by $\mathcal{Q}^O$.

Fix $t \in [0,T]$.  We consider deformations $\boldsymbol{y}\in \mathcal{Y}^O$ of the form $\boldsymbol{y}\coloneqq \boldsymbol{d}_t \circ \boldsymbol{u}$, where $\boldsymbol{u}\in \mathcal{Y}^O_{\boldsymbol{id}}$.   In this case,   $\boldsymbol{y}\in \mathcal{Y}^O_{\boldsymbol{d}_t}$  with $D\boldsymbol{y}=(D\boldsymbol{d}_t \circ \boldsymbol{u})D\boldsymbol{u}$  and $\imt(\boldsymbol{y},\Omega)=\boldsymbol{d}_t(\imt(\boldsymbol{u},\Omega))$ as shown in Proposition \ref{prop:equivalence-formulation}(i) below. Then, we take nematic directors of the form $\boldsymbol{n}\coloneqq \boldsymbol{m}\circ \boldsymbol{d}^{-1}_t$, where $\boldsymbol{m}\in W^{1,2}(\imt(\boldsymbol{u},\Omega);\S^{N-1})$. Hence, $\boldsymbol{p}\coloneqq(\boldsymbol{u},\boldsymbol{m})\in \mathcal{Q}^O_{\boldsymbol{id}}$ and $\boldsymbol{q}\coloneqq(\boldsymbol{y},\boldsymbol{n})\in \mathcal{Q}^O_{\boldsymbol{d}_t}$.
By applying Corollary \ref{cor:cov} with $\boldsymbol{d}_t$, we rewrite the energy in \eqref{eqn:energy-I} as
\begin{equation*}
	I(\boldsymbol{q})=\int_\Omega W((D\boldsymbol{d}_t \circ \boldsymbol{u})D\boldsymbol{u},\boldsymbol{m}\circ \boldsymbol{u})\,\d\boldsymbol{x}+\int_{\imt(\boldsymbol{u},\Omega)} |D\boldsymbol{m}(D\boldsymbol{d}_t)^{-1}|^2 \det D \boldsymbol{d}_t\,\d\boldsymbol{w},
\end{equation*}
while, from \eqref{eqn:functional-L}, we obtain
\begin{equation*}
	\mathcal{L}(t,\boldsymbol{q})=\int_{\Omega} \boldsymbol{f}_t \cdot (\boldsymbol{d}_t \circ \boldsymbol{u})\,\d\boldsymbol{x}+\int_\Sigma\boldsymbol{g}_t\cdot (\boldsymbol{d}_t \circ \boldsymbol{u}\restr{\Sigma})\,\d\boldsymbol{x}+\int_{\imt(\boldsymbol{u},\Omega)} (\boldsymbol{h}_t \circ \boldsymbol{d}_t)\det D \boldsymbol{d}_t \cdot \boldsymbol{m}\,\d\boldsymbol{w}.  
\end{equation*}
Therefore, we introduce the functionals $\mathcal{J}\colon [0,T]\times \mathcal{Q}^O\to [0,+\infty]$ and $\mathcal{M}\colon [0,T]\times \mathcal{Q}^O\to \R$ by setting
\begin{equation}
	 \label{eqn:functional-J}
	\mathcal{J}(t,\boldsymbol{p})\coloneqq \int_\Omega W((D\boldsymbol{d}_t \circ \boldsymbol{u})D\boldsymbol{u},\boldsymbol{m}\circ \boldsymbol{u})\,\d\boldsymbol{x}+\int_{\imt(\boldsymbol{u},\Omega)} |D\boldsymbol{m}(D\boldsymbol{d}_t)^{-1}|^2 \det D \boldsymbol{d}_t\,\d\boldsymbol{w},
\end{equation}
and
\begin{equation}
		\label{eqn:functional-K}
	\mathcal{M}(t,\boldsymbol{p})\coloneqq \int_{\Omega} \boldsymbol{f}_t \cdot (\boldsymbol{d}_t \circ \boldsymbol{u})\,\d\boldsymbol{x}+\int_\Sigma\boldsymbol{g}_t\cdot (\boldsymbol{d}_t \circ \boldsymbol{u}\restr{\Sigma})\,\d\boldsymbol{x}+\int_{\imt(\boldsymbol{u},\Omega)} (\boldsymbol{h}_t \circ \boldsymbol{d}_t)\det D \boldsymbol{d}_t \cdot \boldsymbol{m}\,\d\boldsymbol{w}
\end{equation}
for all $t\in [0,T]$ and $\boldsymbol{p}=(\boldsymbol{u},\boldsymbol{m})\in \mathcal{Q}$. 
Then, we define $\mathcal{F}\colon [0,T]\times \mathcal{Q}^O\to (-\infty,+\infty]$ as
\begin{equation*}
	\mathcal{F}(t,\boldsymbol{p})\coloneqq \mathcal{J}(t,\boldsymbol{p})-\mathcal{M}(t,\boldsymbol{p}) \quad \text{for all $t\in [0,T]$ and all $\boldsymbol{p}\in\mathcal{Q}^O$.}
\end{equation*}

The auxiliary formulation lies within the framework presented in Subsection \ref{subsec:RIS}. Recalling Definition~\ref{def:energetic-solution}, the corresponding existence result reads as follows.

\begin{theorem}[Existence of energetic solutions for the auxiliary formulation]
	\label{thm:existence-a}
	Assume that $W$ satisfies {\rm \ref{it:W1}--\ref{it:W7}}.  Let  $\boldsymbol{d}$, $\boldsymbol{f}$, $\boldsymbol{g}$, and $\boldsymbol{h}$ be as in \eqref{eqn:bd}--\eqref{eqn:bd-inverse} and  \eqref{eqn:forces}--\eqref{eqn:field}, respectively. Then, for every initial datum $\boldsymbol{p}^0\in \mathcal{Q}^O_{\boldsymbol{id}}$ satisfying 
	\begin{equation}
		\label{eqn:ic-a}
	 \mathcal{F}(0,\boldsymbol{p}^0)\leq \mathcal{F}(0,\widehat{\boldsymbol{p}})+\mathcal{D}(\boldsymbol{p}^0,\widehat{\boldsymbol{p}}) \quad \text{for all $\widehat{\boldsymbol{p}}\in \mathcal{Q}^O_{\boldsymbol{id}}$,}
	\end{equation}
	there exists an energetic solution $\boldsymbol{p}\colon t\mapsto \boldsymbol{p}_t$ to the rate-independent system $(\mathcal{Q}^O_{\boldsymbol{id}},\mathcal{F}\restr{\mathcal{Q}^O_{\boldsymbol{id}}},\mathcal{D}\restr{\mathcal{Q}^O_{\boldsymbol{id}}\times \mathcal{Q}^O_{\boldsymbol{id}}})$  which is measurable and satisfies the initial condition $\boldsymbol{p}_0=\boldsymbol{p}^0$.
\end{theorem}

The proof of Theorem~\ref{thm:existence-a} will be given in the next subsection. For the moment,  we begin to   formalize the equivalence between the original formulation and the auxiliary one. 
 A  result in the same spirit of the next one was given in \cite[Proposition~5.1]{MoCo}. 

\begin{proposition}[Equivalence of the two formulations]
	\label{prop:equivalence-formulation}
	Let  $\boldsymbol{d}$ and $\boldsymbol{d}^{-1}$ be as in \eqref{eqn:bd}--\eqref{eqn:bd-inverse}. Set 
	\begin{equation}
		\label{eqn:upsilon}
		\Upsilon_t(\boldsymbol{p})\coloneqq (\boldsymbol{d}_t\circ \boldsymbol{u},\boldsymbol{m}\circ \boldsymbol{d}_t^{-1}) \quad \text{for all $t\in [0,T]$ and $\boldsymbol{p}=(\boldsymbol{u},\boldsymbol{m})\in\mathcal{Q}^O$.}
	\end{equation}
	Then, we have the following:
	\begin{enumerate}[label=(\roman*)]
		\item Let $t\in [0,T]$. Then, the map $\Upsilon_t\colon \mathcal{Q}^O\to \mathcal{Q}^O$ is well defined. In particular, for all $\boldsymbol{u}\in \mathcal{Y}^O$, we have 
		\begin{equation}
			\label{eqn:imtU-equivalence}
			\imt(\boldsymbol{d}_t \circ \boldsymbol{u} ,U)= \boldsymbol{d}_t(\imt(\boldsymbol{u},U)) \qquad \text{for all  $U \in \mathcal{U}_{\boldsymbol{u}}\cap \mathcal{U}_{\boldsymbol{d}_t \circ \boldsymbol{u}}$}
		\end{equation}
		and
		\begin{equation}
			\label{eqn:imt-equivalence}
			\imt(\boldsymbol{d}_t \circ \boldsymbol{u},\Omega)=\boldsymbol{d}_t(\imt(\boldsymbol{u},\Omega)).
		\end{equation}
		\item Let $t\in [0,T]$. Then, the map $\Upsilon_t$ is a bijection and its inverse $\Upsilon^{-1}_t\colon    \mathcal{Q}^O  \to  \mathcal{Q}^O$  takes the form
		\begin{equation*}
			\Upsilon_t^{-1}(\boldsymbol{y},\boldsymbol{n})=(\boldsymbol{d}_t^{-1}\circ \boldsymbol{y},\boldsymbol{n}\circ \boldsymbol{d}_t).
		\end{equation*}
		In particular, for every $\boldsymbol{y}\in\mathcal{Y}^O$, we have 
		\begin{equation}
			\label{eqn:imtU-equivalence-inv}
			\imt(\boldsymbol{d}_t^{-1} \circ \boldsymbol{y} ,U)= \boldsymbol{d}_t^{-1}(\imt(\boldsymbol{y},U)) \qquad \text{for all  $U \in \mathcal{U}_{\boldsymbol{y}}\cap \mathcal{U}_{\boldsymbol{d}_t^{-1} \circ \boldsymbol{y}}$}
		\end{equation}
		and
		\begin{equation}
			\label{eqn:imt-equivalence-inv}
			\imt(\boldsymbol{d}_t^{-1} \circ \boldsymbol{y},\Omega)=\boldsymbol{d}_t^{-1}(\imt(\boldsymbol{y},\Omega)).
		\end{equation}
		\item Let $(t_n)$ be  a sequence in $\subset [0,T]$ such that   $t_n \to t$ for some $t\in [0,T]$. Let $(\boldsymbol{p}_n)$ and $(\boldsymbol{q}_n)$ be two sequences in $\mathcal{Q}^O$. We have 
		\begin{equation}
			\label{eqn:seq-cont}
			\text{$\Upsilon_{t_n}(\boldsymbol{p}_n)\to \Upsilon_t(\boldsymbol{p})$ in $\mathcal{Q}^O$ \qquad whenever \qquad $\boldsymbol{p}_n \to \boldsymbol{p}$ in $\mathcal{Q}^O$ for some $\boldsymbol{p}\in\mathcal{Q}^O$  }
		\end{equation}
		and 
		\begin{equation}
			\label{eqn:seq-cont-inv}
			\text{$\Upsilon_{t_n}^{-1}(\boldsymbol{q}_n)\to \Upsilon_t^{-1}(\boldsymbol{q})$ in $\mathcal{Q}^O$ \qquad whenever \qquad  $\boldsymbol{q}_n \to \boldsymbol{q}$ in $\mathcal{Q}^O$ for some $\boldsymbol{q}\in\mathcal{Q}^O$.}
		\end{equation}
		\item Let $t\in [0,T]$. Then, there holds  $\Upsilon_t(\mathcal{Q}_{\boldsymbol{id}}^O)=\mathcal{Q}_{\boldsymbol{d}_t}^O$.
	\end{enumerate}
\end{proposition} 
\begin{proof}
\emph{Claim (i).} Let $\boldsymbol{u}\in \mathcal{Y}^O$. Since $\boldsymbol{u}\in W^{1,1}(\Omega;\RN)$  satisfies Lusin's condition (N${}^{-1}$),    we know that $\boldsymbol{d}_t \circ \boldsymbol{u}\in W^{1,1}(\Omega;\RN)$ with $D(\boldsymbol{\boldsymbol{d}}_t \circ \boldsymbol{u})=(D\boldsymbol{d}_t \circ \boldsymbol{u})D\boldsymbol{u}$ thanks to the chain rule (see, e.g., \cite[Lemma A.11]{bresciani.friedrich.moracorral}). Using \eqref{eqn:pp}, we estimate
	\begin{equation*}
		\begin{split}
			\int_{\Omega} A(|D(\boldsymbol{d}_t \circ \boldsymbol{u})|)\,\d\boldsymbol{x} &\leq \int_\Omega A \left (\|D\boldsymbol{d}_t\|_{L^\infty(O;\R^{N \times N})}\,|D\boldsymbol{u}|\right )\,\d\boldsymbol{x}\\
			&\leq (1+\|D\boldsymbol{d}_t\|_{L^\infty(O;\R^{N \times N})})^{p_A} \int_{\Omega} A(|D\boldsymbol{u}|)\,\d\boldsymbol{x}<+\infty.
		\end{split}
	\end{equation*} 
	As $\boldsymbol{d}_t \circ \boldsymbol{u}\in L^\infty(\Omega;\RN)$, we have $\boldsymbol{d}_t \circ \boldsymbol{u}\in W^{1,A}(\Omega;\RN)$. Also,  $\det D(\boldsymbol{d}_t \circ \boldsymbol{u})=(\det D \boldsymbol{d}_t \circ \boldsymbol{u})\det D \boldsymbol{u}>0$ almost everywhere thanks to the Lusin's condition (N${}^{-1}$) satisfied by $\boldsymbol{u}$. Moreover
	\begin{equation*}
		\|\det D(\boldsymbol{d}_t \circ \boldsymbol{u})\|_{L^1(\Omega)}\leq C(N)\|D\boldsymbol{d}_t\|_{L^\infty(O;\RNN)}^N \|\det D \boldsymbol{u}\|_{L^1(\Omega)}<+\infty.
	\end{equation*}
	We claim that $\boldsymbol{d}_t \circ \boldsymbol{u}$ satisfies \eqref{DIV}.  To see this, 
	let $\boldsymbol{\psi}\in C^1_{\rm c}(\R^N;\R^N)$. The map $\boldsymbol{\zeta}\coloneqq (\adj D \boldsymbol{d}_t)\boldsymbol{\psi}\circ \boldsymbol{d}_t$ belongs to $L^\infty(O;\R^N)$. Thanks Piola's identity and to the Lipschitz continuity of  $\boldsymbol{\psi}\circ \boldsymbol{d}_t$, we find  $\div \boldsymbol{\zeta}=\div \boldsymbol{\psi}\circ \boldsymbol{d}_t\,\det D \boldsymbol{d}_t \in L^\infty(O)$ where the divergence of $\boldsymbol{\zeta}$ is understood in the sense of distributions. Let $(\rho_\varepsilon)_{\varepsilon>0}$ be a sequence of standard mollifiers and set $\boldsymbol{\zeta}_\varepsilon\coloneqq \overline{\boldsymbol{\zeta}}\ast \rho_\varepsilon$, where $\overline{\boldsymbol{\zeta}}\in L^\infty(\RN;\RN)$ denotes the extension of $\boldsymbol{\zeta}$ to the whole space by zero. Then, $\boldsymbol{\zeta}_\varepsilon \in C^1_{\rm c}(\R^N;\R^N)$ satisfies $\div \boldsymbol{\zeta}_\varepsilon=(\div \boldsymbol{\zeta})\ast \rho_\varepsilon$. Clearly, $\boldsymbol{\zeta}_\varepsilon \to \boldsymbol{\zeta}$ and $\div \boldsymbol{\zeta}_\varepsilon \to \div \boldsymbol{\zeta}$ almost everywhere in $O$ and we have the bound
	\begin{equation*}
		\|\boldsymbol{\zeta}_\varepsilon\|_{L^\infty(O;\R^N)} +  \|	\div \boldsymbol{\zeta}_\varepsilon\|_{L^\infty(O)}\leq \|\boldsymbol{\zeta}\|_{L^\infty(O;\R^N)}+\|	\div \boldsymbol{\zeta}\|_{L^\infty(O)}. 
	\end{equation*}
	Testing \eqref{DIV} for $\boldsymbol{u}$ with $\boldsymbol{\zeta}_\varepsilon$ and some arbitrary $\varphi \in C^1_{\rm c}(\Omega)$, we obtain
	\begin{equation*}
		-\int_{\Omega} \boldsymbol{\zeta}_\varepsilon\circ \boldsymbol{u}\cdot (\cof D \boldsymbol{u})D \varphi\,\d\boldsymbol{x}=\int_{\Omega} \div \boldsymbol{\zeta}_\varepsilon \circ \boldsymbol{u}\,(\det D \boldsymbol{u})\,\varphi\,\d\boldsymbol{x}.
	\end{equation*}
	Passing to the limit, as $\varepsilon \to 0^+$, in the previous identity with the aid of the dominated convergence theorem, we deduce
	\begin{equation*}
		-\int_{\Omega} \boldsymbol{\zeta}\circ \boldsymbol{u}\cdot (\cof D \boldsymbol{u})D \varphi\,\d\boldsymbol{x}=\int_{\Omega} \div \boldsymbol{\zeta} \circ \boldsymbol{u}\,(\det D \boldsymbol{u})\,\varphi\,\d\boldsymbol{x}.
	\end{equation*}
	Then, plugging in the expressions of $\boldsymbol{\zeta}$ and $\div\,\boldsymbol{\zeta}$, we recover 
	\begin{equation*}
		-\int_{\Omega} \boldsymbol{\psi}\circ \boldsymbol{d}_t \circ \boldsymbol{u}\cdot (\cof D \boldsymbol{y})D \varphi\,\d\boldsymbol{x}=\int_{\Omega} (\div \boldsymbol{\psi}) \circ \boldsymbol{d}_t \circ \boldsymbol{u} (\det D \boldsymbol{d}_t \circ \boldsymbol{u})\det D \boldsymbol{u}\,\varphi\,\d\boldsymbol{x}.
	\end{equation*}
	As both $\boldsymbol{\psi}$ and $\varphi$ are arbitrary, this proves the claim. Therefore, $\boldsymbol{d}_t \circ \boldsymbol{u}\in \mathcal{Y}^O$.
	
	We move to the proof of \eqref{eqn:imtU-equivalence}. Let $U \in \mathcal{U}_{\boldsymbol{u}}$, so that the map $\boldsymbol{u}^*\restr{\partial U}$ is continuous. Clearly, the set $\boldsymbol{u}^*(\partial U)$ is compact and $\boldsymbol{d}_t\restr{\boldsymbol{u}^*(\partial U)}$ is a homeomorphism. Let $(W_j)_{j\in J}$ be the connected components of $O \setminus \boldsymbol{u}^*(\partial U)$, where $J\subset \N$.
	By Jordan's separation theorem \cite[Theorem 3.29, Claim 4]{FoGa95book},  the connected components of $O \setminus \boldsymbol{d}_t(\boldsymbol{u}^*(\partial U))$ are given by $(V_j)_{j\in J}$, where $V_j\coloneqq \boldsymbol{d}_t(W_j)$ for all $j\in J$. Let $J_0 \subset J$ denote the set of indices for which $\imt(\boldsymbol{u},\Omega)=\bigcup_{j\in J_0} W_j$. Since $\boldsymbol{d}_t\in W^{1,\infty}(O;O)$ is a homeomorphism and $\det D \boldsymbol{d}_t>0$ almost everywhere, we have
	\begin{equation*}
		\deg(\boldsymbol{d}_t,W,\boldsymbol{\xi})=\begin{cases}
			1 & \text{if $\boldsymbol{\xi}\in \boldsymbol{d}_t(W)$,}\\
			0 & \text{if  $\boldsymbol{\xi}\in \RN \setminus \boldsymbol{d}_t(W)$,}
		\end{cases} \qquad \text{for all open sets $W\subset O$.}
	\end{equation*}
	At this point, by applying the multiplication theorem for the degree \cite[Theorem 2.10]{FoGa95book}, we obtain 
	\begin{equation*}
		\deg(\boldsymbol{d}_t \circ \boldsymbol{u} ,U,\boldsymbol{\xi})=\sum_{j\in J} \deg(\boldsymbol{d}_t,W_j,\boldsymbol{\xi})\deg(\boldsymbol{u},U,\boldsymbol{d}_t^{-1}(\boldsymbol{\xi}))=\deg(\boldsymbol{u},U,\boldsymbol{d}_t^{-1}(\boldsymbol{\xi}))=\begin{cases}
			1 & \text{if $\boldsymbol{\xi}\in \bigcup_{j\in J_0} V_j$,}\\
			0 & \text{if $\boldsymbol{\xi}\in \bigcup_{j\in J \setminus  J_0} V_j$.}
		\end{cases}
	\end{equation*} 
	Therefore
	\begin{equation*}
		\imt(\boldsymbol{d}_t \circ \boldsymbol{u},U)=\bigcup_{j\in J_0} V_j=\bigcup_{j\in J_0} \boldsymbol{d}_t(W_j)=\boldsymbol{d}_t(\imt(\boldsymbol{u},U)),	
	\end{equation*}
	so that \eqref{eqn:imtU-equivalence} is proved.
	
	To show \eqref{eqn:imt-equivalence} we proceed as follows. Let $(U_i)\subset \mathcal{U}_{\boldsymbol{u}}$ be such that $\imt(\boldsymbol{u},\Omega)=\bigcup_{i\in \N} \imt(\boldsymbol{u},U_i)$. For each $i\in \N$, we find a sequence $(H_i^j)_j$ of compact sets  contained in $\imt(\boldsymbol{u},U_i)$  whose union coincides with  $\imt(\boldsymbol{u},U_i)$. Hence $\imt(\boldsymbol{u},\Omega)=\bigcup_{i,j\in \N} H^i_j$.
	By \cite[Proposition~4.7(b)]{HS}, for every $i,j\in \N$, we find $\widetilde{U}^i_j\in \mathcal{U}_{\boldsymbol{u}} \cap \mathcal{U}_{\boldsymbol{d}_t \circ \boldsymbol{u}}$ such that $H^i_j  \subset \imt(\boldsymbol{u},\widetilde{U}^i_j)$. Thus
	\begin{equation*}
		\imt(\boldsymbol{u},\Omega)=\bigcup_{i,j\in \N} H^i_j \subset \bigcup_{i,j\in \N} \imt(\boldsymbol{u},\widetilde{U}^i_j) \subset \imt(\boldsymbol{u},\Omega),
	\end{equation*}
	so that $\imt(\boldsymbol{u},\Omega)=\bigcup_{i,j\in \N} \imt(\boldsymbol{u},\widetilde{U}^i_j)$. Applying $\boldsymbol{d}_t$ at both sides of the previous identity and using \eqref{eqn:imtU-equivalence}, we deduce 
	\begin{equation*}
		\boldsymbol{d}_t (\imt(\boldsymbol{u},\Omega))=\bigcup_{i,j\in \N} \boldsymbol{d}_t(\imt(\boldsymbol{u},\widetilde{U}^i_j))= \bigcup_{i,j\in \N} \imt(\boldsymbol{d}_t \circ \boldsymbol{u},\widetilde{U}^i_j) \subset \imt(\boldsymbol{d}_t \circ \boldsymbol{u},\Omega).
	\end{equation*}
	Exchanging the roles of $\boldsymbol{u}$ and $\boldsymbol{d}_t \circ \boldsymbol{u}$, we prove the opposite inclusion, namely $\imt(\boldsymbol{d}_t \circ \boldsymbol{u},\Omega) \subset \boldsymbol{d}_t (\imt(\boldsymbol{u},\Omega))$. Thus, also \eqref{eqn:imt-equivalence} is proved. 
	
	Now, let $\boldsymbol{m}\in W^{1,2}(\imt(\boldsymbol{u},\Omega))$. By \cite[Theorem~11.53]{Leo} and \eqref{eqn:imt-equivalence}, $\boldsymbol{m}\circ \boldsymbol{d}_t^{-1}\in W^{1,2}(\imt(\boldsymbol{d}_t \circ \boldsymbol{u},\Omega);\SN)$ and $(\boldsymbol{d}_t \circ \boldsymbol{u},\boldsymbol{m}\circ \boldsymbol{d}_t^{-1})\in \mathcal{Q}^O$. We conclude that $\Upsilon_t \colon \mathcal{Q}^O \to \mathcal{Q}^O$ is well defined.
	
	\emph{Claim (ii).} By reasoning as in claim (i), but having $\boldsymbol{d}_t^{-1}$ in place of $\boldsymbol{d}_t$,  we prove that $(\boldsymbol{d}_t^{-1}\circ \boldsymbol{y},\boldsymbol{n}\circ \boldsymbol{d}_t) \in \mathcal{Q}^O$ for all $(\boldsymbol{y},\boldsymbol{n})\in \mathcal{Q}^O$ and that \eqref{eqn:imtU-equivalence-inv}--\eqref{eqn:imt-equivalence-inv} hold. This shows that $\Upsilon_t$ is bijective and provides the expression of its inverse.  
	
	\emph{Claim (iii).} We show only \eqref{eqn:seq-cont} as the proof of \eqref{eqn:seq-cont-inv} is totally analogous.
	Let $\boldsymbol{p}_n=(\boldsymbol{u}_n,\boldsymbol{m}_n)$ and $\boldsymbol{p}=(\boldsymbol{u},\boldsymbol{m})\in\mathcal{Q}^O$ be such that $\boldsymbol{p}_n \to \boldsymbol{p}$ in $\mathcal{Q}^O$.
	 We set $\boldsymbol{q}_n=(\boldsymbol{y}_n,\boldsymbol{n}_n)\coloneqq \Upsilon_{t_n}(\boldsymbol{p}_n)$ and $\boldsymbol{q}=(\boldsymbol{y},\boldsymbol{n})\coloneqq \Upsilon_t(\boldsymbol{p})$, so that $\boldsymbol{q}_n,\boldsymbol{q}\in \mathcal{Q}^O$  by (i). 
	
	Note that $\boldsymbol{u}_n \to \boldsymbol{u}$ in $L^1(\Omega;\RN)$ by the Sobolev embedding theorem. This yields $\boldsymbol{y}_n \to \boldsymbol{y}$ in $L^1(\Omega;\RN)$ in view of the Lipschitz continuity and the boundedness of $\boldsymbol{d}$. Using the chain rule and \eqref{eqn:pp}, we estimate
	\begin{equation*}
		\begin{split}
			\int_\Omega A(|D\boldsymbol{y}_n|)\,\d \boldsymbol{x} &\leq \int_\Omega A(\|D\boldsymbol{d}_{t_n}\|_{L^\infty(O;\RNN)}|D\boldsymbol{u}_n|)\,\d\boldsymbol{x}\\
			&\leq (1+\|D\boldsymbol{d}\|_{C^0([0,T];L^\infty(O;\RNN))})^{p_A} \int_\Omega A (|D\boldsymbol{u}_n|)\,\d\boldsymbol{x},
		\end{split}
	\end{equation*}
	where the right-hand side is uniformly bounded in view of \eqref{eqn:modular-norm}.   Hence, we conclude that $\boldsymbol{y}_n \wks \boldsymbol{y}$ in $W^{1,A}(\Omega;\RN)$. 
	
	Observe that 
	\begin{equation*}
		\int_{\imt(\boldsymbol{y}_n,\Omega)} |\boldsymbol{n}_n|^2\,\d\boldsymbol{\xi}=\leb(\imt(\boldsymbol{y}_n,\Omega))=\leb(\boldsymbol{d}_{t_n}(\imt(\boldsymbol{u}_n,\Omega)))\leq \leb(O),
	\end{equation*}
	where we used \eqref{eqn:imt-equivalence}. Thus,  $(\chi_{\imt(\boldsymbol{y}_n,\Omega)}\boldsymbol{n}_n)$ is bounded in $L^2(\RN;\RN)$.   By assumption 
	\begin{equation*}
		\text{$\chi_{\imt(\boldsymbol{u}_n,\Omega)}\boldsymbol{m}_n\wk \chi_{\imt(\boldsymbol{u},\Omega)}\boldsymbol{m}$ in $L^2(\RN;\RN)$.}
	\end{equation*}
	Also, given $\boldsymbol{\psi}\in C^0_{\rm c}(\RN;\RN)$, we have
	\begin{equation*}
		\text{$\chi_{O} \boldsymbol{\psi}\circ \boldsymbol{d}_{t_n}\det D \boldsymbol{d}_{t_n}\to \chi_{O} \boldsymbol{\psi}\circ \boldsymbol{d}_{t}\det D \boldsymbol{d}_{t} $ in $L^2(\RN;\RN)$}
	\end{equation*}
	thanks to the dominated convergence theorem. By applying  Corollary \ref{cor:cov} with $\boldsymbol{d}_{t_n}$ and $\boldsymbol{d}_t$, and then passing to the limit, as $n\to \infty$, we obtain
	\begin{equation*}
		\begin{split}
			\int_{\imt(\boldsymbol{y}_n,\Omega)} \boldsymbol{n}_n \cdot \boldsymbol{\psi}\,\d\boldsymbol{\xi}&=\int_{\imt(\boldsymbol{u}_n,\Omega)} \boldsymbol{m}_n \cdot \boldsymbol{\psi}\circ \boldsymbol{d}_{t_n} \det D \boldsymbol{d}_{t_n}\,\d\boldsymbol{w}\\
			&\to \int_{\imt(\boldsymbol{u}_n,\Omega)} \boldsymbol{m} \cdot \boldsymbol{\psi}\circ \boldsymbol{d}_t \det D \boldsymbol{d}_t\,\d\boldsymbol{w}=\int_{\imt(\boldsymbol{y},\Omega)} \boldsymbol{n} \cdot \boldsymbol{\psi}\,\d\boldsymbol{\xi}.
		\end{split}
	\end{equation*}
	This shows that $\chi_{\imt(\boldsymbol{y}_n,\Omega)}\boldsymbol{n}_n \wk \chi_{\imt(\boldsymbol{y},\Omega)}\boldsymbol{n}$ in $L^2(\RN;\RN)$. 
	Eventually, using \eqref{eqn:imt-equivalence} and again Corollary \ref{cor:cov} with $\boldsymbol{d}_{t_n}$ and $\boldsymbol{d}_t$, we estimate
	\begin{equation*}
		\begin{split}
			\int_{\imt(\boldsymbol{y}_n,\Omega)} |D\boldsymbol{n}_n|^2\,\d\boldsymbol{\xi}&=\int_{\boldsymbol{d}_{t_n}(\imt(\boldsymbol{u}_n,\Omega))} |D\boldsymbol{m}_n \circ \boldsymbol{d}_{t_n}^{-1}|^2\,\d\boldsymbol{\xi}\\
			&=\int_{\imt(\boldsymbol{u}_n,\Omega)} |D\boldsymbol{m}_n|^2\,\det D \boldsymbol{d}_{t_n}\,\d\boldsymbol{w}\leq \|\det D \boldsymbol{d}\|_{C^0([0,T];L^\infty(O))} \int_{\imt(\boldsymbol{u}_n,\Omega)} |D\boldsymbol{m}_n|^2\,\d\boldsymbol{w}.
		\end{split}
	\end{equation*}
	Thus, $(\chi_{\imt(\boldsymbol{y}_n,\Omega)}D\boldsymbol{n}_n)$ is bounded in $L^2(\RN;\RNN)$ and, with  the same argument as above, we are able to show that $\chi_{\imt(\boldsymbol{y}_n,\Omega)}D\boldsymbol{n}_n \wk \chi_{\imt(\boldsymbol{y},\Omega)}D\boldsymbol{n}$ in $L^2(\RN;\RNN)$. Therefore, $\boldsymbol{q}_n\to \boldsymbol{q}$ in $\mathcal{Q}^O$.
	
	\emph{Claim (iv).} Let $\boldsymbol{u}\in \mathcal{Y}^O$. We show that $(\boldsymbol{d}_t \circ \boldsymbol{u})\restr{\Lambda}=\boldsymbol{d}_t\circ \boldsymbol{u}\restr{\Lambda}$, which yields $\Upsilon_t(\mathcal{Q}_{\boldsymbol{id}}^O)\subset \mathcal{Q}_{\boldsymbol{d}_t}^O$. We argue
	as in \cite[Proposition 5.1]{MoCo}. By Proposition~\ref{prop:embedding}(i), there exists  $(\boldsymbol{u}_n)\subset C^1(\closure{\Omega};\RN)$  such that $\boldsymbol{u}_n \to \boldsymbol{u}$ in $W^{1,N-1}(\Omega;\RN)$. By the continuity of the trace operator, we have that $\boldsymbol{u}_n\restr{\Lambda} \to \boldsymbol{\boldsymbol{u}}\restr{\Lambda}$ in $L^1( \Lambda;\RN)$. Applying the dominated convergence theorem, up to subsequences, we obtain $(\boldsymbol{d}_t \circ \boldsymbol{u}_n)\restr{\Lambda} \to \boldsymbol{d}_t\circ \boldsymbol{u}\restr{\Lambda}$ in $L^1(\Lambda;\RN)$.  Observe that $(\boldsymbol{d}_t \circ \boldsymbol{u}_n)$ is bounded in $W^{N-1,1}(\Omega;\RN)$ and $\boldsymbol{d}_t \circ \boldsymbol{u}_n \to \boldsymbol{d}_t \circ \boldsymbol{u}$ almost everywhere~in $\Omega$. Therefore, $\boldsymbol{d}_t \circ \boldsymbol{u}_n \wk \boldsymbol{d}_t \circ \boldsymbol{u}$ in $W^{1,N-1}(\Omega;\RN)$, so that $(\boldsymbol{d}_t \circ \boldsymbol{u}_n)\restr{\Lambda} \to (\boldsymbol{d}_t \circ \boldsymbol{u})\restr{\Lambda}$ in $L^1(\Lambda;\RN)$ by the weak continuity of the trace operator. This proves the claim. 
	With the same reasoning, but having $\boldsymbol{d}_t^{-1}$ in place of $\boldsymbol{d}_t$, we show that $\Upsilon_t^{-1}(\mathcal{Q}^O_{\boldsymbol{d}_t})\subset \mathcal{Q}^O_{\boldsymbol{id}}$. 
\end{proof}

We aim at computing the time derivative of the auxiliary energy functional. To this end, we first need a preliminary result stating some consequences of our assumptions on $W$. 

\begin{lemma}[Multiplicative estimates]
	\label{lem:gamma}
	Assume that $W$ satisfies {\rm \ref{it:W4}}. Then, there exist $\eta_W=\eta_W(N)>0$ such that, for all $\boldsymbol{F},\boldsymbol{G}\in \RNN_+ $ with $|\boldsymbol{G}-\boldsymbol{I}|<\eta_W$ and $\boldsymbol{z}\in \SN$, there hold:
	\begin{align}
		\label{eqn:gamma1}
		W(\boldsymbol{G}\boldsymbol{F},\boldsymbol{z})+b_W &\leq \frac{N}{N-1} (W(\boldsymbol{F},\boldsymbol{z})+b_W),\\
		\label{eqn:gamma2}
		|\partial_{\boldsymbol{F}}W(\boldsymbol{G}\boldsymbol{F},\boldsymbol{z})\boldsymbol{F}^\top|&\leq \frac{a_W N^2}{N-1} (W(\boldsymbol{F},\boldsymbol{z})+b_W),\\
		\label{eqn:gamma3}
		|W(\boldsymbol{G}\boldsymbol{F},\boldsymbol{z})-W(\boldsymbol{F},\boldsymbol{z})|&\leq \frac{a_W N^2}{N-1} (W(\boldsymbol{F},\boldsymbol{z})+b_W)|\boldsymbol{G}-\boldsymbol{I}|.
	\end{align}
\end{lemma}
\begin{proof}
	The estimates \eqref{eqn:gamma1}--\eqref{eqn:gamma2} are proved in \cite[Lemma 2.4(a)]{Ball.op} (see also \cite[Proposition 5.2]{FraMie}). For the proof of \eqref{eqn:gamma3}, we refer to \cite[Proposition 2.2]{Laz}.
\end{proof}

We proceed with our task. For convenience, we denote by $\mathcal{J}^{\rm e}$ and $\mathcal{J}^{\rm n}$ the functionals determined by the two summands on the right-hand side of \eqref{eqn:functional-J}, and by $\mathcal{M}^{\rm b}$, $\mathcal{M}^{\rm s}$, and $\mathcal{M}^{\rm f}$ the ones determined by the three terms on the right-hand side of \eqref{eqn:functional-K}.

\begin{proposition}[Time derivative of the auxiliary  energy]
	\label{prop:time-derivative-td}
	Assume that $W$ satisfies {\rm \ref{it:W4}} and {\rm \ref{it:W6}}. Let $\boldsymbol{d}$, $\boldsymbol{f}$, $\boldsymbol{g}$, and   $\boldsymbol{h}$ be as in \eqref{eqn:bd}--\eqref{eqn:bd-inverse} and  \eqref{eqn:forces}--\eqref{eqn:field}, respectively. Denote by $P\subset [0,T]$ the complement of the set  of times at which all these functions are differentiable.
	 Then, the following hold:
	\begin{enumerate}[label=(\roman*)]
		\item Let $t\in (0,T)\setminus P$ and $\boldsymbol{p}=(\boldsymbol{u},\boldsymbol{m})\in \mathcal{Q}^O$ with $\mathcal{J}(t,\boldsymbol{p})<+\infty$. Then, we have
		\begin{align}
			\label{eqn:Je-time-derivative}
			\partial_t \mathcal{J}^{\rm e}(t,\boldsymbol{p})&=\int_\Omega \Big (\partial_{\boldsymbol{F}} W \big  ((D\boldsymbol{d}_t \circ \boldsymbol{u}) D\boldsymbol{u},\boldsymbol{m}\circ \boldsymbol{u} \big )(D\boldsymbol{u}  )^\top\Big)  : (D \dot{\boldsymbol{d}}_t\circ \boldsymbol{u})\,\d\boldsymbol{x},\\
			\label{eqn:Jn-time-derivative}
			\begin{split}
				\partial_t \mathcal{J}^{\rm n}(t,\boldsymbol{p})
				&= 	\int_{\imt(\boldsymbol{u},\Omega)}  |D\boldsymbol{m}(D\boldsymbol{d}_t)^{-1}|^2\boldsymbol{I} :\left((D\dot{\boldsymbol{d}}_t)(D\boldsymbol{d}_t)^{-1} \right)\,\det D \boldsymbol{d}_t\,\d\boldsymbol{w}\\
				&-2 	\int_{\imt(\boldsymbol{u},\Omega)} \Big((D\boldsymbol{d}_t)^{-\top}(D\boldsymbol{m})^{\top}(D\boldsymbol{m})(D\boldsymbol{d}_t)^{-1}  \Big):\left((D\dot{\boldsymbol{d}}_t)(D\boldsymbol{d}_t)^{-1} \right)\,\det D \boldsymbol{d}_t\,\d\boldsymbol{w}.
			\end{split}
		\end{align}	
		\item Let  $t\in (0,T)\setminus P$ and $\boldsymbol{p}=(\boldsymbol{u},\boldsymbol{m})\in \mathcal{Q}^O$ with $\mathcal{M}(t,\boldsymbol{p})<+\infty$. Then, we have
		\begin{align}
			\label{eqn:Kb-time-derivative}
			\partial_t \mathcal{M}^{\rm b}(t,\boldsymbol{p})&=\int_\Omega \dot{\boldsymbol{f}}_t \cdot (\boldsymbol{d}_t \circ \boldsymbol{u})\,\d\boldsymbol{x}+\int_\Omega \boldsymbol{f}_t \cdot (\dot{\boldsymbol{d}}_t \circ \boldsymbol{u})\,\d\boldsymbol{x},\\
			\label{eqn:Ks-time-derivative}
			\partial_t \mathcal{M}^{\rm s}(t,\boldsymbol{p})&=\int_\Lambda \dot{\boldsymbol{g}}_t \cdot (\boldsymbol{d}_t \circ \boldsymbol{u}\restr{\Lambda})\,\d\haus+\int_\Lambda \boldsymbol{g}_t \cdot (\dot{\boldsymbol{d}}_t \circ \boldsymbol{u}\restr{\Lambda})\,\d\haus,\\
			\label{eqn:Kf-time-derivative}
			\begin{split}
				\partial_t \mathcal{M}^{\rm f}(t,\boldsymbol{p})&=\int_{\imt(\boldsymbol{u},\Omega)} \left( \dot{\boldsymbol{h}}_t \circ \boldsymbol{d}_t+ (D\boldsymbol{h}_t\circ\boldsymbol{d}_t)\dot{\boldsymbol{d}}_t\right)\det D\boldsymbol{d}_t \cdot \boldsymbol{m}\,\d\boldsymbol{w}\\
				&+ \int_{\imt(\boldsymbol{u},\Omega)} (\cof D \boldsymbol{d}_t : D\dot{\boldsymbol{d}}_t)\boldsymbol{h}_t \circ \boldsymbol{d}_t \cdot \boldsymbol{m}\,\d\boldsymbol{w}. 
			\end{split}
		\end{align}
	\end{enumerate}
\end{proposition}
\begin{proof}
Let us recall the sets $Z_{\boldsymbol{d}}$ and $\widetilde{Z}_{\boldsymbol{d},t}$ in \eqref{eq:Z}--\eqref{eq:Z-tilde}.

\emph{Claim (i).} First, we show \eqref{eqn:Je-time-derivative}.  For convenience,   we set $\boldsymbol{G}_{\tau}\coloneqq D\boldsymbol{d}_\tau\circ \boldsymbol{u}$ for all $\tau \in (0,T)$, $\boldsymbol{F}\coloneqq D\boldsymbol{u}$, and $\boldsymbol{z}\coloneqq \boldsymbol{m}\circ \boldsymbol{u}$. Let $h\in \R$ with $h\neq 0$ be such that $t+h\in (0,T)$. We write 
	\begin{equation}\label{eqn:Je-td}
		\frac{1}{h} \big (\mathcal{J}^{\rm e}(t+h,\boldsymbol{p})- \mathcal{J}^{\rm e}(t,\boldsymbol{p}) \big)=\int_\Omega \frac{W(\boldsymbol{G}_{t+h}\boldsymbol{F},\boldsymbol{z}) - W(\boldsymbol{G}_t\boldsymbol{F},\boldsymbol{z})}{h}\,\d\boldsymbol{x}.
	\end{equation} 
	By \ref{it:W4}, setting $\widetilde{\boldsymbol{G}}_\tau\coloneqq D\dot{\boldsymbol{d}}_\tau \circ \boldsymbol{u}$ for all $\tau\in (0,T)\setminus P$,  we have 
	\begin{equation*}
		\frac{W(\boldsymbol{G}_{t+h}\boldsymbol{F},\boldsymbol{z}) - W(\boldsymbol{G}_t\boldsymbol{F},\boldsymbol{z})}{h}\to \partial_{\boldsymbol{F}}W(\boldsymbol{G}_t\boldsymbol{F},\boldsymbol{z}):\widetilde{\boldsymbol{G}}_t\boldsymbol{F}= \Big( \partial_{\boldsymbol{F}}W(\boldsymbol{G}_t\boldsymbol{F},\boldsymbol{z})\boldsymbol{F}^\top\Big):\widetilde{\boldsymbol{G}}_t \quad \text{as $h\to 0$,}
	\end{equation*}
	where the convergence holds pointwise in $\{ \det D \boldsymbol{u}>0 \} \setminus (\boldsymbol{u}^{-1}(Z_{\boldsymbol{d}}) \cup \boldsymbol{u}^{-1}(\widetilde{Z}_{\boldsymbol{d},t}))$.  As $\boldsymbol{u}$ satisfies Lusin's condition (N${}^{-1}$), the convergence holds almost everywhere in $\Omega$ . Thus, in order to establish \eqref{eqn:Je-time-derivative}, it is sufficient to find a suitable domination for the integrand on the right-hand side of \eqref{eqn:Je-td} and to apply the dominated convergence theorem. 
	We estimate
	\begin{equation*}
		\begin{split}
			|\boldsymbol{G}_{t+h}\boldsymbol{G}_t^{-1}-\boldsymbol{I}|&\leq |\boldsymbol{G}_{t+h}-\boldsymbol{G}_t|\,|\boldsymbol{G}_t^{-1}|\leq \|D \boldsymbol{d}^{-1}_t\|_{L^\infty(O;\RNN)} \|D\boldsymbol{d}_{t+h}-D\boldsymbol{d}_t\|_{L^\infty(O;\RNN)} \\
			&\leq \|D \boldsymbol{d}^{-1}_t\|_{L^\infty(O;\RNN)} \|\dot{\boldsymbol{d}}\|_{L^\infty(0,T;\Lip(\closure{O};\closure{O}))}  |h|, 
		\end{split}
	\end{equation*}
	where the estimate is uniform with respect to $\boldsymbol{x}\in \Omega$. 
	Taking $|h|\ll 1$ so that the right-hand side is smaller than the constant $\eta_W$ in Lemma \ref{lem:gamma}, the estimate \eqref{eqn:gamma3} gives
	\begin{equation*}
		\begin{split}
			|W(\boldsymbol{G}_{t+h}\boldsymbol{F},\boldsymbol{z})-W(\boldsymbol{G}_t\boldsymbol{F},\boldsymbol{z})|&=|W(\boldsymbol{G}_{t+h}\boldsymbol{G}_t^{-1}\boldsymbol{G}_t\boldsymbol{F},\boldsymbol{z})-W(\boldsymbol{G}_t\boldsymbol{F},\boldsymbol{z})|\\
			&\leq C(N,a_W,\dot{\boldsymbol{d}},\boldsymbol{d}^{-1}) (W(\boldsymbol{G}_t\boldsymbol{F},\boldsymbol{z})+b_W)|h|.
		\end{split}
	\end{equation*}
	As $W(\boldsymbol{G}_t\boldsymbol{F},\boldsymbol{z})\in L^1(\Omega)$ by $\mathcal{J}^{\rm e}(t,\boldsymbol{p})<+\infty$, this concludes the proof of \eqref{eqn:Je-time-derivative}.
	
	We move to the proof of \eqref{eqn:Jn-time-derivative}. 
	Given $\boldsymbol{M}\in\RNN$, we consider the smooth function $\phi\colon \RNN_+\to [0,+\infty]$ defined as ${\phi}(\boldsymbol{A})\coloneqq |\boldsymbol{M}\boldsymbol{A}^{-1}|^2\det \boldsymbol{A}$. For all $\boldsymbol{A}\in\RNN_+$ and $\boldsymbol{B}\in\RNN$, a direct computation making use of \eqref{eqn:cof}--\eqref{eqn:inv} yields
	\begin{equation}
		\label{eq:differential-1st1}
		\langle \d \phi(\boldsymbol{A}),\boldsymbol{B}\rangle = -2 (\boldsymbol{M}\boldsymbol{A}^{-1}):(\boldsymbol{M}\boldsymbol{A}^{-1}\boldsymbol{B}\boldsymbol{A}^{-1})  \det \boldsymbol{A}+|\boldsymbol{M}\boldsymbol{A}^{-1}|^2 \cof\boldsymbol{A}:\boldsymbol{B},
	\end{equation}
	so that
	\begin{equation}
		\label{eq:differential-1st}
		|\langle \d \phi(\boldsymbol{A}),\boldsymbol{B}\rangle| \leq C(N)| \left(|\boldsymbol{A}|^N|\boldsymbol{A}^{-1}|^3+|\boldsymbol{A}|^{N-1}|\boldsymbol{A}^{-1}|^2   \right)|\boldsymbol{M}|^2|\boldsymbol{B}|. 
	\end{equation}
	Then, the first-order Taylor expansion yields
	\begin{equation}
		\label{eq:taylor-1st}
		|\phi(\boldsymbol{A}+\boldsymbol{B})-\phi(\boldsymbol{A})|\leq \sup_{\upsilon \in [0,1]} |\langle \d \phi(\boldsymbol{A}+\upsilon\boldsymbol{B}),\boldsymbol{B}\rangle|
	\end{equation}
	as long as $\boldsymbol{A}+\upsilon \boldsymbol{B}\in \RNN_+$ for all $\upsilon\in [0,1]$.
	
	Now, set $\boldsymbol{M}\coloneqq D\boldsymbol{m}$, and $\boldsymbol{D}_\tau\coloneqq D\boldsymbol{d}_\tau$ and $\widetilde{\boldsymbol{D}}_\tau\coloneqq D \dot{\boldsymbol{d}}_\tau$ for all $\tau \in (0,T)$. For $h\in \R$ with $h\neq 0$ such that $t+h\in (0,T)$, we write
	\begin{equation}\label{eqn:Jn-td}
		\frac{1}{h} \big (\mathcal{J}^{\rm n}(t+h,\boldsymbol{p})- \mathcal{J}^{\rm n}(t,\boldsymbol{p}) \big)=\int_{\imt(\boldsymbol{u},\Omega)} \frac{\phi(\boldsymbol{D}_{t+h})-\phi(\boldsymbol{D}_t)}{h}\,\d\boldsymbol{w}.
	\end{equation}
	 Using the chain rule, we compute 
	\begin{equation*}
		\begin{split}
			\frac{\phi(D_{t+h})-\phi(D_t)}{h}&\to  \langle \d \phi(\boldsymbol{D}_t),\widetilde{\boldsymbol{D}}_t\rangle\\ &=-2(\boldsymbol{M}\boldsymbol{D}_t^{-1}):\left( \boldsymbol{M}\boldsymbol{D}_t^{-1}\widetilde{\boldsymbol{D}}_t\boldsymbol{D}_t^{-1}  \right)\det \boldsymbol{D}_t+|\boldsymbol{M}\boldsymbol{D}_t^{-1}|^2\cof \boldsymbol{D}_t:\widetilde{\boldsymbol{D}}_t\\
			&=\left( |\boldsymbol{M}\boldsymbol{D}_t^{-1}|^2\boldsymbol{I}-2\boldsymbol{D}_t^{-\top}\boldsymbol{M}^\top\boldsymbol{M}\boldsymbol{D}_t^{-1} \right):\left( \widetilde{\boldsymbol{D}}_t\boldsymbol{D}_t^{-1} \right)\det \boldsymbol{D}_t, \quad \text{as $h\to 0$,}
		\end{split}
	\end{equation*}
	where the convergence holds pointwise in $\imt(\boldsymbol{u},\Omega) \setminus (Z_{\boldsymbol{d}} \cup \widetilde{Z}_{\boldsymbol{d},t})$ and, in turn, almost everywhere in $\imt(\boldsymbol{u},\Omega)$.
	 
	We seek a domination for the integrand on the right-hand side of \eqref{eqn:Jn-td}. Let $\boldsymbol{A}=\boldsymbol{D}_t$ and $\boldsymbol{B}=\boldsymbol{D}_{t+h}-\boldsymbol{D}_t$. Note that $\boldsymbol{A}$ takes values in a compact subset $K\subset \RNN_+$, while
	\begin{equation}
		\label{eq:B}
		\|\boldsymbol{B}\|_{L^\infty(O,\RNN)}\leq C(\boldsymbol{d})|h|
	\end{equation} 
	thanks to \eqref{eqn:bd}. 
	Given an open  neighborhood  $U\subset \subset \RNN_+$
	of $K$ in $\RNN_+$, we can take $|h|<C(\boldsymbol{d})\dist(K;\partial U)$ according  \eqref{eq:B} in order to have $\boldsymbol{A}+\upsilon \boldsymbol{B}\in U$ for all $\upsilon\in [0,1]$. Then, setting
	\begin{equation}
		\label{eq:m}
		m\coloneqq \max_{\boldsymbol{C}\in \closure{U}}|\boldsymbol{C}|+\max_{\boldsymbol{C}\in \closure{U}}|\boldsymbol{C}^{-1}|,
	\end{equation}
	the estimates in \eqref{eq:differential-1st}--\eqref{eq:taylor-1st} and \eqref{eq:B} yield
	\begin{equation*}
		|\phi(\boldsymbol{D}_{t+h})-\phi(\boldsymbol{D}_t)|\leq C(N,m,\boldsymbol{d})|\boldsymbol{M}|^2|h|.
	\end{equation*}
	Given that $\boldsymbol{M}\in L^2(\imt(\boldsymbol{u},\Omega);\RNN)$ because of $\mathcal{J}^{\rm e}(t,\boldsymbol{p})<+\infty$ and \eqref{eqn:bd}--\eqref{eqn:bd-inverse}, the previous equation gives a suitable domination for the integrand on the right-hand side of \eqref{eqn:Je-td}. Therefore, formula \eqref{eqn:Jn-time-derivative} is established by applying the dominated convergence theorem. 
	
	\emph{Claim (ii).}  
	We look at \eqref{eqn:Kb-time-derivative}. For $h\in \R$ with $h\neq 0$ such that $t+h\in (0,T)$, we write
	\begin{equation*}
		\label{eqn:Lb-td}
		\frac{\mathcal{M}^{\rm b}(t+h,\boldsymbol{p})-\mathcal{M}^{\rm b}(t,\boldsymbol{p})}{h}=\int_\Omega \frac{\boldsymbol{f}_{t+h}-\boldsymbol{f}_t}{h} \cdot (\boldsymbol{d}_{t+h}\circ \boldsymbol{u})\,\d\boldsymbol{x} + \int_\Omega \boldsymbol{f}_t \cdot \frac{\boldsymbol{d}_{t+h}\circ \boldsymbol{u}-\boldsymbol{d}_t \circ \boldsymbol{u}}{h}\,\d\boldsymbol{x}.
	\end{equation*}
	Then, using H\"{o}lder inequality,  we estimate 
	\begin{equation*}
		\left | \int_\Omega \left ( \frac{\boldsymbol{f}_{t+h}-\boldsymbol{f}_t}{h}-\dot{\boldsymbol{f}}_t \right) \cdot (\boldsymbol{d}_{t+h}\circ \boldsymbol{u})\,\d\boldsymbol{x}  \right | \leq C(O) \left \| \frac{\boldsymbol{f}_{t+h}-\boldsymbol{f}_t}{h}-\dot{\boldsymbol{f}}_t  \right \|_{L^1(\Omega;\RN)}  
	\end{equation*}
	and
	\begin{equation*}
		\left | \int_\Omega \boldsymbol{f}_t \cdot \left ( \frac{\boldsymbol{d}_{t+h}\circ \boldsymbol{u}-\boldsymbol{d}_t \circ \boldsymbol{u}}{h}-\dot{\boldsymbol{d}}_t \circ \boldsymbol{u} \right )\,\d\boldsymbol{x}  \right | \leq \|\boldsymbol{f}_t\|_{L^1(\Omega;\RN)} \,\left \| \frac{\boldsymbol{d}_{t+h}-\boldsymbol{d}_t}{h}-\dot{\boldsymbol{d}}_t \right \|_{L^\infty(O;O)},  
	\end{equation*}
	where both right-hand sides go to zero, as $h\to 0$. Thus,  \eqref{eqn:Kb-time-derivative} immediately follows. 
	
	The proof of \eqref{eqn:Ks-time-derivative} is totally analogous and we move to that of \eqref{eqn:Kf-time-derivative}. Define $\boldsymbol{k}\colon\tau \mapsto \boldsymbol{k}_\tau \coloneqq \boldsymbol{h}_\tau \circ \boldsymbol{d}_\tau \det D \boldsymbol{d}_\tau$. By \cite[Theorem 11.53]{Leo}, this defines a function $\boldsymbol{k}\colon [0,T]\to W^{1,1}(O;\RN)$.
	For $h\in \R$ with $h\neq 0$ and $t+h\in (0,T)$, we have
	\begin{equation*}
		\frac{\mathcal{M}^{\rm f}(t+h,\boldsymbol{p})-\mathcal{M}^{\rm f}(t,\boldsymbol{p})}{h}=\int_{\imt(\boldsymbol{u},\Omega)} \frac{\boldsymbol{k}_{t+h}-\boldsymbol{k}_t}{h} \cdot \boldsymbol{m}\,\d\boldsymbol{w}.
	\end{equation*}
	By Lemma~\ref{lem:k}, we have $\boldsymbol{k}\in AC([0,T];L^1(O;\RN))$. Then
	\begin{equation*}
		\left | \int_{\imt(\boldsymbol{u},\Omega)} \left (\frac{\boldsymbol{k}_{t+h}-\boldsymbol{k}_t}{h}-\dot{\boldsymbol{k}}_t \right ) \cdot \boldsymbol{m}\,\d\boldsymbol{w} \right | \leq \left \| \frac{\boldsymbol{k}_{t+h}-\boldsymbol{k}_t}{h}-\dot{\boldsymbol{k}}_t \right\|_{L^1(O;\RN)}, 
	\end{equation*}
	where the right-hand side goes to zero, as $h\to 0$. Thus, given  \eqref{eq:time-derivative-k}, we recover \eqref{eqn:Kf-time-derivative}.
\end{proof}

The next result concludes our proof of the equivalence between the auxiliary formulation and the original one. 

\begin{lemma}[Equivalence of the two formulations]
	\label{lem:equiv}
	Assume that W satisfies {\rm \ref{it:W1}--\ref{it:W7}}. 
	Let $\boldsymbol{d}$, $\boldsymbol{f}$, $\boldsymbol{g}$, and $\boldsymbol{h}$ be as in \eqref{eqn:bd}--\eqref{eqn:bd-inverse} and \eqref{eqn:forces}--\eqref{eqn:field}, respectively. Eventually, consider the map $\Upsilon_t\colon \mathcal{Q}^O \to \mathcal{Q}^O$ for $t\in [0,T]$ defined as in \eqref{eqn:upsilon}. Let $\boldsymbol{p}\colon [0,T] \to \mathcal{Q}_{\boldsymbol{id}}^O$ with $\boldsymbol{p}\colon t\mapsto \boldsymbol{p}_t$ and $\boldsymbol{q}\colon [0,T] \to \mathcal{Q}^O$  with  $\boldsymbol{q}\colon t\mapsto \boldsymbol{q}_t$ be two functions satisfying $\boldsymbol{q}_t= \Upsilon_t(\boldsymbol{p}_t)$ for all $t\in [0,T]$. Then, $\boldsymbol{p}$ is an energetic solution of the rate-independent system $(\mathcal{Q}^O_{\boldsymbol{id}},\mathcal{E}\restr{\mathcal{Q}^O_{\boldsymbol{id}}},\mathcal{D}\restr{\mathcal{Q}^O_{\boldsymbol{id}} \times \mathcal{Q}^O_{\boldsymbol{id}} })$ if and only if $\boldsymbol{q}$ fulfils conditions (i)--(iii) in Theorem \ref{thm:existence-td}.
\end{lemma}
\begin{proof}
	By applying  Corollary \ref{cor:cov} with $\boldsymbol{d}_t$ and  Proposition \ref{prop:equivalence-formulation}(i), we find
	\begin{equation}
		\label{eq:id1}
		\text{$\mathcal{E}(t,\boldsymbol{q}_t)=\mathcal{F}(t,\boldsymbol{p}_t)$ for all $t\in [0,T]$.}
	\end{equation}
	The identities 
	\begin{equation}
		\label{eq:id2}
		\mathcal{D}(\boldsymbol{p}_t,\widehat{\boldsymbol{p}})=\mathcal{D}(\boldsymbol{q}_t,\Upsilon_t(\widehat{\boldsymbol{p}})), \qquad \mathcal{D}(\boldsymbol{q}_t,\widehat{\boldsymbol{q}})=\mathcal{D}(\boldsymbol{p}_t,\Upsilon_t^{-1}(\widehat{\boldsymbol{q}})) \qquad \text{for all $t\in [0,T]$
			and $\widehat{\boldsymbol{p}}, \widehat{\boldsymbol{q}}\in\mathcal{Q}^O$.}
	\end{equation}
	and
	\begin{equation}
		\label{eq:id3}
		\mathrm{Var}_{\mathcal{D}}(\boldsymbol{p};[0,t])=\mathrm{Var}_{\mathcal{D}}(\boldsymbol{q};[0,t]) \quad \text{for all $t\in [0,T]$.}
	\end{equation}
	are immediate. Eventually, combining Corollary \ref{cor:cov} with $\boldsymbol{d}_t$,  Proposition \ref{prop:equivalence-formulation}(i), and Proposition \ref{prop:time-derivative-td}, we get
	\begin{equation}
		\label{eq:id4}
		\partial_t \mathcal{E}(t,\boldsymbol{q}_t)+\mathcal{P}(t,\boldsymbol{q}_t)=\partial_t \mathcal{F}(t,\boldsymbol{p}_t) \quad \text{for all $t\in [0,T]$,}
	\end{equation}
	where we recall \eqref{eqn:P}.  Therefore, the desired equivalence follows \eqref{eq:id1}--\eqref{eq:id4} by taking into account claims (ii) and (iv) of Proposition \ref{prop:equivalence-formulation}.
\end{proof}

In regards to the previous result,  the measurability of $\boldsymbol{p}$ is equivalent to that
 of $\boldsymbol{q}$. This fact will be shown in Lemma \ref{lem:equiv-measurability} below.

\subsection{Proof of Theorem \ref{thm:existence-td} and Theorem~\ref{thm:existence-a}}

We move towards the proof of Theorem~\ref{thm:existence-a}. 
The following result constitutes the analogue of Proposition \ref{prop:energy-ti} for the auxiliary formulation. Some steps are simpler in this case because of the confinement condition given by $O$.

\begin{proposition}[Coercivity and lower semicontinuity of the auxiliary  energy]
	\label{prop:energy-td}
Assume that $W$ satisfies {\rm \ref{it:W1}--\ref{it:W7}}. Let $\boldsymbol{d}$, $\boldsymbol{f}$, $\boldsymbol{g}$, and $\boldsymbol{h}$ be as in \eqref{eqn:bd}--\eqref{eqn:bd-inverse} and \eqref{eqn:forces}--\eqref{eqn:field}, respectively. Then, the following hold:
\begin{enumerate}[label=(\roman*)]
	\item \emph{Coercivity:} There exist two constant $\widetilde{K}_1,\widetilde{K}_2>0$ and two Borel functions $\widetilde{\Gamma},\widetilde{\gamma}\colon (0,+\infty) \to [0,+\infty]$ with 
	\begin{equation}
		\label{eqn:eta-tilde}
		\lim_{\vartheta \to +\infty} \frac{\widetilde{\Gamma}(\vartheta)}{\vartheta}= \lim_{\vartheta \to +\infty}\frac{\widetilde{\gamma}(\vartheta)}{\vartheta}=+\infty, \qquad \lim_{\vartheta \to 0^+} \widetilde{\gamma}(\vartheta)=+\infty
	\end{equation}
	such that
	\begin{equation}\label{eqn:coercivity-td}
		\mathcal{F}(t,\boldsymbol{p})\geq \widetilde{K}_1 \left( \int_\Omega  \left\{ A(|D\boldsymbol{u}|) +  \widetilde{\Gamma}(|\adj D \boldsymbol{u}|)+   \widetilde{\gamma}(\det D \boldsymbol{u}) \right\}\,\d \boldsymbol{x} + \int_{\imt(\boldsymbol{u},\Omega)} |D\boldsymbol{m}|^2 \,\d\boldsymbol{w}     \right) - \widetilde{K}_2
	\end{equation}
	for all $t\in [0,T]$ and $\boldsymbol{p}=(\boldsymbol{u},\boldsymbol{m})\in \mathcal{Q}_{\boldsymbol{id}}^O$. 
	\item \emph{Compactness:} Let $t\in [0,T]$ and let $(\boldsymbol{p}_n)$ be a sequence in $ \mathcal{Q}_{\boldsymbol{id}}^O$ with $\sup_{n\in\N} \mathcal{F}(t,\boldsymbol{p}_n)<+\infty$. Then, there exists a not relabeled subsequence of $(\boldsymbol{p}_n)$ and $\boldsymbol{p}\in\mathcal{Q}^O_{\boldsymbol{id}}$ such that $\boldsymbol{p}_n \to \boldsymbol{p}$ in $\mathcal{Q}^O$ and, in addition, \eqref{eqn:compactness-determinant}--\eqref{eqn:compactness-director} and \eqref{eqn:compactness-composition} hold with $\boldsymbol{p}_n=(\boldsymbol{u}_n,\boldsymbol{m}_n)$ and $\boldsymbol{p}=(\boldsymbol{u},\boldsymbol{m})$ in place of $\boldsymbol{q}_n=(\boldsymbol{y}_n,\boldsymbol{n}_n)$ and $\boldsymbol{q}=(\boldsymbol{y},\boldsymbol{n})$, respectively.
	\item \emph{Lower semicontinuity:} Let $t\in [0,T]$ and let $(\boldsymbol{p}_n)$ be a sequence  in $ \mathcal{Q}_{\boldsymbol{id}}^O$ such that   $\boldsymbol{p}_n\to \boldsymbol{p}$ in $\mathcal{Q}^O$ for some $\boldsymbol{p}\in\mathcal{Q}^O$. Then, there holds
	\begin{equation*}
		\mathcal{F}(t,\boldsymbol{p})\leq \liminf_{n\to \infty} \mathcal{F}(t,\boldsymbol{p}_n).
	\end{equation*}
\end{enumerate}	
\end{proposition}
\begin{proof}
\emph{Claim (i).} Let $t\in [0,T]$ and $\boldsymbol{p}=(\boldsymbol{u},\boldsymbol{m})\in \mathcal{Q}_{\boldsymbol{id}}^O$.  We begin by proving \eqref{eqn:coercivity-td}. We have
\begin{equation*}
	\begin{split}
		|D\boldsymbol{u}|&\leq |(D\boldsymbol{d}_t)^{-1}\circ \boldsymbol{u}|\,|(D\boldsymbol{d}_t \circ \boldsymbol{u})D\boldsymbol{u}|\leq \|(D\boldsymbol{d}_t)^{-1}\|_{L^\infty(O;\RNN)}
		 \,|(D\boldsymbol{d}_t \circ \boldsymbol{u})D\boldsymbol{u}|.
	\end{split}
\end{equation*}
Thanks to the monotonicity of $A$ and \eqref{eqn:pp}, this yields the bound
\begin{equation*}
	A(|D\boldsymbol{u}|)\leq  \left ( C(\boldsymbol{d}^{-1})+1 \right )^{p_A}  A(|(D\boldsymbol{d}_t \circ \boldsymbol{u})D\boldsymbol{u}|),
\end{equation*}
so that
\begin{equation}\label{eqn:A-growth}
	 \int_\Omega   A(|(D\boldsymbol{d}_t \circ \boldsymbol{u})D\boldsymbol{u}|)\,\d\boldsymbol{x} \geq C(\boldsymbol{d}^{-1},p_A) \int_\Omega  A(|D\boldsymbol{u}|) \,\d\boldsymbol{x}.
\end{equation}
The rest of the argument is analogous to the one in \cite[p. 552]{MoCo}. Observe that 
\begin{equation*}
	\cof \big( (D\boldsymbol{d}_t \circ \boldsymbol{u})D\boldsymbol{u}  \big)=\cof (D\boldsymbol{d}_t \circ \boldsymbol{u}) \cof D\boldsymbol{u}, \qquad \det \big( (D\boldsymbol{d}_t \circ \boldsymbol{u})D\boldsymbol{u}  \big)=\det  (D\boldsymbol{d}_t \circ \boldsymbol{u}) \det D\boldsymbol{u}.
\end{equation*}
By \eqref{it:W5}, the function $\Gamma$ is increasing on $(\vartheta_2,+\infty)$. Set
\begin{equation*}
	m_1\coloneqq \sup_{t\in [0,T]} \esssup_{\vphantom{t\in [0,T]}\boldsymbol{w}\in O} |(\cof D \boldsymbol{d}_t(\boldsymbol{w}))^{-1}|.
\end{equation*}
Then
\begin{equation*}
	m_1\leq \|D\boldsymbol{d}\|_{C^0([0,T];L^\infty(O;\RNN))}\,\|\det D \boldsymbol{d}\|_{C^0([0,T];L^\infty(O))}.
\end{equation*}
Clearly
\begin{equation*}
	|\cof D \boldsymbol{u}|\leq |(\cof(D\boldsymbol{d}_t \circ \boldsymbol{u}))^{-1}|\,|(\cof D \boldsymbol{d}_t \circ \boldsymbol{u})(\cof D \boldsymbol{u})|\leq m_1 |(\cof D \boldsymbol{d}_t \circ \boldsymbol{u})\,(\cof D \boldsymbol{u})|
\end{equation*}
from which we deduce
\begin{equation}
	\label{eqn:BB}
	\Gamma \left ( \frac{1}{m_1} |\cof D \boldsymbol{u}| \right ) \leq \Gamma(|(\cof D \boldsymbol{d}_t \circ \boldsymbol{u})\,(\cof D \boldsymbol{u})|) \quad \text{a.e. on $\{|\cof D \boldsymbol{u}|\geq  \vartheta_2 m_1 \}$.}
\end{equation}
Therefore, defining  $\widetilde{\Gamma}\colon (0,+\infty) \to [0,+\infty]$ as
\begin{equation*}
	\widetilde{\Gamma}(\vartheta)\coloneqq \chi_{[\vartheta_2 m_1,+\infty)} \Gamma \left(\frac{\vartheta}{m_1}\right),
\end{equation*}
the superlinear growth condition for $\widetilde{\Gamma}$ in \eqref{eqn:eta-tilde} is fulfilled and, from \eqref{eqn:BB}, we obtain
\begin{equation*}
	\int_\Omega \Gamma(|\cof (D \boldsymbol{d}_t \circ \boldsymbol{u}) \cof D \boldsymbol{u}|)\,\d\boldsymbol{x}\geq \int_\Omega \widetilde{\Gamma}(|\cof D \boldsymbol{u}|)\,\d\boldsymbol{x}. 
\end{equation*}
Again by \eqref{it:W5}, the function $\gamma$ is decreasing on $(0,\vartheta_1)$ and increasing on $(\vartheta_2,+\infty)$. Set 
\begin{equation}
	\label{eqn:m45}
	m_2\coloneqq \inf_{t\in [0,T]} \essinf_{\boldsymbol{w}\in O} \det D \boldsymbol{d}_t(\boldsymbol{w}), \qquad m_3\coloneqq \sup_{t\in [0,T]} \esssup_{\boldsymbol{w}\in O} \det D \boldsymbol{d}_t(\boldsymbol{w})=\|\det D\boldsymbol{d}\|_{C^0([0,T];L^\infty(O))}. 
\end{equation}
Observe that $m_2>0$ by \eqref{eq:md}. Taking into account the monotonicity properties of $\gamma$ and the Lusin's condition (N${}^{-1}$) satisfied by $\boldsymbol{u}$, we find
\begin{equation}
	\label{eqn:b1}
	\text{$\gamma(m_2 \det D \boldsymbol{u})\leq \gamma(\det D \boldsymbol{d}_t\circ \boldsymbol{u}\det D \boldsymbol{u})$ \quad  a.e. in $\left \{\det D \boldsymbol{u}\geq \frac{\vartheta_2}{m_2}\right \}$}
\end{equation}
and
\begin{equation}
	\label{eqn:b2}
	\text{$\gamma(m_3\det D \boldsymbol{u})\leq \gamma(\det D \boldsymbol{d}_t \circ \boldsymbol{u}\det D \boldsymbol{u})$ \quad a.e. in $\left \{ \det D \boldsymbol{u}\leq \frac{\vartheta_1}{m_3} \right \}$}. 
\end{equation}
Define $\widetilde{\gamma}\colon (0,+\infty) \to [0,+\infty]$ by setting
\begin{equation*}
	\widetilde{\gamma}(\vartheta) \coloneqq   \chi_{(0,\vartheta_1/m_3]}(\theta) \gamma(m_3 \vartheta)+\chi_{ [\vartheta_2/m_2,+\infty)}(\vartheta) \gamma(m_2 \vartheta). 
\end{equation*}
It is immediate to check that $\widetilde{\gamma}$ satisfies the limit conditions in \eqref{eqn:eta-tilde}, while \eqref{eqn:b1}--\eqref{eqn:b2} yield
\begin{equation}
	\label{eqn:b-growth}
	\int_\Omega \gamma(\det D \boldsymbol{d}_t \circ \boldsymbol{u}\det D \boldsymbol{u})\,\d\boldsymbol{x}\geq \int_\Omega \widetilde{\gamma}(\det D \boldsymbol{u})\,\d\boldsymbol{x}.
\end{equation} 
Therefore, combining \ref{it:W5} and \eqref{eqn:A-growth}--\eqref{eqn:b-growth}, we get
\begin{equation}
	\label{eqn:coervitiy-Je}
	\mathcal{J}^{\rm e}(t,\boldsymbol{u},\boldsymbol{m})\geq C (\boldsymbol{d}^{-1},p_A) \int_\Omega A(|D\boldsymbol{u}|)\,\d\boldsymbol{x}+\int_\Omega \widetilde{\Gamma}(\cof D  \boldsymbol{u})\,\d\boldsymbol{x}+\int_\Omega \widetilde{\gamma}(\det D  \boldsymbol{u})\,\d\boldsymbol{x}.
\end{equation}
Next, from \eqref{eqn:bd}--\eqref{eqn:bd-inverse}, we easily deduce the estimate
\begin{equation}
	\label{eqn:coercivity-Jn}
	\mathcal{J}^{\rm n}(t,\boldsymbol{u},\boldsymbol{m}) \geq C(\boldsymbol{d})\int_{\imt(\boldsymbol{u},\Omega)}|D\boldsymbol{m}|^2\,\d\boldsymbol{w}.
\end{equation}
Recalling \eqref{eqn:bd}, we have	
\begin{equation*}
	|\mathcal{M}^{\rm b}(t,\boldsymbol{u},\boldsymbol{m})| \leq \|\boldsymbol{f}_t\|_{L^1(\Omega;\RN)}\|\boldsymbol{d}_t\|_{\Lip(\closure{O};\closure{O})}\leq \|\boldsymbol{f}\|_{C^0([0,T];L^1(\Omega;\RN))} \|\boldsymbol{d}\|_{C^0([0,T];\Lip(\closure{O};\closure{O}))}
\end{equation*}
thanks to H\"{o}lder's inequality. 	Similarly
\begin{equation*}
	|\mathcal{M}^{\rm s}(t,\boldsymbol{u},\boldsymbol{m})| \leq \|\boldsymbol{g}_t\|_{L^1(\Lambda;\RN)} \|\boldsymbol{d}_t\|_{\Lip(\closure{O};\closure{O})} \leq \|\boldsymbol{g}\|_{C^0([0,T];L^1(\Lambda;\RN))} \|\boldsymbol{d}\|_{C^0([0,T];\Lip(\closure{O};\closure{O}))}.
\end{equation*}
Using H\"{o}lder's inequality and Corollary \ref{cor:cov} with $\boldsymbol{d}_t$, we get
\begin{equation*}
	\begin{split}
		|\mathcal{M}^{\rm f}(t,\boldsymbol{u},\boldsymbol{m})|&\leq \int_{\imt(\boldsymbol{u},\Omega)} |\boldsymbol{h}_t \circ \boldsymbol{d}_t|\det D\boldsymbol{d}_t\,\d\boldsymbol{w}\\
		&=\int_{\boldsymbol{d}_t(\imt(\boldsymbol{u},\Omega))} |\boldsymbol{h}_t|\,\d \boldsymbol{\xi} \leq \|\boldsymbol{h}_t\|_{L^1(O;\RN)}\leq \|\boldsymbol{h}\|_{C^0([0,T];L^1(O;\RN))}.
	\end{split}
\end{equation*}
Therefore
\begin{equation}
	\label{eqn:coervitity-K}
	|\mathcal{M}(t,\boldsymbol{u,\boldsymbol{m}})|\leq C(\boldsymbol{d},\boldsymbol{f},\boldsymbol{g},\boldsymbol{h}) .
\end{equation}
Putting together \eqref{eqn:coervitiy-Je}--\eqref{eqn:coervitity-K}, we obtain \eqref{eqn:coercivity-td} for $\widetilde{K}_1=\widetilde{K}_1(\boldsymbol{d},\boldsymbol{d}^{-1},p_A)>0$ and $\widetilde{K}_2\coloneqq K_2(\boldsymbol{d},\boldsymbol{f},\boldsymbol{g},\boldsymbol{h}) >0$.

\emph{Claim (ii).} From \eqref{eqn:coercivity-td} and the confinement condition, we see that
\begin{equation*}
	\sup_{n\in \N} \left\{ \|\boldsymbol{u}_n\|_{W^{1,A}(\Omega;\RN)} + \|\widetilde{\Gamma}(|\cof D \boldsymbol{u}_n|)\|_{L^1(\Omega)} +\|\widetilde{\gamma}(\det D \boldsymbol{u}_n)\|_{L^1(\Omega)}+\|D\boldsymbol{m}_n\|_{L^2(\imt(\boldsymbol{u}_n,\Omega);\RNN)}  \right\}<+\infty, 
\end{equation*}
where $\boldsymbol{p}_n=(\boldsymbol{u}_n,\boldsymbol{m}_n)$. Therefore, the conclusion follows by applying Theorem \ref{thm:compactness}.

\emph{Claim (iii).}  Recall \eqref{eqn:upsilon} and define $\boldsymbol{q}_n\coloneqq \Upsilon_t(\boldsymbol{p}_n)$ and $\boldsymbol{q}\coloneqq \Upsilon_t(\boldsymbol{p})$. 
By Proposition~\ref{prop:equivalence-formulation}(ii), we have $\boldsymbol{q}_n \to \boldsymbol{q}$ in $\mathcal{Q}^O$. As in Proposition \ref{prop:energy-ti}(ii), this entails 
\begin{equation*}
	\mathcal{F}(t,\boldsymbol{p})=\mathcal{E}(t,\boldsymbol{q})\leq \liminf_{n\to \infty} \mathcal{E}(t,\boldsymbol{q}_n)=\liminf_{n\to \infty} \mathcal{F}(t,\boldsymbol{p}_n),
\end{equation*}
where the first and the last equality are justified by applying Corollary \ref{cor:cov} with $\boldsymbol{d}_t$.
\end{proof}

 The next result is the analogue of Proposition~\ref{prop:power-ti} for the auxiliary formulation.

\begin{proposition}[Estimates for the power of the auxiliary energy]
\label{prop:power-td}
Assume that $W$ satisfies {\rm \ref{it:W1}--\ref{it:W7}}. Let $\boldsymbol{d}$, $\boldsymbol{f}$, $\boldsymbol{g}$, and $\boldsymbol{h}$ be as in \eqref{eqn:bd}--\eqref{eqn:bd-inverse} and  \eqref{eqn:forces}--\eqref{eqn:field}, respectively. Denote by $P\subset [0,T]$ the complement of the set of times at which the functions $\boldsymbol{d}$, $\boldsymbol{f}$, $\boldsymbol{g}$, and $\boldsymbol{h}$ are all differentiable. 
 Then, the following hold:
\begin{enumerate}[label=(\roman*)]
	\item \emph{Control of the power:} There exist a  function $\widetilde{\lambda} \in L^1(0,T)$ and a constant $\widetilde{K}>0$ such that
	\begin{equation*}
		|\partial_t \mathcal{F}(t,\boldsymbol{p})|\leq \widetilde{\lambda}(t) \left( \mathcal{F}(t,\boldsymbol{p}) + \widetilde{K} \right) \quad \text{for all $t\in (0,T)\setminus P$ and $\boldsymbol{p}\in \mathcal{Q}_{\boldsymbol{id}}^O$.}
	\end{equation*}
	\item \emph{Modulus of continuity of the power:} Let $M>0$  and $t\in (0,T)\setminus P$. Then, for each $\varepsilon>0$, there exists a constant $\widetilde{\delta}=\widetilde{\delta}(M,t,\varepsilon)>0$ 
	such that
	\begin{equation*}
		\left | \frac{\mathcal{F}(t+h,\boldsymbol{p})-\mathcal{F}(t,\boldsymbol{p})}{h}-\partial_t \mathcal{F}(t,\boldsymbol{p}) \right |\leq \varepsilon \quad \text{for all $\boldsymbol{p}\in\mathcal{Q}_{\boldsymbol{id}}^O$ with $\mathcal{F}(t,\boldsymbol{p})\leq M$ and $h\in (-\widetilde{\delta},\widetilde{\delta})$.}
	\end{equation*}
\end{enumerate}
\end{proposition}
\begin{proof}	
We refer to the sets $Z_{\boldsymbol{d}}$ and $\widetilde{Z}_{\boldsymbol{d},t}$  in \eqref{eq:Z}--\eqref{eq:Z-tilde}. 	
	
\emph{Claim (i).} 
Let $t\in (0,T)\setminus P$ and $\boldsymbol{p}=(\boldsymbol{u},\boldsymbol{m})\in \mathcal{Q}_{\boldsymbol{id}}^O$ with $\mathcal{F}(t,\boldsymbol{p})<+\infty$. For convenience, set $\boldsymbol{G}_t\coloneqq D\boldsymbol{d}_t \circ \boldsymbol{u}$, $\widetilde{\boldsymbol{G}}_t\coloneqq D  \dot{\boldsymbol{d}}_t \circ \boldsymbol{u}$, $\boldsymbol{F}\coloneqq D\boldsymbol{u}$, and $\boldsymbol{z}\coloneqq \boldsymbol{m}\circ \boldsymbol{u}$. By Proposition \ref{prop:time-derivative-td}(i), we have
\begin{equation*}
	\partial_t \mathcal{J}^{\rm e}(t,\boldsymbol{p})=\int_\Omega \Big (\partial_{\boldsymbol{F}} W(\boldsymbol{G}_t \boldsymbol{F},\boldsymbol{z})\boldsymbol{F}^\top \Big): \widetilde{\boldsymbol{G}}_t\,\d\boldsymbol{x}=\int_\Omega \boldsymbol{K}(\boldsymbol{G}_t\boldsymbol{F}):\big(\widetilde{\boldsymbol{G}}_t \boldsymbol{G}_t^{-1}\big)\,\d\boldsymbol{x}.
\end{equation*}
With the aid of \ref{it:W4}, we can estimate
\begin{equation}
	\label{eqn:dJe}
	\begin{split}
		|\partial_t \mathcal{J}^{\rm e}(t,\boldsymbol{p})|&\leq C(\dot{\boldsymbol{d}},\boldsymbol{d}^{-1}) \int_\Omega \left( a_W W(\boldsymbol{G}_t\boldsymbol{F},\boldsymbol{z}) + b_W  \right)\,\d\boldsymbol{x}\\ &\leq C_1(a_W,\dot{\boldsymbol{d}},\boldsymbol{d}^{-1}) \left( \mathcal{J}^{\rm e}(t,\boldsymbol{p}) + C_2(\Omega,a_W,b_W) \right).
	\end{split}
\end{equation}
Set $\boldsymbol{M}\coloneqq D\boldsymbol{m}$, $\boldsymbol{D}_t\coloneqq D \boldsymbol{d}_t$, and $\widetilde{\boldsymbol{D}}_t\coloneqq D\dot{\boldsymbol{d}}_t$. Then, Proposition \ref{prop:time-derivative-td}(i) yields
\begin{equation*}
	\partial_t \mathcal{J}^{\rm n}(t,\boldsymbol{p})=\int_{\imt(\boldsymbol{u},\Omega)} \left (  |\boldsymbol{M}\boldsymbol{D}_t^{-1}|^2\boldsymbol{I}-2\boldsymbol{D}_t^{-\top}\boldsymbol{M}^\top\boldsymbol{M}\boldsymbol{D}_t^{-1}  \right) : \left (\widetilde{\boldsymbol{D}}_t \boldsymbol{D}_t^{-1}\right)\det \boldsymbol{D}_t\,\d\boldsymbol{w}.
\end{equation*}
This expression can be bounded as
\begin{equation}
	\label{eqn:dJn}
	|\partial_t \mathcal{J}^{\rm e}(t,\boldsymbol{p})|\leq  C(\boldsymbol{d},\dot{\boldsymbol{d}}) \int_{\imt(\boldsymbol{u},\Omega)} |\boldsymbol{M}\boldsymbol{D}_t^{-1}|^2 \det \boldsymbol{D}_t \,\d\boldsymbol{w} = C(\boldsymbol{d},\dot{\boldsymbol{d}}) \mathcal{J}^{\rm n}(t,\boldsymbol{p}).
\end{equation}
Using H\"{o}lder's inequality and the confinement condition, we estimate
\begin{equation}
	\label{eqn:dKbs}
	\begin{split}
		|\partial_t \mathcal{K}^{\rm b}(t,\boldsymbol{p})+\partial_t \mathcal{K}^{\rm s}(t,\boldsymbol{p})| &\leq C(O) \left(  \|\dot{\boldsymbol{f}}_t\|_{L^1(\Omega;\RN)} + \|\dot{\boldsymbol{g}}_t\|_{L^1(\Sigma;\RN)} \right)\\ &+ \|\dot{\boldsymbol{d}}\|_{L^\infty(0,T;L^\infty(O;O))} \left ( \|\boldsymbol{f}\|_{C^0([0,T];L^1(\Omega;\RN))} + \|\boldsymbol{g}\|_{C^0([0,T];L^1(\Sigma;\RN))}\right ),
	\end{split}
\end{equation}
while applying Corollary \ref{cor:cov} with $\boldsymbol{d}_t$, we obtain
\begin{equation}
	\label{eqn:dKf}
	|\partial_t \mathcal{K}^{\rm f}(t,\boldsymbol{p})|\leq \int_{\imt(\boldsymbol{u},\Omega)} |\dot{\boldsymbol{h}}_t \circ \boldsymbol{d}_t| \det D\boldsymbol{d}_t\,\d\boldsymbol{w}=\int_{\boldsymbol{d}_t(\imt(\boldsymbol{u},\Omega))} |\dot{\boldsymbol{h}}_t|\,\d\boldsymbol{\xi}\leq C(O) \|\dot{\boldsymbol{h}}_t\|_{L^1(O;\RN)}. 
\end{equation}
Therefore, the combination of \eqref{eqn:dJe}--\eqref{eqn:dKf} gives the desired estimate.
\vskip 3pt

\emph{{Claim (ii).}} Let $M>0$,  $t\in (0,T)\setminus P$, $\varepsilon>0$, and $\boldsymbol{p}=(\boldsymbol{u},\boldsymbol{m})\in\mathcal{Q}^O_{\boldsymbol{id}}$ with $\mathcal{F}(t,\boldsymbol{p})\leq M$. First, note that 
\begin{equation*}
	\mathcal{J}(t,\boldsymbol{p})\leq C(\boldsymbol{d},\boldsymbol{f},\boldsymbol{g},\boldsymbol{h},M)
\end{equation*}
as a consequence of \eqref{eqn:coervitity-K}. For convenience, we divide the rest of the proof into four steps.

{\bf \em Step 1 (Elastic term).} 
We employ the notation  $\boldsymbol{F}\coloneqq D \boldsymbol{u}$,  $\boldsymbol{z}\coloneqq \boldsymbol{m}\circ \boldsymbol{u}$, and $\boldsymbol{G}_\tau\coloneqq D\boldsymbol{d}_\tau \circ \boldsymbol{u}$ and $\widetilde{\boldsymbol{G}}_\tau\coloneqq D\dot{\boldsymbol{d}}_\tau \circ \boldsymbol{u}$ for all $\tau\in [0,T]$. We consider the case $0<h<T-t$, the case $h<0$ being analogous.

For $\upsilon \in [0,1]$, we set $\boldsymbol{H}_\upsilon\coloneqq \boldsymbol{G}_t+\upsilon (\boldsymbol{G}_{t+h}-\boldsymbol{G}_t)$.
Observe that 
\begin{equation}
	\label{eq:H}
	\|\boldsymbol{H}_\upsilon-\boldsymbol{G}_t\|_{L^\infty(\Omega \setminus \boldsymbol{u}^{-1}(Z_{\boldsymbol{d}}))}\leq C(\boldsymbol{d})h
\end{equation}
thanks to \eqref{eqn:bd}. 
For every $\tau\in [0,T]$, the map $\boldsymbol{G}_\tau$ takes values in a compact set $K\subset \RNN_+$. Fix an open neighborhood $U\subset \subset \RNN_+$  of $K$ in $\RNN_+$. By \eqref{eq:H}, for $h\ll 1$, we  have that  $\boldsymbol{H}_\upsilon$ takes values in $U$ for all $\upsilon\in [0,1]$. Thus, using again \eqref{eq:H} together with the Lipschitz continuity of the determinant and the inversion of matrices in $\closure{U}$, we find 
\begin{align}
 	\label{eq:HH1}
 	\sup_{\upsilon \in [0,1]} \|\det \boldsymbol{H}_\upsilon-\det \boldsymbol{G}_t\|_{L^\infty(\Omega\setminus \boldsymbol{u}^{-1}(Z_{\boldsymbol{d}}))}&\leq C(\boldsymbol{d},U)h,\\
 	\label{eq:HH2}
 	 \sup_{\upsilon \in [0,1]} \| \boldsymbol{H}_\upsilon^{-1}-\boldsymbol{G}_t^{-1}\|_{L^\infty(\Omega\setminus \boldsymbol{u}^{-1}(Z_{\boldsymbol{d}});\RNN)}&\leq C(\boldsymbol{d},U)h.
\end{align}

Now, fix $\boldsymbol{x}\in \{ \det D \boldsymbol{u}>0\} \setminus \left( \boldsymbol{u}^{-1}(Z_{\boldsymbol{d}}\cup \widetilde{Z}_{\boldsymbol{d},t})  \right)$. Applying the mean value theorem to $\upsilon \mapsto W(\boldsymbol{H}_\upsilon(\boldsymbol{x})\boldsymbol{F}(\boldsymbol{x}))$, we find $\bar{\upsilon}=\bar{\upsilon}(t,h,\boldsymbol{x})\in [0,1]$ such that
\begin{equation*}
	\frac{W(\boldsymbol{G}_{t+h}(\boldsymbol{x})\boldsymbol{F}(\boldsymbol{x}),\boldsymbol{z}(\boldsymbol{x}))-W(\boldsymbol{G}_t(\boldsymbol{x})\boldsymbol{F}(\boldsymbol{x}),\boldsymbol{z}(\boldsymbol{x}))}{h}=\partial_{\boldsymbol{F}} W(\boldsymbol{H}_{\bar{\upsilon}}(\boldsymbol{x})\boldsymbol{F}(\boldsymbol{x}),\boldsymbol{z}(\boldsymbol{x}))\boldsymbol{F}^\top(\boldsymbol{x}):\frac{\boldsymbol{G}_{t+h}(\boldsymbol{x})-\boldsymbol{G}_t(\boldsymbol{x})}{h}.
\end{equation*}
Henceforth, the point $\boldsymbol{x}$ remains fixed, but we omit it from the formulas. 
We estimate
\begin{equation}
	\label{eq:W1}
	\begin{split}
		\Big |\partial_{\boldsymbol{F}} W(\boldsymbol{H}_{\bar{\upsilon}}\boldsymbol{F},\boldsymbol{z})\boldsymbol{F}^\top :& \frac{\boldsymbol{G}_{t+h}-\boldsymbol{G}_t}{h} - \partial_{\boldsymbol{F}} W (\boldsymbol{G}_t\boldsymbol{F},\boldsymbol{z}) \boldsymbol{F}^\top : \widetilde{\boldsymbol{G}}_t \Big |\\
		&\leq \left |\partial_{\boldsymbol{F}} W(\boldsymbol{H}_{\bar{\upsilon}}\boldsymbol{F},\boldsymbol{z})\boldsymbol{F}^\top-\partial_{\boldsymbol{F}} W(\boldsymbol{G}_t\boldsymbol{F},\boldsymbol{z})\boldsymbol{F}^\top \right |\,\bigg | \frac{\boldsymbol{G}_{t+h}-\boldsymbol{G}_t}{h} \bigg |\\
		&+ \Big |  \partial_{\boldsymbol{F}} W(\boldsymbol{G}_t\boldsymbol{F},\boldsymbol{z})\boldsymbol{F}^\top \Big |\, \bigg | \frac{\boldsymbol{G}_{t+h}-\boldsymbol{G}_t}{h} - \widetilde{\boldsymbol{G}}_t \bigg |.
	\end{split}
\end{equation}
We look at the first term on the right-hand side of \eqref{eq:W1}. Setting $\boldsymbol{J}_{\bar{\upsilon}}\coloneqq \boldsymbol{H}_{\bar{\upsilon}} \boldsymbol{G}_t^{-1}$, we find
\begin{equation*}
	\begin{split}
		\left |\partial_{\boldsymbol{F}} W(\boldsymbol{H}_{\bar{\upsilon}}\boldsymbol{F},\boldsymbol{z})\boldsymbol{F}^\top-\partial_{\boldsymbol{F}} W(\boldsymbol{G}_t\boldsymbol{F},\boldsymbol{z})\boldsymbol{F}^\top \right | &= \left |\boldsymbol{K}(\boldsymbol{J}_{\bar{\upsilon}} \boldsymbol{G}_t\boldsymbol{F})\boldsymbol{H}_{\bar{\upsilon}}^{-\top} -\boldsymbol{K}(\boldsymbol{G}_t\boldsymbol{F})\boldsymbol{G}_t^{-\top} \right |\\
		&\leq  |\boldsymbol{K}(\boldsymbol{J}_{\bar{\upsilon}} \boldsymbol{G}_t\boldsymbol{F})-\boldsymbol{K}(\boldsymbol{G}_t\boldsymbol{F})|\,\left | \boldsymbol{H}_{\bar{\upsilon}}^{-1}  \right |+  |\boldsymbol{K}(\boldsymbol{G}_t\boldsymbol{F})|\,\left|\boldsymbol{H}_{\bar{\upsilon}}^{-1}-\boldsymbol{G}_t^{-1}\right|.
	\end{split}
\end{equation*}
By \eqref{eq:H}, we have
\begin{equation*}
	|\boldsymbol{J}_{\bar{\upsilon}}- \boldsymbol{I}|\leq C(\boldsymbol{d}^{-1}) |\boldsymbol{H}_{\bar{\upsilon}} - \boldsymbol{G}_t|\leq C(\boldsymbol{d},\boldsymbol{d}^{-1})h.
\end{equation*}
Let $\delta_W=\delta_W(\varepsilon)$ be as in \ref{it:W5}. For a suitable $\bar{h}=\bar{h}(\delta_W)>0$, the right-hand side of the previous equation is smaller than $\delta_W$ for all $h\leq \bar{h}$. In this case, \ref{it:W5} gives
\begin{equation*}
	|\boldsymbol{K}(\boldsymbol{J}_{\bar{\upsilon}} \boldsymbol{G}_t\boldsymbol{F})-\boldsymbol{K}(\boldsymbol{G}_t\boldsymbol{F})|\,\left | \boldsymbol{H}_{\bar{\upsilon}}^{-1}  \right | \leq C(\boldsymbol{d})  \left (W(\boldsymbol{G}_t\boldsymbol{F},\boldsymbol{z})+b_W \right )\varepsilon .
\end{equation*}
Using \ref{it:W4} and \eqref{eq:HH2}, we estimate
\begin{equation*}
	|\boldsymbol{K}(\boldsymbol{G}_t\boldsymbol{F})|\,\left|\boldsymbol{H}_{\bar{\upsilon}}^{-1}-\boldsymbol{G}_t^{-1}\right|\leq a_W C(\boldsymbol{d},\boldsymbol{d}^{-1},U) \left (W(\boldsymbol{G}_t\boldsymbol{F},\boldsymbol{z})+b_W \right ) h.
\end{equation*}
Thus, for the first term on the right-hand side of \eqref{eq:W1}, we obtain
\begin{equation*}
	\left |\partial_{\boldsymbol{F}} W(\boldsymbol{H}_{\bar{\upsilon}}\boldsymbol{F},\boldsymbol{z})\boldsymbol{F}^\top-\partial_{\boldsymbol{F}} W(\boldsymbol{G}_t\boldsymbol{F},\boldsymbol{z})\boldsymbol{F}^\top \right |\,\bigg | \frac{\boldsymbol{G}_{t+h}-\boldsymbol{G}_t}{h} \bigg |\leq a_W C(\boldsymbol{d},\boldsymbol{d}^{-1}) \left ( W(\boldsymbol{G}_t\boldsymbol{F},\boldsymbol{z})+b_W\right ) (\varepsilon+h).
\end{equation*}
Looking at the second term on the right-hand side of \eqref{eq:W1}, this is estimated with the aid of \ref{it:W4} as
\begin{equation*}
	\begin{split}
		 \Big |  \partial_{\boldsymbol{F}} W(\boldsymbol{G}_t\boldsymbol{F},\boldsymbol{z})\boldsymbol{F}^\top \Big |\, \bigg | \frac{\boldsymbol{G}_{t+h}-\boldsymbol{G}_t}{h} - \widetilde{\boldsymbol{G}}_t \bigg |&\leq C(\boldsymbol{d}^{-1})|\boldsymbol{K}(\boldsymbol{G}_t\boldsymbol{F})|\,\left | \frac{D\boldsymbol{d}_{t+h}\circ \boldsymbol{u}-D\boldsymbol{d}_t\circ \boldsymbol{u}}{h}- D\dot{\boldsymbol{d}}_t \circ \boldsymbol{u}  \right |\\
		 &\leq a_W C(\boldsymbol{d}^{-1}) \left ( W(\boldsymbol{G}_t\boldsymbol{F},\boldsymbol{z})+b_W \right)\omega_{\boldsymbol{d},t}(h), 
	\end{split}
\end{equation*}
where
\begin{equation}
	\label{eq:omomega}
	\omega_{t}(h)\coloneqq \left \|  \frac{\boldsymbol{d}_{t+h}-\boldsymbol{d}_t}{h} - \dot{\boldsymbol{d}}_t \right \|_{W^{1,\infty}(O;O)}.
\end{equation}
All these estimates give  the inequality
\begin{equation}
	\label{eq:Je}
	\left | \frac{\mathcal{J}^{\rm e}(t+h,\boldsymbol{p})-\mathcal{J}^{\rm e}(t,\boldsymbol{p})}{h}-\partial_t \mathcal{J}^{\rm e}(t,\boldsymbol{p}) \right | \leq C(\boldsymbol{d},\boldsymbol{d}^{-1},\boldsymbol{f},\boldsymbol{g},\boldsymbol{h},U,M,a_W,b_W,\Omega) (\varepsilon + h + \omega_{t}(h)). 
\end{equation}

{\bf \em Step 2 (Nematic term).} Fix 
 $\boldsymbol{M}\in\RNN$. As in the proof of Proposition \ref{prop:time-derivative-td}(ii), we consider the smooth function  $\phi\colon \RNN_+ \to (0,+\infty)$ given by $\phi(\boldsymbol{A})\coloneqq |\boldsymbol{M}\boldsymbol{A}^{-1}|^2\det \boldsymbol{A}$.  By a second-order Taylor expansion, we find
\begin{equation}
	\label{eq:taylor}
	\left | \phi(\boldsymbol{A}+\boldsymbol{B})-\phi(\boldsymbol{A})-\langle \d \phi(\boldsymbol{A}),\boldsymbol{B} \rangle  \right | \leq \frac{1}{2} \sup_{\upsilon \in [0,1]} |\langle \d^2 \phi \left(\boldsymbol{A}+\upsilon \boldsymbol{B}\right),\boldsymbol{B}\rangle | 
\end{equation}
for all $\boldsymbol{A}\in \RNN_+$ and $\boldsymbol{B}\in\RNN$ for which $\boldsymbol{A}+\upsilon \boldsymbol{B}\in\RNN_+$ for all $\upsilon\in [0,1]$. In the previous equation,   $\d^2 \phi(\boldsymbol{A}+\upsilon\boldsymbol{B})$ denotes a quadratic form over $\RNN$. In addition to \eqref{eq:differential-1st1}, using \eqref{eqn:cof}--\eqref{eqn:inv} , for all $\boldsymbol{A}\in\RNN_+$ and $\boldsymbol{B}\in\RNN$ we compute
\begin{equation*}
	\begin{split}
		\langle \d^2 \phi(\boldsymbol{A}),\boldsymbol{B}\rangle &=2|\boldsymbol{C}
		|^2\det\boldsymbol{A}+4  (\boldsymbol{M}\boldsymbol{A}^{-1}):(\boldsymbol{C}\boldsymbol{B}\boldsymbol{A}^{-1}) \det\boldsymbol{A}-4 \left( (\boldsymbol{M}\boldsymbol{A}^{-1}):\boldsymbol{C} \right)(\cof \boldsymbol{A}:\boldsymbol{B})\\
		&+|\boldsymbol{M}\boldsymbol{A}^{-1}|^2 \left( (\cof \boldsymbol{A}:\boldsymbol{B})\boldsymbol{I}-(\cof\boldsymbol{A})\boldsymbol{B}^\top  \right):(\boldsymbol{B}\boldsymbol{A}^{-1}),
	\end{split}
\end{equation*}
 where, for convenience, we set $\boldsymbol{C}\coloneqq \boldsymbol{M}\boldsymbol{A}^{-1}\boldsymbol{B}\boldsymbol{A}^{-1}$. Therefore
\begin{equation}
	\label{eq:second}
	\begin{split}
		|\langle \d^2\phi(\boldsymbol{A}),\boldsymbol{B}\rangle|&\leq C(N) \left( |\boldsymbol{A}|^N|\boldsymbol{A}^{-1}|^4+|\boldsymbol{A}|^2|\boldsymbol{A}^{-1}|^3+ |\boldsymbol{A}|^{N-1}|\boldsymbol{A}^{-1}|^3  \right) |\boldsymbol{M}|^2|\boldsymbol{B}|^2
	\end{split}
\end{equation}
for all $\boldsymbol{A}\in \RNN_+$ and $\boldsymbol{B}\in\RNN$.

Now, we choose $\boldsymbol{M}= D\boldsymbol{m}$, and we set $\boldsymbol{D}_\tau\coloneqq D \boldsymbol{d}_\tau$ and $\widetilde{\boldsymbol{D}}_\tau\coloneqq D \dot{\boldsymbol{d}}_\tau$ for $\tau\in [0,T]$. Observe that
\begin{equation*}
	\|\boldsymbol{M}\|_{L^2(\imt(\boldsymbol{u},\Omega);\RNN)}\leq C(M,K_1,K_2)	
\end{equation*}
by \eqref{eqn:coercivity-td} and the assumption $\mathcal{F}(t,\boldsymbol{p})\leq M$. 
We consider only $h>0$ since the case $h<0$ works analogously. 
We write
\begin{equation*}
	\begin{split}
		\left | \frac{\mathcal{J}^{\rm n}(t+h,\boldsymbol{p})-\mathcal{J}^{\rm n}(t,\boldsymbol{p})}{h}-\partial_t \mathcal{J}^{\rm n}(t,\boldsymbol{p}) \right |&\leq \int_{\imt(\boldsymbol{u},\Omega)} \left |  \frac{\phi(\boldsymbol{D}_{t+h})-\phi(\boldsymbol{D}_t)}{h} - \left  \langle \d \phi(\boldsymbol{D}_t),\widetilde{\boldsymbol{D}}_t \right  \rangle   \right |\,\d \boldsymbol{w}\\
		&\leq \frac{1}{h}\int_{\imt(\boldsymbol{u},\Omega)}  \left | {\phi(\boldsymbol{D}_{t+h})-\phi(\boldsymbol{D}_t)} - \left \langle \d \phi(\boldsymbol{D}_t),{\boldsymbol{D}_{t+h}-\boldsymbol{D}_t} \right  \rangle\right |\,\d\boldsymbol{w}\\
		&+\int_{\imt(\boldsymbol{u},\Omega)} \left | \left \langle \d \phi(\boldsymbol{D}_t), \frac{\phi(\boldsymbol{D}_{t+h})-\phi(\boldsymbol{D}_t)}{h} - \widetilde{\boldsymbol{D}}_t \right \rangle  \right | \,\d\boldsymbol{w}.
	\end{split}
\end{equation*}
Let $\boldsymbol{A}=\boldsymbol{D}_t$ and $\boldsymbol{B}=\boldsymbol{D}_{t+h}-\boldsymbol{D}_t$. The function $\boldsymbol{A}$ takes values in a compact set $K\subset \RNN_+$. Given an open neighborhood $U\subset \subset \RNN_+$ of $K$, from \eqref{eqn:bd} we see that $\boldsymbol{A}+\upsilon \boldsymbol{B}\in U$ for all $\upsilon\in [0,1]$ as long as $h\ll 1$ by arguing as in the proof of Proposition \ref{prop:time-derivative-td}(ii).

Defining $m>0$ as in \eqref{eq:m} and using \eqref{eq:B} and \eqref{eq:taylor}--\eqref{eq:second}, the first term on the right-hand side is estimated as
\begin{equation*}
	\begin{split}
		\frac{1}{h}\int_{\imt(\boldsymbol{u},\Omega)} \big | {\phi(\boldsymbol{D}_{t+h})-\phi(\boldsymbol{D}_t)} &-  \langle \d \phi(\boldsymbol{D}_t),{\boldsymbol{D}_{t+h}-\boldsymbol{D}_t}   \rangle\big |\,\d\boldsymbol{w}\\
		&\leq \frac{C(N,\boldsymbol{d},\boldsymbol{d}^{-1},r)}{h} \int_{\imt(\boldsymbol{u},\Omega)} |\boldsymbol{M}|^2|\boldsymbol{D}_{t+h}-\boldsymbol{D}_t|^2\,\d\boldsymbol{w}\\
		&\leq C(N,\boldsymbol{d},\boldsymbol{d}^{-1},m,M,K_1,K_2) h.
	\end{split}
\end{equation*}
For the second term on the right-hand side, we employ \eqref{eq:differential-1st} and \eqref{eq:m} to find
\begin{equation*}
	\begin{split}
		\int_{\imt(\boldsymbol{u},\Omega)} \Big | \Big \langle \d \phi(\boldsymbol{D}_t), &\frac{\phi(\boldsymbol{D}_{t+h})-\phi(\boldsymbol{D}_t)}{h} - \widetilde{\boldsymbol{D}}_t \Big \rangle  \Big | \,\d\boldsymbol{w}\\
		&\leq C(N,\boldsymbol{d},\boldsymbol{d}^{-1},r) \int_{\imt(\boldsymbol{u},\Omega)} |\boldsymbol{M}|^2 \left | \frac{\phi(\boldsymbol{D}_{t+h})-\phi(\boldsymbol{D}_t)}{h} - \widetilde{\boldsymbol{D}}_t \right |\,\d \boldsymbol{w}\\
		&\leq C(N,\boldsymbol{d},\boldsymbol{d}^{-1},m,M,K_1,K_2) \omega_{t}(h), 
	\end{split}
\end{equation*} 
where we recall \eqref{eq:omomega}
 Hence, we obtain
 \begin{equation}
 	\label{eq:Jn}
 		\left | \frac{\mathcal{J}^{\rm n}(t+h,\boldsymbol{p})-\mathcal{J}^{\rm n}(t,\boldsymbol{p})}{h}-\partial_t \mathcal{J}^{\rm n}(t,\boldsymbol{p}) \right | \leq  C(N,\boldsymbol{d},\boldsymbol{d}^{-1},m,M,K_1,K_2)  \left( h+\omega_{\boldsymbol{d},t}(h) \right).
 \end{equation}
 
{\bf \em Step 3 (Applied loads).} Define $\boldsymbol{k}\colon \tau \mapsto \boldsymbol{h}_\tau \circ \boldsymbol{d}_\tau \det D\boldsymbol{d}_\tau$ and recall Lemma \ref{lem:k}. 
By repeating the estimates in the proof of Proposition \ref{prop:time-derivative-td}(ii), we obtain
\begin{equation}
	\label{eq:M}
	\left | \frac{\mathcal{M}(t+h,\boldsymbol{p})-\mathcal{M}(t,\boldsymbol{p})}{h}-\partial_t \mathcal{M}(t,\boldsymbol{p})   \right |\leq C_1(O)\widetilde{\omega}_t(h)+C_2(\boldsymbol{f},\boldsymbol{g},\boldsymbol{h})\omega_t(h),
\end{equation}
where
\begin{equation*}
	\widetilde{\omega}_t(h)\coloneqq \left \| \frac{\boldsymbol{f}_{t+h}-\boldsymbol{f}_t}{h}-\dot{\boldsymbol{f}}_t  \right \|_{L^1(\Omega;\RN)}+\left \| \frac{\boldsymbol{g}_{t+h}-\boldsymbol{g}_t}{h}-\dot{\boldsymbol{g}}_t  \right \|_{L^1(\Lambda;\RN)}+\left \| \frac{\boldsymbol{k}_{t+h}-\boldsymbol{k}_t}{h}-\dot{\boldsymbol{k}}_t  \right \|_{L^1(O;\RN)}.
\end{equation*}

{\bf \em Step 4 (Conclusion).} By putting together \eqref{eq:Je}, \eqref{eq:Jn}, and \eqref{eq:M}, we find
\begin{equation*}
	\begin{split}
		\left | \frac{\mathcal{F}(t+h,\boldsymbol{p})-\mathcal{F}(t,\boldsymbol{p})}{h}-\partial_t \mathcal{F}(t,\boldsymbol{p})   \right |&\leq C(\varepsilon+2h+3\omega_t(h)+\widetilde{\omega}_t(h))
	\end{split}
\end{equation*}
for some constant $C=C(N,\Omega,\boldsymbol{d},\boldsymbol{d}^{-1},\boldsymbol{f},\boldsymbol{g},\boldsymbol{h},U,M,K_1,K_2)>0$. Therefore, the claim follows up to replacing $\varepsilon$ with $\varepsilon/C$ and choosing $\widetilde{\delta}(\varepsilon)$ depending on $\omega_t$ and $\widetilde{\omega}_t$.

\end{proof}

We are ready to  give the proof of Theorem~\ref{thm:existence-a}.

\begin{proof}[Proof of Theorem~\ref{thm:existence-a}]
As for the one of Theorem \ref{thm:existence-ti}, also this proof  relies on the application of Theorem~\ref{thm:MR}. We consider the rate-independent system $(\mathcal{Q}^O_{\boldsymbol{id}},\mathcal{F}\restr{\mathcal{Q}^O_{\boldsymbol{id}}},\mathcal{D}\restr{\mathcal{Q}^O_{\boldsymbol{id}}\times \mathcal{Q}^O_{\boldsymbol{id}}})$. The space $\mathcal{Q}^O_{\boldsymbol{id}}$ satisfies \ref{it:Q}, assumption \ref{it:E1} follows from Proposition~\ref{prop:energy-td} by arguing as in the proof of Theorem \ref{thm:existence-ti}, and   \ref{it:E2}--\ref{it:E3} are shown in Proposition \ref{prop:power-td}. Eventually, conditions \ref{it:D1}--\ref{it:D2} are immediate and \ref{it:D3} is checked by using Proposition \ref{prop:energy-td}(ii).
\end{proof}

Next, we discuss the measurability of the functions $\boldsymbol{p}$ and $\boldsymbol{q}$ in Lemma~\ref{lem:equiv}. 

\begin{lemma}[Measurability]
	\label{lem:equiv-measurability}
		Assume that W satisfies {\rm \ref{it:W1}--\ref{it:W7}}.
	Let $\boldsymbol{d}$, $\boldsymbol{f}$, $\boldsymbol{g}$, and $\boldsymbol{h}$ be as in \eqref{eqn:bd}--\eqref{eqn:field}. 
	Eventually, consider $\Upsilon_t\colon \mathcal{Q}^O \to \mathcal{Q}^O$ for $t\in [0,T]$ be the map introduced in Proposition~\ref{prop:equivalence-formulation}. Let $\boldsymbol{p}\colon [0,T]\to \mathcal{Q}^O_{\boldsymbol{id}}$ with $ \boldsymbol{p}\colon t\mapsto \boldsymbol{p}_t$ and $\boldsymbol{q}\colon [0,T]\to \mathcal{Q}^O$ with $\boldsymbol{q}\colon t\mapsto \boldsymbol{q}_t$ satisfy $\boldsymbol{q}_t=\Upsilon_t(\boldsymbol{p}_t)$ for all $t\in [0,T]$. Then, $\boldsymbol{p}$ is measurable with $\sup_{t\in [0,T]}\mathcal{F}(t,\boldsymbol{p}_t)<+\infty$ if and only if $\boldsymbol{q}$ is measurable with $\sup_{t\in [0,T]}\mathcal{E}(t,\boldsymbol{q}_t)<+\infty$.
\end{lemma}
\begin{proof}
We assume that  $\boldsymbol{p}$ is measurable and we prove that so is $\boldsymbol{q}$.The proof for the reverse implication works the same. 
Using \ref{eqn:coercivity-td} and the confinement condition, we see that $\boldsymbol{p}$  takes values in a metrizable subset $\mathcal{A}$ of $\mathcal{Q}^O$. Similarly, using the estimate \eqref{eqn:coercivity-ti}, which holds true also in the setting of Section \ref{sec:td}, we see that $\boldsymbol{q}$ takes also values in a metrizable subset $\mathcal{B}$ of $\mathcal{Q}^O$.
 Define $\Upsilon^{\mathcal{A},\mathcal{B}} \colon [0,T] \times \mathcal{A} \to  [0,T] \times \mathcal{B}$ by setting $\Upsilon^{\mathcal{A},\mathcal{B}}(t,\widehat{\boldsymbol{p}})\coloneqq (t,\Upsilon_t\restr{\mathcal{A}}(\widehat{\boldsymbol{p}}))$. Observe that $\Upsilon^{\mathcal{A},\mathcal{B}}$ is continuous as a consequence of \eqref{eqn:seq-cont} and the metrizability of $\mathcal{A}$ and $\mathcal{B}$. Therefore, writing $\boldsymbol{q}$ as the composition of the projection $(t,\widehat{\boldsymbol{q}})\mapsto \widehat{\boldsymbol{q}}$ from $[0,T] \times \mathcal{B}$ to $\mathcal{B}$, the function $\Upsilon^{\mathcal{A},\mathcal{B}}$, and the measurable map $t\mapsto (t,\boldsymbol{p}_t)$ from $[0,T]$ to $[0,T]\times \mathcal{A}$, we conclude that $\boldsymbol{q}$ is measurable.
\end{proof}

At this point, the proof Theorem~\ref{thm:existence-td} becomes immediate.

\begin{proof}[Proof of Theorem~\ref{thm:existence-td}]
Recall \eqref{eqn:upsilon}. Setting $\boldsymbol{p}^0\coloneqq \Upsilon_0(\boldsymbol{q}^0)\in \mathcal{Q}_{\boldsymbol{id}}^O$, From the identities \eqref{eq:id1}--\eqref{eq:id3} in the proof of Lemma \ref{lem:equiv} and from claims (ii) and (iv) of Proposition \ref{prop:equivalence-formulation}, we see that $\boldsymbol{p}^0$ verifies \eqref{eqn:ic-a}. Hence,	
by Theorem~\ref{thm:existence-a}, there exists an energetic solution $\boldsymbol{p}\colon t \mapsto \boldsymbol{p}_t$ to $(\mathcal{Q}_{\boldsymbol{id}}^O,\mathcal{F}\restr{\mathcal{Q}_{\boldsymbol{id}}^O},\mathcal{D}\restr{\mathcal{Q}_{\boldsymbol{id}}^O \times \mathcal{Q}_{\boldsymbol{id}}^O})$ which is measurable and satisfies $\boldsymbol{p}_0=\boldsymbol{p}^0$. Then, the function $\boldsymbol{q}\colon t\mapsto \boldsymbol{q}_t\coloneqq \Upsilon_t(\boldsymbol{p}_t)$ satisfies all the desired properties in view of Lemma \ref{lem:equiv} and Lemma \ref{lem:equiv-measurability}.
\end{proof}

\section*{Appendix}

\setcounter{section}{0}
\setcounter{subsection}{0}
\setcounter{theorem}{0}

\setcounter{equation}{0}

\renewcommand{\thesubsection}{A.\arabic{subsection}}

\renewcommand{\thetheorem}{A.\arabic{theorem}}

\renewcommand{\theequation}{A.\arabic{equation}}

This appendix collects a few subsidiary results that have been employed through the paper.

\subsection{An embedding inequality for Orlicz-Sobolev maps on spheres}
In this subsection, we register an embedding inequality  for Orlicz-Sobolev maps on spheres. This estimate is a simple reformulation of a result from \cite{CarozzaCianchi16}  which has been employed in the proof of Theorem~\ref{thm:reg-approx-diff}.  

Let $A$ be an N-function satisfying ($\Delta_2$) and \eqref{eqn:growth-infinity}. Define  $A_{N-1}\colon [0,+\infty) \to [0,+\infty)$ to be the conjugate function of $B_{N-1}\colon [0,+\infty) \to [0,+\infty)$  given by
\begin{equation}
	\label{eqn:A_N-1}
	B_{N-1}(s)\coloneqq \begin{cases}
		\bar{A}(s) & \text{if $N=2$,}\\
		s^{\frac{N-1}{N-2}} \displaystyle \int_s^{+\infty}  \frac{\bar{A}(\sigma)}{\sigma^{1+\frac{N-1}{N-2}}}\,\d\sigma & \text{if $N\geq 3$.}
	\end{cases} 
\end{equation}
By \cite[Lemma~3.2]{CarozzaCianchi16}, we have 
\begin{equation*}
	B_{N-1}(s)\geq \bar{A}(s) \quad \text{for all $s\geq 0$}
\end{equation*}
so that
\begin{equation*}
	\label{eqn:A_N-1-domination}
	A(s)\geq A_{N-1}(s) \quad \text{for all $s\geq 0$}
\end{equation*}
and
\begin{equation}
	\label{eq:A_N-1-inverse}
	A^{-1}(s)\leq A_{N-1}^{-1}(s) \quad \text{for all $s>0$.}
\end{equation}
In principle, $A_{N-1}$ is only a Young function (see, e.g.,  \cite[Definition 2.2]{MSZ}). From \cite[Proposition 4.3]{CarozzaCianchi19}, we see that  $A_{N-1}$ satisfies ($\Delta_2$). Then, $A_{N-1}$ satisfies also condition ($\nabla_2$) (see, e.g., \cite[Definition 2, p. 22]{rao.ren}) in view of \cite[Equations (2.8) and (4.17)]{CarozzaCianchi19}. Thus, \cite[Proposition 2.16]{MSZ} ensures that $A_{N-1}$ is actually an N-function.
In particular, $A_{N-1}$ is strictly increasing, so that
\begin{equation}
	\label{eq:aa}
	\text{$A_{N-1}^{-1}$ is strictly increasing and concave.}
\end{equation}
We now present the embedding inequality. Without seeking for generality nor optimality, we  derive an estimate which is tailored to our needs.

\begin{proposition}[Embedding inequality  for Orlicz-Sobolev maps on spheres]
	\label{thm:embedding-appendix}
Let $A$ be an N-function satisfying ($\Delta_2$) and \eqref{eqn:growth-infinity}. Then, there exists a constant $C_{\rm em}=C_{\rm em}(N)>0$ such that, given $\boldsymbol{x}_0\in\RN$ and $r\in(0,1)$, there holds
 \begin{equation*}
 	\esssup_{S(\boldsymbol{x}_0,r)} |\boldsymbol{v}|  \leq {C_{\rm em}} A_{N-1}^{-1} \left ( \dashint_{S(\boldsymbol{x}_0,r)} A(|D^{S(\boldsymbol{x}_0,r)}\boldsymbol{v}|)\,\d\haus +  \dashint_{S(\boldsymbol{x}_0,r)} A(|\boldsymbol{v}|)\,\d\haus\right )
 \end{equation*}
 for all $\boldsymbol{v}\in W^{1,A}(S(\boldsymbol{x}_0,r);\RN)$.
\end{proposition}
\begin{proof}
Let	$\boldsymbol{v}\in W^{1,A}(S(\boldsymbol{x}_0,r);\RN)$. Its oscillation is defined as
\begin{equation*}
	\osc_{S(\boldsymbol{x}_0,r)} \boldsymbol{v}\coloneqq \sum_{i=1}^N \left ( \esssup_{S(\boldsymbol{x}_0,r)} v^i - \essinf_{S(\boldsymbol{x}_0,r)} v^i\right) .
\end{equation*}
As a particular case of \cite[Theorem~4.1]{CarozzaCianchi16}, we have the estimate
\begin{equation}
	\label{eqn:car-cia}
	\osc_{S(\boldsymbol{x}_0,r)} \boldsymbol{v} \leq \widetilde{C}_{\rm em} A_{N-1}^{-1} \left(  \dashint_{S(\boldsymbol{x}_0,r)} A(|D^{S(\boldsymbol{x}_0,r)}\boldsymbol{v}|)\,\d\haus\right)
\end{equation}
for some constant $\widetilde{C}_{\rm em}=\widetilde{C}_{\rm em}(N)>1$.
Let $\boldsymbol{x},\boldsymbol{z}\in S(\boldsymbol{x}_0,r)$. Then
\begin{equation*}
	|\boldsymbol{v}(\boldsymbol{x})|\leq \osc_{S(\boldsymbol{x}_0,r)} \boldsymbol{v}+|\boldsymbol{v}(\boldsymbol{z})|.
\end{equation*}
By taking the integral average with respect to $\boldsymbol{z}\in S(\boldsymbol{x}_0,r)$, we obtain 
\begin{equation}
	\label{eq:cc1}
	|\boldsymbol{v}(\boldsymbol{x})|\leq \osc_{S(\boldsymbol{x}_0,r)} \boldsymbol{v}+\dashint_{S(\boldsymbol{x}_0,r)} |\boldsymbol{v}|\,\d\haus.
\end{equation}
Also
\begin{equation}
	\label{eq:cc2}
	\dashint_{S(\boldsymbol{x}_0,r)} |\boldsymbol{v}|\,\d\haus\leq A^{-1}\left( \dashint_{S(\boldsymbol{x}_0,r)} A(|\boldsymbol{v}|)\,\d\haus \right)\leq A^{-1}_{N-1}\left( \dashint_{S(\boldsymbol{x}_0,r)} A(|\boldsymbol{v}|)\,\d\haus \right)
\end{equation}
thanks to Jensen's inequality and \eqref{eq:A_N-1-inverse}. Therefore, combining \eqref{eqn:car-cia}--\eqref{eq:cc2}, we obtain
\begin{equation*}
	\begin{split}
		|\boldsymbol{v}(\boldsymbol{x})|&\leq \widetilde{C}_{\rm em} \left \{ A_{N-1}^{-1} \left(  \dashint_{S(\boldsymbol{x}_0,r)} A(|D^{S(\boldsymbol{x}_0,r)}\boldsymbol{v}|)\,\d\haus\right)+A^{-1}_{N-1}\left( \dashint_{S(\boldsymbol{x}_0,r)} A(|\boldsymbol{v}|)\,\d\haus \right) \right \}\\
		&\leq 2\widetilde{C}_{\rm em}  A_{N-1}^{-1} \left(  \dashint_{S(\boldsymbol{x}_0,r)} A(|D^{S(\boldsymbol{x}_0,r)}\boldsymbol{v}|)\,\d\haus+ \dashint_{S(\boldsymbol{x}_0,r)} A(|\boldsymbol{v}|)\,\d\haus \right), 
	\end{split}
\end{equation*}
where in the last line we exploited \eqref{eq:aa}.
\end{proof}

\subsection{Differential of tensor functions}
In this subsection, we recall the expressions for the differentials of basic tensor functions. We consider 
\begin{align*}
	\adj &\colon \RNN \to \RNN; \: \boldsymbol{A}\mapsto \adj \boldsymbol{A}, \qquad \cof \colon \RNN \to \RNN;\:\boldsymbol{A}\mapsto \cof \boldsymbol{A},\\
	  \mathrm{inv}&\colon \RNN_+ \to \RNN_+;\:\boldsymbol{A} \mapsto \boldsymbol{A}^{-1}, \qquad \hspace{3,5pt} \det \colon \RNN \to \R;\:\boldsymbol{A}\mapsto \det \boldsymbol{A}.
\end{align*}
For every $\boldsymbol{A}\in \RNN_+$, the differential of these maps can be computed using the formulas
\begin{align}
	\label{eqn:adj}
	\langle \d(\adj)(\boldsymbol{A}), \boldsymbol{B} \rangle &= ((\cof \boldsymbol{A}):\boldsymbol{B})\boldsymbol{A}^{-1}-(\adj \boldsymbol{A})\boldsymbol{B}\boldsymbol{A}^{-1},\\
	\label{eqn:cof}
	\langle \d(\cof)(\boldsymbol{A}), \boldsymbol{B} \rangle &= ((\cof \boldsymbol{A}):\boldsymbol{B})\boldsymbol{A}^{-\top}-(\cof \boldsymbol{A})\boldsymbol{B}^\top \boldsymbol{A}^{-\top},\\
	\label{eqn:inv}
	\langle \d (\mathrm{inv})(\boldsymbol{A}),\boldsymbol{B}\rangle &=-\boldsymbol{A}^{-1}\boldsymbol{B}\boldsymbol{A}^{-1},\\
	\label{eqn:det} 
	\langle \d(\det)(\boldsymbol{A}), \boldsymbol{B} \rangle &= \cof
	\boldsymbol{A}:\boldsymbol{B},
\end{align}
which are valid for all $\boldsymbol{B}\in\RNN$. Formula \eqref{eqn:det} is classical. The expression in \eqref{eqn:adj} can be deduced from \eqref{eqn:det} the identity $\boldsymbol{A}(\adj \boldsymbol{A})=(\det \boldsymbol{A})\boldsymbol{I}$. Then, \eqref{eqn:cof}--\eqref{eqn:inv} follows by the linearity of the transposition and the formula $\boldsymbol{A}^{-1}=(\det \boldsymbol{A})^{-1}(\adj \boldsymbol{A})$.

\subsection{Chain rule for the external field}
In this subsection, we derive simple chain rule for the composition of Sobolev space-valued and Lipschitz maps. This result has been employed in the proof of Proposition~\ref{prop:time-derivative-td}.

\begin{lemma}[Chain rule for the external field]
	\label{lem:k}
	Let $\boldsymbol{d}$, $\boldsymbol{d}^{-1}$, and $\boldsymbol{h}$ be as in \eqref{eqn:bd}--\eqref{eqn:bd-inverse} and \eqref{eqn:field}, respectively. Then, the function  
	$\boldsymbol{k}\colon t\mapsto \boldsymbol{k}_t\coloneqq \boldsymbol{h}_t \circ \boldsymbol{d}_t \det D \boldsymbol{d}_t$ belongs to the space $AC([0,T];L^1(O;\RN))$ and its time derivative is given by 
	\begin{equation}
		\label{eq:time-derivative-k}
		\dot{\boldsymbol{k}}_t=\left( (D\boldsymbol{h}_t\circ \boldsymbol{d}_t)\dot{\boldsymbol{d}}_t+\dot{\boldsymbol{h}}_t \circ \boldsymbol{d}_t\right) \det D \boldsymbol{d}_t+ \left( \cof D \boldsymbol{d}_t:D\dot{\boldsymbol{d}}_t \right)(\boldsymbol{h}_t \circ \boldsymbol{d}_t).
	\end{equation}
\end{lemma}
\begin{proof}
	First, define $\boldsymbol{l}\colon t\mapsto \boldsymbol{l}_t\coloneqq \boldsymbol{h}_t \circ \boldsymbol{d}_t$.
	As both $\boldsymbol{d}_t$ and $\boldsymbol{d}_t^{-1}$ are Lipschitz, $\boldsymbol{l}_t\in W^{1,1}(O;\RN)$ for all $t\in [0,T]$ by \cite[Theorem~11.53]{Leo}. 
	Define ${\boldsymbol{\delta}}\colon (0,T)\times O \to \RN$ and ${\boldsymbol{\eta}}\colon (0,T)\times O \to \RN$ by setting ${\boldsymbol{\delta}}(t,\boldsymbol{w})\coloneqq \boldsymbol{d}_t(\boldsymbol{w})$ and ${\boldsymbol{\eta}}(t,\boldsymbol{w})\coloneqq \boldsymbol{h}_t(\boldsymbol{w})$. Given the regularity of $\boldsymbol{d}$ and $\boldsymbol{h}$, using Fubini's theorem together with standard properties of the Bochner integral, we see that $\boldsymbol{\delta}\in W^{1,\infty}((0,T)\times O;O)$ and ${\boldsymbol{\eta}}\in W^{1,1}((0,T)\times O;\RN)$.
	
	Define $\boldsymbol{P}\colon (0,T) \times O \to (0,T)\times O$ as $\boldsymbol{P}(t,\boldsymbol{w})\coloneqq (t,{\boldsymbol{\delta}}(t,\boldsymbol{w}))$. Then, $\boldsymbol{P}$ is injective, Lipschitz and so is its inverse. Thus, setting ${\boldsymbol{\lambda}}\coloneqq {\boldsymbol{\eta}}\circ \boldsymbol{P}$, we have ${\boldsymbol{\lambda}}\in W^{1,1}((0,T)\times O;\RN)$ and the chain rule gives 
	\begin{equation*}
		\partial_t {\boldsymbol{\lambda}}(t,\boldsymbol{w})=({\boldsymbol{\eta}} (t,{\boldsymbol{\delta}}(t,\boldsymbol{w})))\partial_t {\boldsymbol{\delta}}(t,\boldsymbol{w})+\partial_t {\boldsymbol{\eta}}(t,{\boldsymbol{\delta}}(t,\boldsymbol{w}))
	\end{equation*}
	again by \cite[Theorem~11.53]{Leo}. 
	From this, observing that $\boldsymbol{\lambda}(t,\boldsymbol{w})=\boldsymbol{l}_t(\boldsymbol{w})$ and using the same argument based on Fubini's theorem as above, we deduce that $\boldsymbol{l}\in AC([0,T];L^1(\Omega;\RN))$ with $\dot{\boldsymbol{l}}_t=(D\boldsymbol{h}_t \circ \boldsymbol{d}_t)\dot{\boldsymbol{d}}_t+\dot{\boldsymbol{h}}_t\circ \boldsymbol{d}_t$. At this point, $\boldsymbol{k}\in AC([0,T];L^1(\Omega;\RN))$ follows and the formula for its time derivative is obtained by using \eqref{eqn:det}. 
\end{proof}

\section*{Acknowledgements}
The First Author acknowledges the support of the Alexander von Humboldt Foundation  through the \mbox{Humboldt} Research Fellowship for Postdocs. The second Author was supported by PRIN Project 2022E9CF89 “Geometric Evolution Problems and Shape Optimizations”. PRIN projects are part of PNRR Italia Domani, financed by the European Union through NextGenerationEU. Both Authors are members  of Gruppo Nazionale per l’Analisi Matematica, la Probabilit\`a e le loro Applicazioni (GNAMPA) of Instituto Nazionale di Alta Matematica (\mbox{INdAM}). Part of this work has been developed during the summer school “Variational and PDE Methods in Nonlinear Science” in Cetraro (Italy) hosted by Fondazione Centro Internazionale Matematico Estico (CIME), and the “XXIV Symposium on Trends in Applications of Mathematics to Mechanics” (STAMM) in W\"{u}rzburg (Germany)  organized by the International Society for the Interaction of Mechanics and Mathematics (ISIMM). Both Institutions are gratefully acknowledged by the Authors.


\begin{thebibliography} {99}

\bibitem{Adams} 
{R. A. Adams}, \emph {Sobolev Spaces}, Academic Press, New York (1975).

\bibitem{alberico.cianchi}
{A. Alberico, A. Cianchi}, \emph{Differentiabiliy properties of Orlicz-Sobolev functions}, Ark. Mat. {43}, 1--28 (2005).

\bibitem{ball.convexity}
{J. M. Ball}, \emph{Convexity conditions and existence theorems in nonlinear elasticity}, Arch. Ration. Mech. Anal. 63, 337–403 (1976/1977).

\bibitem{ball.lc}
J. M. Ball, \emph{Mathematics and liquid crystals},  Mol. Cryst. Liq. Cryst. 674, 1--27 (2017).

\bibitem{Ball.op} 
{J.~M.~Ball}, \emph{Some open problems in elasticity}. In {P. Newton, P. Holmes,  A. Weinstein}, \emph{Geometry,
Mechanics, and Dynamics}, Springer, New York, 3--59  (2002).

\bibitem{ball.currie.olver}
{J.~M.~Ball, J.~C.~Currie, P.~J.~Olver},  \emph{Null Lagrangians, Weak Continuity, and Variational Problems of Arbitrary Order}, J. Funct. Anal. {41}, 135--174  (1981).


\bibitem{BDesimone} 
{M. Barchiesi, A. DeSimone}, \emph{Frank energy for nematic elastomers: a nonlinear model}, {ESAIM Control Optim. Calc. Var.} {21}, 277--372  (2015).

\bibitem{BHM17} 
{M. Barchiesi, D. Henao, C. Mora-Corral}, \emph{Local Invertibility in Sobolev Spaces with Applications to Nematic Elastomers and Magnetoelasticity}, {Arch. Rational Mech. Anal.} { 224}, 743--816 (2017).


\bibitem{bresciani} {M. Bresciani}, \emph{Quasistatic evolution in magnetoelasticity under subcritical coercivity assumptions}, Calc. Var. {62}, no. 7, Article no. 181 (2023).

\bibitem{bresciani.davoli.kruzik}
{M. Bresciani, E. Davoli, M. Kru\v{z}\'{i}k}, \emph{Existence results in large-strain magnetoelasticity}, Ann. Inst. H. Poincar\'{e} Analyse Non Lin\'{e}aire 40, 557--592 (2023).

\bibitem{bresciani.friedrich.moracorral}
{M. Bresciani, M. Friedrich, C. Mora-Corral}, \emph{Variational models with Eulerian-Lagrangian formulation allowing for material failure}, arXiv preprint (2024), available at \url{https://arxiv.org/abs/2402.12870}.



\bibitem{CarozzaCianchi19} {M. Carozza,  A. Cianchi}, \emph{Continuity properties of weakly monotone Orlicz-Sobolev functions}, {Adv. Calc.  Var.} {14}, no. 1, 107--126 (2021).

\bibitem{CarozzaCianchi16} {M. Carozza,  A. Cianchi}, \emph{Smooth approximation of Orlicz-Sobolev maps between manifolds}, 
{Potential Anal.} {45}, 557--578 (2016).


\bibitem{Cianchi96} {A. Cianchi}, \emph{Continuity properties of functions from {O}rlicz-{S}obolev spaces and embedding theorems},
{Ann. Scuola Norm. Sup. Pisa Cl. Sci.}  23, no. 4, 575--608 (1996).


\bibitem{cianchi.embedding} 
A. Cianchi, \emph{A sharp embedding theorem for Orlicz-Sobolev spaces}, Indiana Univ. Math. J.. 45, no. 1,  39--65 (1996).

\bibitem{cianchi.trace}
A. Cianchi, \emph{Orlicz–Sobolev boundary trace embeddings}, Math. Z. 266, 431--449 (2010).







\bibitem{Daco} { B. Dacorogna}, \emph{Direct methods in the calculus of variations}, Applied Mathematical Sciences Vol. 78, Second Edition, Springer, New York (2008).



\bibitem{dalmaso.lazzaroni}
{ G. Dal Maso, G. Lazzaroni}, \emph{Quasistatic crack growth in finite elasticity with non-interpenetration}, Ann. Inst. H. Poincar\'{e} Anal. Non Lin\'{e}aire {27}, 257--290 (2010).




\bibitem{FoLe} {I. Fonseca, G. Leoni}, \emph{Modern Methods in the Calculus of Variations: $L^p$ spaces}, Springer Monographs in Mathematics, Springer, New York (2007).

\bibitem{FoGa95book}{I.~Fonseca, W.~Gangbo}, {\em {Degree theory in analysis and
  applications}}, Oxford University Press, New York (1995).


\bibitem{FraMie} {G. Francfort, A. Mielke}, \emph{Existence results for a class of rate-independent material models with nonconvex elastic energies}, J. Reine Angew. Math. {595}, 55--91 (2006).

\bibitem{giaquinta.modica.soucek}
{M. Giaquinta, G. Modica, J. Sou\v{c}ek}, \emph{Cartesian Currents in the Calculus of Variations I. Cartesian Currents}, Ergebnisse der Mathematik und ihrer Grenzgebiete, 3. Folge / A Series of Modern Surveys in Mathematics Vol. 37, Springer-Verlag, Berlin - Heidelberg (1998).

\bibitem{goffman.ziemer}
{C. Goffman, W. P. Ziemer}, \emph{Higher dimensional mappings for which the area formula holds}, {Ann. Math.} {92}, no. 2, 482--488 (1970).

%


%
\bibitem{hh}
{P. Harjuletho, P. H\"{a}st\"{o}}, \emph{Orlicz spaces and generalized Orlicz spaces}, Springer Nature Switzerland AG (2019).
%




\bibitem{HeMo11} {D. Henao, C. Mora-Corral}, \emph{Fracture surface and regularity of inverses for $BV$ deformations}, {Arch. Rat. Mech. Anal.} {201}, 575--629 (2011).

%

\bibitem{HS}{ D. Henao, B. Stroffolini}, \emph{On Sobolev-Orlicz nematic elastomers}, {Nonlinear Anal.} {194}, 111513 (2020). 



\bibitem{KKM} J. Kauhanen, P. Koskela, J. Mal\'{y}, \emph{On functions with derivatives in a Lorentz space}, Manuscripta Math. {100}, 87--101  (1999).

\bibitem{Kras}  {M. A. Krasnosel'ski\v{i}, Ya. B. Ruticki\v{i}}, \emph{Convex Functions and Orlicz Spaces}, P. Noordhoff LTD., Groningen, The Netherlands  (1961).

\bibitem{kruzik.roubicek}
M. Kru\v{z}\'ik, T. Roub\'{i}\v{c}ek, \emph{Mathematical methods in continuum mechanics of solids}, Interaction of mechanics and mathematics, Springer Nature Switzerland AG (2019).

\bibitem{kruzik.stefanelli.zeman} M. Kru\v{z}\'ik, U. Stefanelli, J. Zeman, \emph{Existence results for incompressible magnetoelasticity}, {Discrete Cont. Dyn. Syst. A} { 35} , 2615--2623 (2015).

\bibitem{Kufn} A. Kufner, O. John and S. Fucik, \emph{Function Spaces}, Springer Netherlands (1977).

\bibitem{Laz} {G. Lazzaroni}, \emph{Quasistatic crack growth in finite elasticity with Lipscthiz data}, Annali di Matematica {190}, 165--194  (2011).

\bibitem{Leo}
{G. Leoni}, \emph{A first course on Sobolev spaces}, Second Edition, Graduate Studies in Mathematics Vol. 181, American Mathematical Society (2017), Providence (Rodhe Island).


\bibitem{Maly.ac}
J. Mal\'{y}, \emph{Absolutely continuous functions of several variables}, J. Math. Anal. Appl. 231, 492–508  (1999).

\bibitem{Maly} J. Mal\'{y},  \emph{Coarea property of Sobolev maps}. In B. Opic and J. R\'{a}kosn\'{i}k,  \emph{Nonlinear Analysis, Function Spaces
and Applications Vol. 7. Proceedings of the Spring School Held in Prague, July 17--22, 2002,}, Math. Inst. Acad. Sci. of the Czech Republic,
Praha,  142--192 (2003).

\bibitem{MSZ} {J. Maly, D. Swanson, W. P. Ziemer}, \emph{Fine behavior of functions whose gradients are in an Orlicz space}, {Studia Math.} { 190}, no. 1, 33--70 (2009).


\bibitem{MieRu}
{A. Mielke, T. Roub\'{i}\v{c}ek}, \emph{Rate-Independent Systems. Theory and Application}, Springer, New York (2015).

\bibitem{MoCo} {C. Mora-Corral}, \emph{Quasistatic evolution of cavities in nonlinear elasticity}, SIAM J. Math. Anal. {46}, no. 1, 523--571 (2014).

\bibitem{moracorral.oliva}
{C. Mora-Corral, M. Oliva}, \emph{Relaxation of nonlinear elastic energies involving the deformed configuration and applications to nematic elastomers}, ESAIM  Control Optim. Calc. Var. 25, Article no. 19 (2019).



\bibitem{mueller}
{S.~M\"{u}ller}, \emph{ Weak continuity of determinants and nonlinear elasticity}, C. R. Acad. Sci. Paris Sér. I Math.
{307}, no. 9, 501--506 (1988).

\bibitem{mueller.tang.yan}
S. M\"{u}ller, Tang Qi, B. S. Yan, \emph{On a new class of elastic deformations not allowing for cavitation}, Ann. Inst. Henri Poincar\'{e} Analyse Non Lin\'{e}aire  ll, no. 2,  217-243 (1994).

  
\bibitem{O}
  J. Onninen, \emph{Differentiability of monotone Sobolev functions}, Real Anal. Exchange 26, 761–772. 



%




\bibitem{rao.ren}
{M. M. Rao, Z. D. Ren}, \emph{Theory of Orlicz Spaces}, Monographs Textbooks Pure
and Applied Mathematics Vol. 146, Marcel Dekker Inc., New York - Basel - Hong Kong (1991).

\bibitem{scilla.stroffolini}
{G. Scilla, B. Stroffolini}, \emph{Relaxation of nonlinear elastic energies related to Orlicz–Sobolev nematic elastomers}, Rend. Lincei Mat. Appl. 31, 349--389 (2020).

\bibitem{scilla.stroffolini.inv}
{G. Scilla, B. Stroffolini}, \emph{Invertibility of Orlicz-Sobolev maps}. In { M. I. Español, M. Lewicka, L. Scardia, A. Schlömerkemper}, \emph{Research in Mathematics of Materials Science},  Association for Women in Mathematics Series Vol. 31, Springer, Cham, 297--317 (2022).



\end{thebibliography}
\end{document}